\documentclass[11pt]{amsart}
\usepackage{amssymb}

\newcommand\datver[1]{\def\datverp%
 {\par\boxed{\boxed{\text{#1; Run: \today}}}}}

\newcommand\Kf{\kappa}

\newcommand\fib{\pi_0}
\newcommand\tfib{\widetilde{\fib}}
\newcommand\Vfz{\Vf_0}

\newcommand\Lin{A}

\newcommand{\hilbert}{\mathsf{H}}

\newcommand{\rcal}{{\mathcal{R}}}

\newcommand\xib{{\underline\xi}}
\newcommand\etab{{\underline\eta}}
\newcommand\zetab{{\underline\zeta}}
\newcommand\taub{{\underline\tau}}
\newcommand\mub{{\underline\mu}}

\newcommand\xibs{{\underline{\hat\xi}}}
\newcommand\etabs{{\underline{\hat\eta}}}
\newcommand\zetabs{{\underline{\hat\zeta}}}
\newcommand\taubs{{\underline{\hat\tau}}}

\newcommand\xis{{\hat\xi}}
\newcommand\etas{{\hat\eta}}
\newcommand\zetas{{\hat\zeta}}
\newcommand\taus{{\hat\tau}}

\newcommand{\Tmstar}{{}^{\mo}T^*}
\newcommand{\Smstar}{{}^{\mo}S^*}

\newcommand{\Te}{{}^{{\eo}}T}
\newcommand{\Testar}{{}^{\eo}T^*}

\newcommand{\Sestar}[1][\mbox{}]{{}^{\eo}S^*_{#1}}
\newcommand{\Tbstar}{{}^{\bo} T^*}
\newcommand{\Tbdotstar}{{}^{\bo}\overset{\!\! .}{T^*}}
\newcommand{\Sbdotstar}{{}^{\bo}\overset{\!\! .}{S^*}}
\newcommand{\Sbstar}{{}^{\bo} S^*}
\newcommand{\abs}[1]{{\left\lvert{#1}\right\rvert}}
\newcommand{\norm}[1]{{\left\lVert{#1}\right\rVert}}
\renewcommand{\Box}{{\square}}
\newcommand{\tBox}{\tilde{\square}}
\newcommand\h{\frac{1}{2}}

\newcommand{\bH}[1]{{H^{#1}_{\bo}}}

\newcommand{\eH}[1]{{H^{#1}_{\eo}}}
\newcommand{\eHloc}[1]{{H^{#1}_{{\eo},\text{loc}}}}

\DeclareMathOperator{\Op}{\text{Op}}
\newcommand{\trv}{{\operatorname{trv}}}
\DeclareMathOperator{\Card}{\text{Card}}

\DeclareMathOperator{\eWF}{{{WF}_{\eo}}}

\newcommand{\dcal}{{\mathcal{D}}}

\newcommand{\hcal}{{\mathcal{H}}}
\newcommand{\ecal}{{\mathcal{E}}}

\newcommand{\fcal}{{\mathcal{F}}}
\newcommand{\fcalbar}{{\overline{\mathcal{F}}}}
\newcommand{\fcaldot}{{\dot{\mathcal{F}}}}

\newcommand{\ang}[1]{{\left\langle{#1}\right\rangle}}

\newcommand{\ep}{{\epsilon}}
\newcommand{\tensor}{{\otimes}}

\newcommand{\gcal}{{\mathcal{G}}}

\newcommand{\ePs}[1]{{\Psi^{#1}_{\eo}}}

\newcommand{\VPs}[1]{{\Psi^{#1}_{\Vf}}}

\newcommand{\dom}{\dcal}
\newcommand{\domt}{{\tilde\dcal}}

\numberwithin{equation}{section}
\newtheorem{theorem}{Theorem}[section]
\newtheorem{lemma}[theorem]{Lemma}
\newtheorem{proposition}[theorem]{Proposition}
\newtheorem{corollary}[theorem]{Corollary}

\newtheorem{non-theorem}[theorem]{Non-Theorem}

\theoremstyle{remark}
\newtheorem{definition}[theorem]{Definition}
\newtheorem{remark}[theorem]{Remark}

\newcommand{\esigma}[1][\mbox{}]{{{}^{\eo}\sigma_{#1}}}
\newcommand{\Vsigma}[1][\mbox{}]{{{}^{\Vf}\sigma_{#1}}}

\newcommand{\bpi}{\pi_{\bo}}

\newcommand\bo{\operatorname{b}}
\newcommand\bco{\operatorname{b\infty}}

\newcommand\eo{\operatorname{e}}
\newcommand\mo{x\operatorname{e}}

\newcommand\Vb{\mathcal{V}_{\bo}}
\newcommand\Ve{\mathcal{V}_{\eo}}

\newcommand\Vf{\mathcal{V}}

\newcommand\VT{{}^{\Vf}T}
\newcommand\VTstar{{}^{\Vf}T^*}
\newcommand\VSstar{{}^{\Vf}S^*}

\newcommand\bT{{}^{\bo}T}

\newcommand\eT{{}^{\eo}T}

\newcommand\cFTs{{}^{\Phi}\overline{T}\kern-1pt{}^*}

\newcommand{\kbar}{{\bar{K}}}

\newcommand\WF{\operatorname{WF}}

\newcommand\dWF{\WF_{\bo,\domt}}
\newcommand\mWF{\WF_{\bo,\domt'}}
\newcommand\VWF{\WF_{\Vf}}
\newcommand\bWF{\WF_{\bo}}

\newcommand{\tx}{{\tilde{x}}}
\newcommand{\ty}{{\tilde{y}}}
\newcommand{\tz}{{\tilde{z}}}

\newcommand\cA{\mathcal{A}}
\newcommand\cB{\mathcal{B}}
\newcommand\cC{\mathcal{C}}

\newcommand\cL{\mathcal{L}}
\newcommand\cO{\mathcal{O}}

\newcommand\hyp{\mathcal{H}}
\newcommand\cG{\mathcal{G}}
\newcommand\cE{\mathcal{E}}

\newcommand\dhyp{\dot{\mathcal{H}}}
\newcommand\dcG{\dot{\mathcal{G}}}
\newcommand\dcE{\dot{\mathcal{E}}}

\newcommand\NN{\mathbb N}
\newcommand\bbN{\mathbb N}

\newcommand\RR{\mathbb R}

\newcommand\CIc{{\mathcal{C}}^{\infty}_c}

\newcommand\CI{{\mathcal{C}}^{\infty}}

\newcommand\CmI{{\mathcal{C}}^{-\infty}}

\newcommand\Diffb[1]{\operatorname{Diff}^{#1}_{\bo}}

\newcommand\Diffe[1]{\operatorname{Diff}^{#1}_{\eo}}
\newcommand\DiffV[1]{\operatorname{Diff}^{#1}_{\Vf}}

\newcommand\Psib{\Psi_{\bo}}
\newcommand\Psibc{\Psi_{\bco}}

\newcommand\cFNs{{}^{\Phi}\overline N\kern-1pt{}^*}

\newcommand\Id{\operatorname{Id}}
\newcommand{\lag}{\mathcal{L}}
\newcommand\Lap{\varDelta}

\newcommand\Span{\operatorname{sp}}

\newcommand\Sym{\operatorname{Sym}}

\newcommand\codim{\operatorname{codim}}

\newcommand\dCI{\dot{\mathcal{C}}^{\infty}}

\newcommand\loc{{\text{loc}}}

\newcommand\pa{\partial}

\newcommand\restrictedto{\!\!\upharpoonright}

\newcommand\sgn{\operatorname{sgn}}

\newcommand\supp{\operatorname{supp}}

\renewcommand\Re{\operatorname{Re}}

\newcommand{\tbar}{\overline{t}}
\newcommand{\ybar}{\overline{y}}
\newcommand{\xibar}{\overline{\xi}}
\newcommand{\etabar}{\overline{\eta}}

\newcommand{\taubar}{\overline{\tau}}
\newcommand{\zbar}{\overline{z}}

\newcommand\Mand{\text{ and }}

\newcommand\Mhaving{\text{ having }}

\newcommand\Mif{\text{ if }}

\newcommand\Min{\text{ in }}

\newcommand\Mof{\text{ of }}
\newcommand\Mon{\text{ on }}

\newcommand\Msatisfies{\text{ satisfies }}

\newcommand\Msatisfying{\text{ satisfying }}

\newcommand\Mwith{\text{ with }}

\newcommand{\F}{{\mathcal{F}}}
\newcommand{\coiso}{{\mathcal{A}}}
\newcommand{\module}{{\mathcal{M}}}
\author{Richard Melrose}
\author{Andr\'as Vasy}
\author{Jared Wunsch}

\thanks{The authors gratefully acknowledge financial support for this
project, the first from the National Science Foundation under grant
DMS-0408993, the second under grant DMS-0201092, from a Clay Research
Fellowship and a Sloan Fellowship and the third from the National Science
Foundation under grants DMS-0401323 and DMS-0700318.  They thank two
anonymous referees for many suggestions that improved the exposition.}

\title[Propagation on edge manifolds]%
{Propagation of singularities for the wave equation on edge manifolds}
\datver{1.0A; Revised: 22-12-2006}
\date{\today}
\begin{document}

\begin{abstract}
We investigate the geometric propagation and diffraction of singularities
of solutions to the wave equation on manifolds with edge singularities.
This class of manifolds includes, and is modelled on, the product of a smooth
manifold and a cone over a compact fiber. Our main results are a general
`diffractive' theorem showing that the spreading of singularities at the
edge only occurs along the fibers and a more refined `geometric' theorem
showing that for appropriately regular (nonfocusing) solutions, the main
singularities can only propagate along geometrically determined rays.
Thus, for the fundamental solution with initial pole sufficiently close to
the edge, we are able to show that the regularity of the diffracted front
is greater than that of the incident wave.
\end{abstract}

\maketitle

\tableofcontents

\section{Introduction}

\subsection{Main Results}\label{intro:results}

In this paper, we investigate the geometric propagation and diffraction of
singularities of solutions to the wave equation on manifolds with edge
singularities. The main results are extension of those in
\cite{Melrose-Wunsch1} for the particular case of conic metrics, namely a
general `diffractive' theorem limiting the possible spreading of
singularities at the boundary and a more refined `geometric' theorem
showing that for appropriately regular (nonfocusing) solutions the main
singularities can only propagate along geometrically determined rays.

Let $X$ be an $n$-dimensional manifold with boundary, where the boundary,
$\pa X,$ is endowed with a fibration
$$
Z \to \pa X \overset{\fib}{\to} Y,
$$ where $Y,Z$ are without boundary.  Let $b$ and $f$ respectively denote
the dimensions of $Y$ and $Z$ (the `base' and the `fiber'). By an
\emph{edge metric} $g$ on $X$ we shall mean a metric, $g,$ on the interior
of $X$ which is a smooth 2-cotensor up to the boundary but which
degenerates there in a way compatible with the fibration. To state this
condition precisely, consider the Lie algebra $\Vf(X)$ of smooth vector
fields on $X$ which are tangent to the boundary, hence have well-defined
restrictions to the boundary which are required to be tangent to the fibers
of $\fib.$ Let $x$ be a boundary defining function for $X,$ then we require
that $g$ be degenerate in the sense that
\begin{equation}
g(V,V)\in x^2\CI(X)\ \forall\ V\in\Vf(X).
\label{24.8.2006.2}\end{equation}
This of course only fixes an `upper bound' on $g$ near the boundary and we
require the corresponding lower bound and also a special form of the
leading part of the metric near the boundary. The lower bound is just the
requirement that for any boundary point,
\begin{multline}
V\in\Vf(X)\Mand x^{-2}g(V,V)(p)=0\\
\Longrightarrow V=\sum\limits_{i}f_iV_i,\
V_i\in\Vf(X),\ f_i\in\CI(X),\ f_i(p)=0;
\label{24.8.2006.6}\end{multline}
that is the vanishing of $g(V,V)$ to higher order than $2$ at a boundary
point means that $V$ vanishes at that point as an element of $\Vf(X).$
To capture the required special form near $\pa X,$ consider within $\Vf(X)$
the subalgebra of vector fields which vanish at the boundary in the
ordinary sense, we denote this Lie subalgebra
\begin{equation}
\Vfz(X)\subset\Vf(X).
\label{24.8.2006.3}\end{equation}
In addition to \eqref{24.8.2006.2} we require that there be a product
decomposition of a neighborhood of the boundary 
\begin{equation}
[0,\epsilon)\times\pa X\overset{p}{\underset{\sim}\longrightarrow} U\supset\pa X
\label{24.8.2006.4}\end{equation}
and a smooth metric $g_0$ (in the usual sense of extension across the boundary)
on $[0,\epsilon )\times Y$ such that
\begin{equation}
\begin{gathered}
(g-\tfib^*g_0)(V,W)\in x^2\CI(X),\ \forall\ V,\ W\in\Vf(X)\Mand\\
(g-\tfib^*g_0)(V,W)\in x^3\CI(X),\ \forall\ V\in\Vfz(X),\ W\in\Vf(X);
\end{gathered}
\label{24.8.2006.5}\end{equation}
here $\tfib:U\longrightarrow [0,\epsilon)\times Y$ is the product extension
of $\fib$ using \eqref{24.8.2006.4}.

We note that the metric $g_0$ on $[0,\epsilon )\times Y$ in \eqref{24.8.2006.5}
can always be brought to the product form
\begin{equation}\label{edgemetric1}
g_0= dx^2+ h(x)
\end{equation}
near $x=0$ for some boundary defining function $x,$ with $h \in
\CI([0,\ep') \times Y; \Sym^2 T^* ([0,\ep) \times Y)),$ $\ep'>0,$ i.e.\
$h(x)$ an $x$-dependent metric on $Y.$ Thus the global restriction
\eqref{24.8.2006.5} means that in local coordinates $x,$ $y,$ $z$ near each
boundary point,
\begin{equation}\label{edgemetric2}
g=dx^2+h(x,y,dy)+xh'(x,y,z,dx,dy)+x^2k(x,y,z,dx,dy,dz)
\end{equation}
where the $h'\in\CI(U;\Sym^2 T^*([0,\epsilon)\times Y)$ and $k\in\CI(U;
\Sym^2 T^*X)$ and the restriction of $k$ to each fiber of the boundary is
positive-definite. This latter condition is equivalent to \eqref{24.8.2006.6}.

A manifold with boundary equipped with such an edge metric will also be
called an \emph{edge manifold} or a manifold with \emph{edge structure}. We
draw the reader's particular attention to the two extreme cases: if $Z$ is
a point, then an edge metric on $X$ is a simply a metric in the usual
sense, smooth up to the boundary, while if $Y$ is a point, $X$ is a conic
manifold (cf\@. \cite{Melrose-Wunsch1}).  A simple example of a more
general edge metric is obtained by performing a real blowup on a
submanifold $B$ of a smooth, boundaryless manifold $A.$ The blowup
operation simply introduces polar `coordinates' near $B$, i.e.\ replaces
$B$ by its spherical normal bundle, thus yielding a manifold $X$ with
boundary.  The pullback of a smooth metric on $A$ to $X$ is then an edge
metric; here we have $Y=B$ and $F=S^{\codim(B)-1}.$ A \emph{non-example}
that is nonetheless quite helpful in visualizing many of the constructions
used in this paper is a manifold with a codimension-two corner, equipped
with an incomplete metric: if we blow up the corner, we obtain a manifold
of the type considered in this paper, \emph{but where $F$ is now a manifold
with boundary,} given by a closed interval representing an angular variable
at the corner.  More generally, any manifold with corners with an
incomplete metric can be considered to have an iterated edge-structure:
passing to polar coordinates near a corner yields an edge manifold with
fiber given by another manifold with corners, more precisely the
intersection of a sphere with an orthant.  The authors intend to consider
this situation in a subsequent paper.  Note that in every case the boundary
of $X$ is given \emph{geometrically} by $Y,$ as the fibers become
metrically small as we approach the boundary.

Since we work principally with the wave equation, set $M=\RR\times X.$
Thus, the boundary of $M$ has a fibration with fiber $Z$ and base
$\RR\times Y;$ it is an edge manifold with a metric of Lorentzian signature.
We consider solutions $u$ to the wave equation 
\begin{equation}
\Box u = (D_t^2-\Lap_g)u = 0 \Mon M
\label{30.8.2006.8}\end{equation}
with respect to this Lorentzian metric. In the simplest case, in which $X$
decomposes as the product $\RR_+ \times Y \times Z$ and the metric is a warped
product,
$$
\Box = D_t^2-\left(D_x^2+\frac cx D_x +\Lap_Y +\frac{1}{x^2}\Lap_Z\right).
$$
In general the form of the operator is a little more complicated than this,
but with similar leading part at the boundary.

We shall consider below only solutions of \eqref{30.8.2006.8} lying in some
`finite energy space.' Thus, if $\dcal_\alpha$ is 
the domain of 
$\Lap^{\alpha/2},$ where $\Lap$ is the Friedrichs extension from the space
$\dot{\mathcal{C}}^\infty(X),$ of smooth functions vanishing to infinite
order at the boundary, we require that a solution be \emph{admissible} in
the sense that it lies in $\mathcal{C}(\RR; \dcal_\alpha)$ for \emph{some}
$\alpha\in\RR.$

In terms of adapted coordinates $x,$ $y,$ $z$ near a boundary point, an
element of $\Vf(X)$ is locally an arbitrary smooth combination of the basis
vector fields 
\begin{equation}
x\pa_x,\ x\pa_{y_j},\ \pa_{z_k}
\label{30.8.2006.9}\end{equation}
and hence $\Vf(X)$ is equal to the space of all sections of a vector bundle
(determined by $\fib$), which we call the \emph{edge tangent bundle} and
denote $\Te X.$ This bundle is canonically isomorphic to the usual tangent
bundle over the interior (and non-canonically isomorphic to it globally)
with a well-defined bundle map $\Te X\longrightarrow TX$ which has rank $f$
over the boundary. Correspondingly, on $M,$ the symbol of the wave operator
extends to be of the form $x^{-2}p$ with $p$ a smooth function on the dual
bundle to $\Te X,$ called the \emph{edge cotangent bundle} and denoted
$\Testar M.$ The vector fields in \eqref{30.8.2006.9}, and $x\pa_t$, define
linear functions $\tau,$ $\xi,$ $\eta$ and $\zeta$ on this bundle, in terms
of which the canonical one-form on the cotangent bundle lifts under the
dual map $\Testar M\longrightarrow T^*M$ to the smooth section
$$
\tau \frac{dt}x + \xi \frac{dx}x + \eta \cdot \frac{dy}x + \zeta \cdot dz
\Mof \Testar M.
$$ We will let $\Sestar M$ denote the unit cosphere bundle of $\Testar M.$

The characteristic variety of $\Box,$ i.e.\ the zero set of the symbol, will be
interpreted as a subset $\Sigma\subset\Testar M\setminus0;$ over the
boundary it is given by the vanishing of $b.$ Note that one effect of
working on this `compressed' cotangent bundle is that $\Sigma$ is smooth
and is given in each fiber by the vanishing of a non-degenerate Lorentzian
quadratic form.

The variables $\xi$ and $\zeta$ dual to $x$ and $z$ respectively play a
different role to the duals of the base variables $(t,y).$ We will thus
define a new bundle, denoted $\pi(\Sestar(M)),$ whose sections are
required to have $\xi=\zeta=0$ at $\pa M.$ (For an explanation of the
notation, see \S\ref{section:edgeb}).  Over $\pa M,$ we may then decompose
this bundle into elliptic, glancing, and hyperbolic sets much as in the
usual case of manifolds with boundary:
$$
\pi(\Sestar M) = \ecal \cup \gcal \cup \hcal
$$
where $\gcal$ is the light cone, and $\ecal$ its exterior, i.e.\ $\gcal \cup
\hcal$ is the projection of the characteristic set to $\pi (\Sestar M).$

For each normalized point
\begin{equation*}
p=(\tbar,\ybar,\zbar,\taubar=\pm 1,\etabar)\in\hcal,\ |\etabar|<1,
\label{30.8.2006.11}\end{equation*}
it is shown below that there are two line segments of `normal' null
bicharacteristics in $\Sigma,$ each ending at one of the two points above $p$
given by the solutions $\xibar$ of $\xibar^2+\abs{\etabar}^2 =1.$ These
will be denoted
$$
\fcal_{\bullet,p},
$$
where $\bullet$ is permitted to be $I$ or $O,$ for `incoming' or
`outgoing,' as $\sgn\xi=\pm \sgn \taubar$ ($+$ for $I$
and $-$ for $O).$  In the special case that the fibration and metric are of
true product form
$$
dx^2+ h(y,dy) + x^2 k(z,dz),
$$
these bicharacteristics are simply
\begin{gather*}
\fcal_{I,p} = \{t\leq \tbar, x=(\tbar-t)\abs\xibar,
y=y(t), z=\zbar, \tau=\taubar, \xi=\xibar,\eta=\eta(t),\zeta=0\}\\
\Mand\\
\fcal_{O,p}= \{t\geq \tbar, x=(t-\tbar)\abs\xibar,
y=y(t), z=\zbar, \tau=\taubar, \xi=\xibar,\eta=\eta(t),\zeta=0\};
\end{gather*}
where $(y(t), \eta(t))$ evolves along a geodesic in $Y$ which passes
through $(\ybar,\etabar)$ at time $t=\tbar,$ and where
$\taubar^2=\xibar^2+\abs{\eta}^2=1,$ and we have chosen the sign of
$\xibar$ to agree with the sign of $\taubar.$

As it is $Z$-invariant over the boundary, we may write $\hyp$ as the
pull-back to $\pa M$ via $\pi_0$ of a corresponding set $\dhyp.$ We may
therefore consider all the bicharacteristic meeting the boundary in in a
single fiber, with the same `slow variables' $(t,y)$ and set
$$
\fcaldot_{\bullet,q}= \bigcup_{p \in
    \pi_0^{-1}(q)} \fcal_{\bullet,p},\ q\in \dhyp.
$$
These pencils of bicharacteristics touching the boundary
at a given location in the `slow' spacetime variables $(t,y),$ with given
momenta in those variables, form smooth coisotropic (involutive) manifolds
in the cotangent bundle near the boundary. The two main theorems of this
paper show how microlocal singularities incoming along a bicharacteristic
$\fcal_{I,p}$ are connected to those along bicharacteristics $\fcal_{O,p'}$
for various values of $p,$ $p'.$

Our first main result is \emph{global} in the fiber, and corresponds to the
fact that the propagation of singularities is microlocalized in the slow
variables $(t,y)$ and their duals while in general being global in the
fiber. Let $\fcaldot^\circ_{\bullet,q}$ denote the part of
$\fcaldot_{\bullet,q}$ over the interior of $M$ and in an appropriately
small neighborhood of the boundary.

\begin{theorem}\label{theorem:1} For an admissible solution, $u,$ to the
wave equation and any $q\in\dhyp,$
\begin{equation*}
\fcaldot^\circ_{I,q} \cap \WF^k (u) =
\emptyset\Longrightarrow \fcaldot^\circ_{O,q} \cap \WF^k (u) =
\emptyset.
\label{30.8.2006.13}\end{equation*}
\end{theorem}

Thus singularities interact with the boundary by specular reflection,
preserving momentum in the slow variables.  The lack of localization in
the fibers is reflected by the fact that only the sets
$\fcaldot_{\bullet,q}$ rather than $\fcal_{\bullet,p}$ enter into the
statement. Refinements of this result giving appropriate regularity at the
boundary, not just in the interior, are discussed below.

We also prove a result regarding the behavior at glancing rays which
touch the boundary at $\gcal.$ In this case, singularities can propagate
only along continuations of rays as generalized broken bicharacteristics
much as in \cite{MR83h:35120}. These curves are defined in
\S\ref{section:edgeb} and in the non-glancing case are just unions of
incoming and outgoing bicharacteristics associated to the same point in
$\hcal,$ while in the glancing case the definition is subtler.  This allows
Theorem~\ref{theorem:1} to be extended to the following concise statement.

\begin{theorem} Singularities for admissible solutions
propagate only along generalized broken bicharacteristics.
\end{theorem}

We also obtain a result which is microlocal in the fiber variable as well
as in the slow variables but that necessarily has additional `nonfocusing'
hypotheses, directly generalizing that of \cite{Melrose-Wunsch1}, whose
detailed explanation we postpone to \S\ref{section:coisotropic}. This
condition, which is stated relative to a Sobolev space $H^s,$ amounts to
the requirement (away from, but near, the boundary) that the solution lie
in the image of $H^s$ under the action of a sum of products of first-order
pseudodifferential operators with symbols vanishing along $\fcaldot_{I,q}.$
Now consider two points $p,p'$ in $\dhyp$ and lying above $q,$ so in the
same fiber $Z_y.$ They are said to be \emph{geometrically related} if they
are the endpoints of a geodesic segment of length $\pi$ in $Z_y.$

\begin{theorem}\label{theorem:2} For an admissible solution $u$ satisfying
the nonfocusing condition relative to $H^s$ for $\fcaldot^\circ_{I,q},$ $q
\in\dhyp,$ and for $p\in\hcal$ projecting to $q,$
\begin{multline*}
\fcal^\circ_{I,p'} \cap \WF^s (u)=\emptyset \
\forall\ p' \text{ geometrically related to } p\Longrightarrow \\
\fcal^\circ_{O,p} \cap \WF^{r}(u) =\emptyset\ \forall\ r<s.
\label{30.8.2006.14}\end{multline*}
\end{theorem}

So, if we define `geometric generalized bicharacteristics' as unions of
bi\-characteristic segments, entering and leaving $\pa M$ only at
geometrically related points, then solutions satisfying an appropriate
nonfocusing condition have the property that the \emph{strongest}
singularities propagate only along geometric generalized
bicharacteristics. This may be visualized as follows: Consider an example
in which only \emph{one} singularity arrives at $\pa X,$ propagating along
$\fcal_{I,p'}$ Theorem~\ref{theorem:2} shows that this singularity is
diffracted into singularities which may emerge along the whole surface
$\fcaldot_{O,\fib(p')};$ these outgoing singularities are \emph{weaker}
than the incident singularity at all but the special family of
geometrically related bicharacteristics. The geometric generalized
bicharacteristics, along which stronger singularities can propagate, are
precisely those which can be obtained locally as limits of families of
bicharacteristics missing the edge entirely.  In the special case in which
$X$ is simply a smooth manifold $A$ blown up at a submanifold $B,$ the
geometric generalized bicharacteristics are just the lifts of
bichacteristics in the usual sense.  Condensation of singularities
arguments and the uniqueness in H\"ormander's propagation theorem show
that, at least locally, full-strength singularities do indeed propagate
along the geometric generalized bicharacteristics for some solutions.

The nonfocusing condition does hold for the forward fundamental solution
with pole $o$ sufficiently near the boundary
$$
u_o(t) =  \frac{\sin t\sqrt{\Lap}}{\sqrt\Lap} \delta_o.
$$ We have $u_o(t) \in H^{s}_{\loc}(M)$ for all $s<-n/2+1,$ while the
nonfocusing condition holds relative to $H^{s'}$ for all $s'<-n/2+1+f/2.$
If $p\in\hcal$ is sufficiently close to $o$ and $\fcal_{I,p}\cap\WF
U_o(t)\not=\emptyset,$ so the incoming singularity strikes the boundary at
$p,$ then for each $p'$ projecting to the same point in $\dhyp$ as $p,$ but
not geometrically related to it,
$$
\fcal_{O,p'}\cap \WF^{r} (u) = \emptyset\ \forall\ r<-n/2+1+f/2.
$$
Thus, the diffracted wave is almost $f/2$ derivatives smoother than the primary
singularities of the fundamental solution, so we obtain the following
description of the structure of $u_o:$
\begin{corollary}
For all $o \in X^\circ$ let $\mathcal{L}_o$ denote the flowout of
$SN^*(\{o\})$ along bicharacteristics lying over $X^\circ.$ If $o$ is
sufficiently close to $\pa X,$ then for short time, the
fundamental solution $u_o$ is a Lagrangian distribution along
$\mathcal{L}_o$ lying in $H^s$ for all $s<-n/2+1$ together with a
diffracted wave, singular only at $\fcal_O,$ that lies in $H^{r}$ for
all $r<-n/2+1+f/2,$ away from its intersection with $\mathcal{L}_o.$
\end{corollary}
The geometric generalized bicharacteristics in this case are just those at
the intersection of the diffracted wave and $\mathcal{L}_o.$

\subsection{Previous results}

It has been known since the work of the first author \cite{MR58:24409} and
Taylor \cite{Taylor1} that $\CI$-wavefront set of a solution to the wave
equation on a manifold with concave boundary and Dirichlet boundary
conditions does not propagate into the classical shadow region, that is,
that wavefront set arriving tangent to the boundary does not `stick to' the
boundary, but rather continues past it.\footnote{In the analytic category
the contrary is the case --- see \cite{MR55:10859},
\cite{Sjostrand2,Sjostrand3,Sjostrand4}.}  By contrast, diffractive effects
have long been known to occur for propagation of singularities on more
singular spaces.  The first rigorous example is due to Sommerfeld,
\cite{Sommerfeld1}, who analyzed the diffraction into the shadow region
behind a straight edge in two dimensions.  Many more such examples were
studied by Friedlander \cite{MR20:3703}, and a general analysis of the
fundamental solution to the wave equation on product cones was carried out
by Cheeger and Taylor \cite{Cheeger-Taylor2, Cheeger-Taylor1}.  Borovikov
\cite{Borovikov} has analyzed the structure of the fundamental solution on
polyhedra.  All of these works rely in an essential way on the technique of
separation of variables, using the product structure of the cone. There is
a quite general, but heuristic, geometric theory of diffraction due to
Keller \cite{Keller1}. This has been confirmed in a few special cases.

Non-product situations and more singular manifolds have been less explored.
G\'erard and Lebeau in \cite{MR93f:35130} explicitly analyzed the problem
of an analytic conormal wave incident on an analytic corner in $\RR^2,$
obtaining a $1/2$-derivative improvement of the diffracted wavefront over
the incident one, which is the analog of Theorem~\ref{theorem:2}
above. Also in the analytic setting, Lebeau \cite{Lebeau4, Lebeau5}
obtained a diffractive theorem analogous to Theorem~\ref{theorem:1} for a
broad class of manifolds, including manifolds with corners; the second
author \cite{Vasy5} has recently obtained such a theorem in the $\CI$
setting on manifolds with corners. As already noted, the first and third authors
\cite{Melrose-Wunsch1} have previously studied (non-product) conic
manifolds in the $\CI$ setting and obtained both Theorem~\ref{theorem:1}
and Theorem~\ref{theorem:2} in that setting.

\subsection{Proofs and plan of the paper}

The proofs of Theorems~\ref{theorem:1} and \ref{theorem:2} use two distinct
but related pseudodifferential calculi.  We have stated the results in
terms of regularity measured by the \emph{edge calculus,} which is best
suited for studying the fine propagation of regularity into and out of the
boundary along different bicharacteristics.  For arguments that are global
in the fiber, however, the \emph{b-calculus} proves to be a better tool,
and the proof of Theorem~\ref{theorem:1} is really a theorem about
propagation of b-regularity, following the argument in the corners setting
in \cite{Vasy5}. Essential use is made of the fact that the product $x\xi,$
with $\xi$ the the dual variable to $\pa_x,$ is increasing along the
bicharacteristic flow, and the test operators used are \emph{fiber
constant} so as not to incur the large error terms which would otherwise
arise from commutation with $\Lap.$

By contrast, as in \cite{Melrose-Wunsch1}, the extension of the proof of
H\"ormander's propagation theorem for operators with real principal symbols
to the edge calculus necessarily runs into the obstruction presented by
manifolds of \emph{radial points,} at which the Hamilton flow vanishes. The
subprincipal terms then come into play, and propagation results are subject
to the auxiliary hypothesis of \emph{divisibility}.  In particular,
propagation results into and out of the boundary along bicharacteristics in
the edge cotangent bundle up to a given Sobolev order are restricted by the
largest power of $x$ by which $u$ is divisible, relative to the
corresponding scale of edge Sobolev spaces.  The divisibility given by
energy conservation turns out to yield no useful information, as it yields
a regularity for propagated singularities that is less than the regularity
following directly from energy conservation. To prove
Theorem~\ref{theorem:2} we initially settle for less information. As
already noted, $\fcaldot_{I,q}$ is a \emph{coisotropic} submanifold of the
cotangent bundle and as such `coisotropic regularity' with respect to it
may be defined in terms of iterated regularity under the application of
pseudodifferential operators with symbols vanishing along $\fcaldot_{I,
\bullet}.$ Thus, we begin by showing that coisotropic regularity in this
sense, of any order, propagates through the boundary, with a fixed loss of
derivatives.  By interpolation with the results of Theorem~\ref{theorem:1}
it follows that coisotropic regularity propagates through the boundary with
epsilon derivative loss.  Finally, imposing the nonfocusing condition
allows Theorem~\ref{theorem:2} to be proved by a pairing argument, since
this condition is a microlocal characterization of the dual to the space of
distributions with coisotropic regularity.

\section{Edge manifolds}\label{section:edge}

Let $X$ be an edge manifold, as defined in \S\ref{intro:results}.  In this
section, we analyze the geodesic flow of $X$ in the edge cotangent bundle,
and use this flow to describe a normal form for edge metrics.

As remarked above, by using the product decomposition near the boundary
given by the distance along the normal geodesics in the base manifold
$[0,\ep) \times Y,$ we may write an edge metric in the form
\eqref{edgemetric2}. The following result allows us to improve the form of
the metric a bit further, so that $dx$ arises only in the $dx^2$ term:
\begin{proposition}\label{prop:normform}
On an edge manifold, there exist choices for $x$ and for the product
decomposition $U \equiv [0,\ep) \times Y$ such that
$$
g = dx^2 + h(x,y,dy)+ x h'(x,y,z,dy) + x^2 k(x,y,z,dy,dz)
$$
\end{proposition}

To prove Proposition~\ref{prop:normform} (and other results to come)
it is natural to study the Hamilton flow associated to $g$ in the edge
cotangent bundle of $X.$ Let
$$
\alpha = \xi \frac{dx}x + \eta \cdot \frac{dy}x + \zeta \cdot dz
$$ denote the canonical one-form, lifted to $\Testar X.$ The dual basis of
vector fields of $\Te X$ is $x \pa_x,$ $x \pa_y,$ $\pa_z$ and in terms of this
basis the metric takes the form
$$
g = x^2 \begin{pmatrix}
1 & O(x) & O(x) \\
O(x) & H+O(x) & O(x) \\
O(x) & O(x) & K
\end{pmatrix}
$$ where $H$ and $K$ (which are nondegenerate) are defined respectively as
the $dy \tensor dy$ and $dz \tensor dz$ parts of $h$ and $k$ at $x=0.$ Here
and henceforth we employ the convention that an $O(x^k)$ term denotes $x^k$
times a function in $\CI(X).$

Thus, 
\begin{equation}\label{ginverse}
g^{-1} = x^{-2} \begin{pmatrix}
1+O(x) & O(x) & O(x) \\ O(x) & H^{-1} +O(x) & O(x) \\ O(x) &
O(x) & K^{-1}+ O(x)
\end{pmatrix}
\end{equation}

In fact, it is convenient to work on $M=X\times\RR_t$, which has an induced
boundary fibration, with base $Y\times\RR_t$. Correspondingly, elements of
$\Testar M$ are written as
\begin{equation*}
\xi\frac{dx}{x}+\tau\frac{dt}{x}+\eta\cdot\frac{dy}{x}+\zeta\cdot dz,
\end{equation*}
so $(x,t,y,z,\xi,\tau,\eta,\zeta)$ are local coordinates on $M$.
The canonical symplectic form on $\Testar M^\circ$ is $d(\tau
dt/x+ \alpha)$, and
the Hamilton vector field of the symbol of the wave operator $\Box$,
$p=(\tau/x)^2-g^{-1}$, is then $-2 x^{-2}H,$
with
\begin{multline}\label{eq:H-form}
H=-\tau x \pa_t + \left(\tau \xi + O(x) \right)\pa_\tau +
\left(x\xi+O(x^2)\right)\pa_x \\ + \left(\xi^2 + \zeta_i \zeta_j
K^{ij}+O(x)\right)\pa_\xi + \left(\zeta_i K^{ij} + O(x)\right) \pa_{z_j}
\\+ \left(-\h \zeta_i\zeta_j \frac{\pa K^{ij}}{\pa z_k} +
O(x)\right)\pa_{\zeta_k} +\left(x \eta_j H^{ij}+ O(x^2)\right) \pa_{y_i} 
+ \left(\xi\eta_i + O(x)\right)\pa_{\eta_i}.
\end{multline}
(We employ the superscript $\circ$ to denote the interior of a manifold.)

Note that while it has not been emphasized in the expression above, this
vector field is homogeneous of degree $1$ in the fiber variables. We will
be interested in the restriction of this flow to $\Sigma,$ the
characteristic variety of $\Box,$ where $(\tau/x)^2= g^{-1}.$

\begin{lemma}\label{30.8.2006.16} Inside $\Sigma$,
$H$ is {\em radial}, i.e.\ is tangent to the orbits of
the fiber dilations on $\Testar M\setminus 0$,
if and only if  $x=0$ and $\zeta=0$.
\end{lemma}

\begin{proof} To be radial $H$ must be a multiple of $\xi\pa_\xi+\tau\pa_\tau
+\eta\cdot\pa_\eta+\zeta\cdot\pa_\zeta.$ In particular, in the expression above,
the $\pa_t$ component must vanish, which inside $\Sigma$ implies $x=0.$
The vanishing of the $\pa_z$ component further implies $\zeta=0$, and
conversely $H$ is indeed radial at $x=0,$ $\zeta=0.$
\end{proof}

As usual, it is convenient to work with the cosphere bundle, $\Sestar M,$
viewed as the boundary `at infinity' of the radial compactification
of $\Testar M.$ Introducing the new variable
$$
\sigma = \abs\tau^{-1},
$$
it follows that $H\sigma = -\xi\sigma+O(x).$
Setting
\begin{equation*}
\xis = \xi\sigma,\ \etas = \eta\sigma,\ \zetas =\zeta \sigma
\label{30.8.2006.15}\end{equation*}
we note that $\sigma \xis,$ $\etas,$ $\zetas$ are coordinates on the fiber
compactification of $\Testar M \cap \Sigma,$ with $\sigma$ a defining
function for $\Sestar M$ as the boundary at infinity. Then
Lemma~\ref{30.8.2006.16} becomes
$$
\text{Inside }\Sestar M \cap \Sigma,\ \sigma H\text{ vanishes exactly at
}x=0,\ \zetas=0.
$$

Let the linearization of $\sigma H$ at $q\in \Sigma\cap\Sestar M$ (where
$x=0,$ $\zetas=0$) be $\Lin_q.$

\begin{lemma}\label{lemma:linearization}
For $q\in\dhyp$, i\@.e\@. such that $\xis(q)\neq 0,$ the eigenvalues of
$\Lin_q$ are $-\xis$, $0$ and $\xis$, with $dx$ being an eigenvector of
eigenvalue $\xis$, and $d\sigma$ being an eigenvector of eigenvalue
$-\xis$. Moreover, modulo the span of $dx$, the $-\xis$-eigenspace is
spanned by $d\sigma$ and the $d\zetas_j,$ and the $0$-eigenspace is spanned
by $dt$, $dy_j$, $d\etas_j$, $d\xis,$ and $dz_j+\xis^{-1}\sum_i
K^{ij}(q)\,d\zetas_i.$
\end{lemma}

\begin{remark}\label{rem:linearization}
This shows in particular that the span of the $d\zetas_j$ (plus a suitable
multiple of $dx$) is 
invariantly given as the stable/unstable eigenspace of $\Lin_q$ inside
$T^*_q\Sestar M$ according to $\xis>0$ or $\xis<0$. We denote this subspace
of $T^*_q\Sestar M$ by $T^{*,-}_q(\Sestar M).$
\end{remark}

\begin{proof}[Proof of Lemma~\ref{lemma:linearization}]
A straightforward calculation using \eqref{eq:H-form} yields,
in the coordinates $(x,y,t,z,
\xis,\etas,\sigma,\zetas)$,
\begin{equation*}
\sigma H=\xis(x\pa_x-\sigma\pa_\sigma-\zetas\pa_{\zetas})
+K^{ij}\zetas_i\pa_{z_j}+K^{ij}\zetas_i\zetas_j\pa_{\xis}
-\frac{1}{2}\,\frac{\pa K^{ij}}{\pa z_k}\zetas_i\zetas_j\pa_{\zetas_k}
+xH',
\end{equation*}
with $H'$ tangent to the boundary. Thus,
\begin{equation*}
\sigma H=\xis(x\pa_x-\sigma\pa_\sigma-\zetas\pa_{\zetas})
+K^{ij}\zetas_i\pa_{z_j}
+xH'+H'',
\end{equation*}
where $H'$ is tangent to the boundary, and $H''$ vanishes quadratically
at $q$ (as a smooth vector field). In particular, the linearization
$\Lin_q$ is independent of the $H''$ term, and is only affected by
$xH'$ as an operator with (at most) one-dimensional range, lying in
$\Span\{dx\}$. Moreover, $xH'x$ vanishes quadratically at $\pa M$,
so $dx$ is an eigenvector of $\Lin_q$, and by the one-dimensional range
observation, modulo $\Span\{dx\}$, the other eigenspaces can be read
off from the form of $\sigma H$, proving the lemma.
\end{proof}

Therefore, under the flow of $\sigma H$ we have
\begin{equation}\label{hyppoint}
\begin{aligned}
\zetas' &= -\xis \zetas+ O(\zetas^2) + O(x),\\
x' &= \xis x + O(x^2).
\end{aligned}
\end{equation}
The Stable/Unstable Manifold Theorem can then be applied to the flow under
$H$ and shows that for each $(y_0,z_0) \in \pa X,$ $\xi_0,$ $\eta_0,$ with
$\xi_0 \neq 0,$ there exists a unique bicharacteristic with the point
$x=\zetas=0,$ $(y,z)=(y_0,z_0),$ $\xis=\xi_0,$ $\etas=\eta_0$ in its
closure and that a neighborhood of the boundary is foliated by such
trajectories. Choosing $\xi_0=1,$ $\etas=0$ gives a foliation by normal
trajectories arriving perpendicular to $\pa X.$ Let $\tx$ be the distance
to the boundary along the corresponding normal geodesics; it is a boundary
defining function. Gauss's Lemma implies that in terms of the induced
product decomposition the metric takes the form
$$
g=d\tx^2 + h(\tx,\ty, d\ty)+\tx h'(\tx,\ty,\tz,d\ty,d\tz) +
\tx^2 k(\tx,\ty,\tz, d\ty, d\tz).
$$
This concludes the proof of Proposition~\ref{prop:normform}.

This argument in fact gives a little more than
Proposition~\ref{prop:normform}, as we may take $\xi_0 \neq 1.$ In fact,
for any $q=(t,y,z,\bar\tau, \bar\xi, \bar\eta)\in \hyp$ there exist two
bicharacteristics having this point as its limit in the
boundary,\footnote{The notation $\hyp,$ $\dhyp$ will be further elucidated
in \S\ref{section:edgeb}.} which we denote
$$
\fcal_{I/O, q},
$$ with the choice of $I/O$ determined by $\sgn \xi/\tau.$ For $p =
(t,y,\bar\tau, \bar\xi,\bar\eta)\in \dhyp,$ we let
$$
\pi_0^{-1}(p) = \{ (t,y,z,\bar\tau, \bar\xi, \bar\eta): z \in Z_y\}
$$
and
$$
\fcaldot_{I/O,p}= \bigcup_{q \in \pi_0^{-1}(p)}\fcal_{I/O,p},
$$
and we collect all these bicharacteristics into
$$
\fcal_{I/O} = \bigcup_{p} \fcal_{I/O,p},\quad \fcal = \fcal_I\cup \fcal_O.
$$

For later use, we note that when the metric has been reduced to the normal
form guaranteed by Proposition~\ref{prop:normform}, the rescaled Hamilton
vector field on $\Testar (M)$ reduces to
\begin{multline}\label{hamvf.normform}
H=-\tau x \pa_t + \tau \xi \pa_\tau + \xi x \pa_x + (\xi^2 +
\abs{\zeta}_{\kbar}^2)\pa_\xi + (\zeta_i \kbar^{ij} + O(x)) \pa_{z_j}\\ +
(-\h \zeta_i\zeta_j \frac{\pa \kbar^{ij}}{\pa z_k} + O(x))\pa_{\zeta_k} +(x \eta_j H^{ij} + O(x^2)) \pa_{y_i} + (\xi\eta_i +O(x))\pa_{\eta_i};
\end{multline}
here we have let $\kbar^{ij}$ denote the whole of the $K^{-1}+ O(x)$
term in the bottom right block of $g^{-1}$ as expressed in
\eqref{ginverse}.

\section{Calculi}\label{section:calculi}
In order to simplify the descriptions below, we assume that $M$ is compact.
If $M$ is non-compact, we must insist on all operators having properly
supported Schwartz kernels. For such operators all statements but the
$L^2$-boundedness remain valid; when we restrict ourselves to operators
with \emph{compactly} supported Schwartz kernels, we have $L^2$-boundedness
as well.

We use two algebras of pseudodifferential operators. Each arises
from a Lie algebra $\Vf$ of smooth vector fields on $M$ tangent to $\pa M$,
which is a $\CI(M)$-module,
and they have rather similar properties:

\begin{itemize}
\item
For b-ps.d.o's, we take
$\Vf=\Vb(M)$, consisting of all $\CI$ vector fields on $M$ that are
tangent to $\pa M$.
\item
For edge ps.d.o's, $\pa M$ has a fibration $\phi:\pa M\to Y,$ and we take
$\Vf=\Ve(M)$, consisting of all $\CI$ vector fields on $M$ that are
tangent to the fibers of $\phi$ (hence in particular to $\pa M$).
\end{itemize}
In terms of coordinates $x,y,z$ adapted to a local trivialization of the
fibration on an edge manifold, with $x$ a boundary defining function and
$z$ a coordinate in the fiber, it is easy to verify that $\Vb(M)$ is
spanned over $\CI(M)$ by $x \pa_x,$ $\pa_{y_j}$, $\pa_{z_k}$ while $\Ve(M)$
is spanned over $\CI(M)$ by $x \pa_x,$ $x \pa_{y_j}$, $\pa_{z_k}.$

In both cases, $\Vf$ is the space of all smooth sections of a vector bundle
$\VT M\to M$.  Thus, in the two cases, these bundles are denoted by $\bT M$
and $\eT M$ respectively; this notation extends in an obvious manner to all
$\Vf$-objects we define below.  The dual vector bundle is denoted by
$\VTstar M$, and the corresponding cosphere bundle, i.e.\ $(\VTstar
M\setminus 0)/\RR^+$ is denoted by $\VSstar M$. As the space of all smooth
vector fields on $M$ includes $\Vf$, there is a canonical bundle map $\VT
M\to T M$, and a corresponding dual map, $T^*M\to \VTstar M$. These are
isomorphisms over the interior of $M$. Also, being a $\CI$ section of $\VT
M$, every $V\in\Vf$ defines a linear functional on each fiber of $\VTstar
M$.  This is the principal symbol map:
\begin{equation*}
\Vsigma_1:\Vf\to\ \CI\ \text{fiber-linear functions on}\ \VTstar M.
\end{equation*}

We let $\DiffV{k}(M)$ be the space of differential operators generated
by $\Vf$ over $\CI(M)$, so elements of $\DiffV{k}(M)$ are finite
sums of terms of the form $aV_1\ldots V_l$, $l\leq k$, $V_j\in\Vf$,
$a\in\CI(M)$. Thus, $\DiffV{}(M)=\bigcup_k\DiffV{k}(M)$ is a filtered
algebra over $\CI(M)$. Defining
$\Vsigma_0(a)=\pi^* a$ for $a\in\CI(M)$,  $\pi$ being the bundle projection
$\VTstar M\to M$, $\Vsigma$ extends to a map
\begin{multline*}
\Vsigma_m:\DiffV{m}(M)\longrightarrow\\
\big\{\CI\ \text{fiber-homogeneous polynomials of
degree}\ m\ \text{on}\ \VTstar M\big\},
\end{multline*}
that is a filtered ring-homomorphism in the sense that
\begin{equation*}
\Vsigma_{m+l}(AB)=\Vsigma_m(A)
\Vsigma_l(B)\ \text{for}\ A\in\DiffV{m}(M),\ B\in\DiffV{l}(M).
\end{equation*}

We also fix a non-degenerate $b$-density $\nu$ on $M$, hence $\nu$ is of
the form $x^{-1}\nu_0$, $\nu_0$ a non-degenerate $\CI$ density on $M$,
i.e.\ a nowhere-vanishing section of the density bundle $\Omega M:=
\abs{\bigwedge^n}(M).$ The density gives an inner product on $\dCI(M)$. When
below we refer to adjoints, we mean this relative to $\nu$, but the
statements listed below not only do not depend on $\nu$ of the stated form,
but would even hold for any non-degenerate density $x^{-l}\nu_0$, $\nu_0$
as above, $l$ arbitrary, as the statements listed below imply that
conjugation by $x^l$ preserves the calculi.

For each of these spaces $\Vf$ there is an algebra of
pseudodifferential operators, denoted by $\VPs{}(M)$
which is a bifiltered $\CI(M)$
*-algebra of operators acting on $\dCI(M)$,
$\VPs{}(M)=\bigcup_{m,l}\VPs{m,l}(M)$, satisfying
\renewcommand{\theenumi}{\Roman{enumi}}
\begin{enumerate}
\item
$\DiffV{m}(M)\subset\VPs{m,0}(M)$
(with the latter also denoted by $\VPs{m}(M)$),
\item
if $x$ is a boundary defining function of $M$ then $x^l\in\VPs{0,l}(M)$, and
$x^l\VPs{m}(M)=\VPs{m,l}(M)$,
\item
the principal symbol map, already defined above
for $\DiffV{}(M)$, extends to a bifiltered *-algebra homomorphism
\begin{equation*}
\Vsigma_{m,l}:\VPs{m,l}(M)\to x^l S^m_{\text{hom}}(\VTstar M\setminus 0),
\end{equation*}
\item
the principal symbol sequence is exact
\begin{equation*}
0\to\VPs{m-1,l}(M)\hookrightarrow\VPs{m,l}(M)\to
x^l S^m_{\text{hom}}(\VTstar M\setminus 0)\to 0,
\end{equation*}
\item
for $A\in\VPs{m,l}(M)$, $B\in\VPs{m',l'}(M)$, $[A,B]\in\VPs{m+m'-1,l+l'}(M)$
satisfies
\begin{equation*}
\Vsigma_{m+m'-1,l+l'}([A,B])
=\frac{1}{i}\{\Vsigma_{m,l}(A),\Vsigma_{m',l'}(B)\},
\end{equation*}
with the Poisson bracket defined a priori as a homogeneous function
on $T^*M^\circ\setminus 0$, but extends to $x^{l+l'}$ times a smooth
function on $\VTstar M\setminus 0$,
\item
every $A\in\VPs{0,0}(M)$ extends from $\dCI(M)$ by continuity
to define a continuous linear map on $L^2(M)$, and
for each $A\in\VPs{0,0}(M)$ there exists $A'\in\VPs{-1,0}(M)$ such that
for all $u\in L^2(M)$,
\begin{equation*}
\|Au\|_{L^2}\leq (2\sup|a|)\|u\|_{L^2}+\|A'u\|_{L^2},
\end{equation*}
where $a=\Vsigma_{0,0}(A)$.
\end{enumerate}

In addition,
there is an operator wave front set (or microsupport)
$\VWF'$ such that for $A\in\VPs{}(M)$,
$\VWF'(A)\subset\VSstar M$ is closed, and satisfies
\renewcommand{\theenumi}{\Alph{enumi}}
\begin{enumerate}
\item
$\VWF'(AB)\subset\VWF'(A)\cap\VWF'(B),$ 
\item $\VWF'(A+B)\subset\VWF'(A)\cup
\VWF'(B),$ 
\item $\VWF'(a)=\pi^{-1}\supp a$ for $a\in\CI(M),$ where $\pi$ is
the bundle projection,
\item\label{operatorexists}
for any $K\subset\VSstar M$ closed and $U\subset\VSstar M$ open with
$K\subset U$, there exists $A\in\VPs{0,0}(M)$ with
$\VWF'(A)\subset U$ and $\Vsigma_{0,0}(A)=1$ on $K$ (so $\VPs{}(M)$
is microlocal on $\VSstar M$),
\item
if $A\in\VPs{m,l}(M)$, $\Vsigma_{m,l}(A)(q)\neq 0$, $q\in\VSstar M$, then
there exists a microlocal parametrix $G\in\VPs{-m,-l}(M)$ such that
\begin{equation*}
q\notin\VWF'(GA-\Id),\ q\notin\VWF'(AG-\Id),
\end{equation*}
\item
if $A\in\VPs{m,l}(M)$ and $\VWF'(A)=\emptyset$ then $A\in\VPs{-\infty,l}
=\bigcap_{m'}\VPs{m',l}$ -- thus, $\VWF'$ captures $A$ modulo operators
of order $-\infty$ in the smoothing sense, just as $\Vsigma_{m,l}$
captures $A$ modulo operators of order $(m-1,l).$
\end{enumerate}
\renewcommand{\theenumi}{\arabic{enumi}}

Note that if $b\in x^l S^m_{\text{hom}}(\VTstar M\setminus 0)$ then for
any $U\subset\VSstar M$ open satisfying $U\supset\supp b$ there
exists $B\in\VPs{m,l}(M)$ with $\VWF'(B)\subset U$ and
$\Vsigma_{m,l}(B)=b$. Indeed, with $K=\supp b$, let $A\in\VPs{0,0}(M)$
be as in \eqref{operatorexists} above,
so $\VWF'(A)\subset U$ and $\Vsigma_{0,0}(A)=1$ on $K$. By the exactness
of the symbol sequence, there is $B_0\in\VPs{m,l}(M)$ and
$\Vsigma_{m,l}(B)=b$. Then $B=AB_0$ satisfies all requirements.

Beyond these rather generic properties, shared by the b- and edge-calculi,
we will use one further, more specialized property of the b-calculus:
if $V\in\Vb(M)$ is such that $V|_{\pa M}=xD_x$ (restriction in
the sense of $\Vb(M)$, i.e.\ as sections of $\bT M$), then for
$A\in\Psib^m(M)$, $[A,V]\in x\Psib^m(M)$, i.e.\ there is a gain of $x$ over
the a priori statements. (Note that $x\pa_x$, the `radial vector field', is
a well-defined section of $\bT M$ at $\pa M$.)  

The above property in fact follows from the analysis of the \emph{normal
operator} associated to $V.$ More generally, for a differential operator in
the $b$-calculus, we obtain a model, or \emph{normal} operator by freezing
coefficients of b-vector fields at the boundary: in coordinates $(x,w) \in
\RR_+ \times \RR^k$ with $x$ a boundary defining function, if
$$
P=\sum a_{i,\alpha}(x,w) (x D_x)^i D_w^\alpha
$$
then
$$
N(P)= \sum a_{i,\alpha}(0,w) (x D_x)^i D_w^\alpha
$$
is now an operator on $\RR_+\times\pa M.$  Likewise, if $Q$ is an
\emph{edge} operator given by
$$
Q=\sum b_{i,\alpha,\beta}(x,y,z) (x D_x)^i (x D_y)^\alpha D_z^\beta
$$
on a manifold with fibered boundary
then
$$
N(Q)= \sum b_{i,\alpha,\beta}(0,y,z) (x D_x)^i (x D_y)^\alpha D_z^\beta
$$
is an operator on $\RR^{b+1}_+\times F,$ where $F$ is the fiber and $b$ the
dimension of the base of the boundary fibration.  These normal operators
are in fact homomorphisms from the respective algebras, and they extend to
act on the algebras $\VPs{m}(M).$  Their principal utility
is the following:
$$
A \in x\VPs{m}(M) \Longleftrightarrow N(A)=0.
$$
We will make very little use of the normal operators; for more detailed
discussion, see \cite{MR93d:58152}.

We will also require a conormal variant of the b-calculus described
above, where polyhomogeneous symbols ($S^m_{\text{hom}}$) are replaced by
those satisfying Kohn-Nirenberg type symbol estimates ($S^m$).  We denote
this calculus $\Psibc.$ It satisfies all the properties above, except that
the symbol map now takes values in $S^m(\Tbstar M)/S^{m-1}(\Tbstar M).$
This calculus will be important in making approximation arguments in our
b-calculus-based positive commutator estimates.

A unified treatment of the two calculi discussed in this section can be
found in Mazzeo \cite{MR93d:58152}.  This was the original treatment of the
edge calculus, while the $b$-calculus originates with the work of the first
author \cite{MR83f:58073}.

\section{Coisotropic regularity}\label{sec:coiso-reg}

\begin{theorem}
Away from glancing rays, the set $\fcal$ is a coisotropic submanifold of
the symplectic manifold $\Testar M,$ i.e.\ contains its symplectic
orthocomplement.
\end{theorem}
\begin{proof}
We split $\fcal$ into its components $\fcal_I$ and $\fcal_O;$ by symmetry,
it suffices to work on $\fcal_I.$

It suffices to exhibit a set of functions $\alpha_i$ vanishing on $\fcal_I$
whose differentials yield a basis of sections of $N^* \fcal_I$ such that
the Poisson bracket of any two also vanishes on $\fcal_I.$

To begin, we let $H_{p_0}$ denote the Hamilton vector field of
$p_0=\esigma_{2,-2}(\Box)$ (and more generally, $H_q$ the Hamilton vector
field of any symbol $q$); let $\varphi = x^2/\abs{\tau}$ so $\varphi
H_{p_0} \in \mathcal{V}_b(\Testar (M)\backslash 0),$ \emph{and is
homogeneous of degree zero}.  The vector field $\varphi H_{p_0}$ is not
itself a Hamilton vector field: we have
\begin{equation}\label{nonham}\varphi H_{p_0} = H_{\varphi p_0} - p_0
H_\varphi.\end{equation} This last vector field has the virtue of vanishing
on $\fcal_I,$ however.

Now $\fcal_I$ is a conic submanifold of $\Testar M,$ hence has a set of
defining functions that are homogeneous of degree $0.$ Letting $\zetas_i =
\zeta_i/\abs{\tau}$ we note that it follows from \eqref{hyppoint} that we
may take defining functions for $\fcal_I$ to be homogeneous functions of
degree zero of the form
$$
\alpha_i = \zetas_i+O(x), \quad i=1,\dots, f,
$$ with $d\alpha_i(q)$ an eigenvector of the linearization $A_q$ of $\varphi H_{p_0}$
of eigenvalue
$-\hat\xi$,
plus one extra to keep us in the characteristic set: $\alpha_0 =
x^2p_0/\abs{\tau}=\varphi p_0.$  We immediately note that for $i\geq 1,$ we
have
\begin{equation}\label{eq:Ham-a_0-a_i}
\{\alpha_0, \alpha_i\} = \varphi H_{p_0} \alpha_i + p_0 H_\varphi \alpha_i = p_0 H_\varphi (\alpha_i),
\end{equation}
which is a smooth multiple of $p_0,$ so that Poisson brackets with $\alpha_0$
never present difficulties.  It remains to show that $\{\alpha_i,\alpha_j\}$ is in
the span of the $\alpha_i$'s for $i=1,\dots,f.$

The Jacobi identity yields
\begin{equation}\label{Jacobi}
H_{\varphi p_0} \{x\alpha_i,x\alpha_j\} = \{H_{\varphi p_0}(x\alpha_i), x\alpha_j \} - \{H_{\varphi
  p_0}(x\alpha_j), x\alpha_i\}, \quad i,j=1,\dots, f.
\end{equation}
Using \eqref{nonham} and the tangency of $\varphi H_{p_0}$ to $\fcal_I$ and
to $x=0$ we find
that $H_{\varphi p_0}(x\alpha_i)$ is a smooth function vanishing at both $x=0$
and on $\fcal_I,$ hence
we have
$$
H_{\varphi p_0}(x \alpha_i) =  \sum b_{ik}' x \alpha_k
$$ where $b_{ik}' \in \CI.$ We can in fact say a little more: since
$\varphi H_{p_0}$ is a smooth vector field tangent both to $x=0$ and to
$\fcal_I,$ and as $d\alpha_i(q)$, resp.\ $dx$ are eigenvectors of the
linearization of $\varphi H_{p_0}$ at $\pa \fcal_I$ with eigenvalue $-\hat\xi$,
resp.\ $\hat\xi$, we have
\begin{equation*}
\varphi H_{p_0} x=\hat\xi x+x^2 b+\sum c_j x\alpha_j,
\ \varphi H_{p_0} \alpha_i=-\hat\xi \alpha_i+x\sum r_{ij}\alpha_j+r_i,
\end{equation*}
with $b,c_i,r_{ij},r_j$ all smooth, $r_i$ vanishing quadratically at
$\fcal_I.$  Thus,
$$
\varphi H_{p_0}(x \alpha_i) =  \sum b_{ik} (x^2 \alpha_k) + r'_i
$$
for smooth functions $b_{ik}$ and $r'_i$ and with $r'_i$ quadratically
vanishing on $\fcal_I.$

Observe now that $p_0H_{\varphi}(x\alpha_i)=(\varphi p_0)xc_i$ with $c_i$ smooth,
and $\{(\varphi p_0)xc_i,x\alpha_j\}$ vanishes at $p_0=0$ (cf.\
\eqref{eq:Ham-a_0-a_i}).  As a result, by \eqref{nonham} and
\eqref{Jacobi}, we see that
$$
\varphi H_{p_0} \{x\alpha_i,x\alpha_j\} = \sum x c_{ijkl} \{x\alpha_k, x\alpha_l\}+r
$$
where $r$ is a smooth term vanishing on $\fcal_I.$

Now we restrict to $\fcal_I.$ By construction of $\fcal_I,$ the vector
field $\varphi H_{p_0}$ vanishes identically at $\fcal_I \cap \pa M,$ hence
it is divisible by a factor of $x$ when restricted to $\fcal_I,$ with the
quotient being a vector field on $\fcal_I$ transverse to $\pa \fcal_I.$
Thus if we let $s$ be a parameter along the flow generated by $\varphi
H_{p_0},$ with $s=0$ at $\pa \fcal_I,$ we have with $c'_{ijkl}$ smooth,
$$
\frac{d}{ds} \{x\alpha_i,x\alpha_j\}|_{\fcal_I} = \sum c'_{ijkl} \{x\alpha_k, x\alpha_l\}|_{\fcal_I}
$$
Now $x \alpha_k$ and $x \alpha_l$ Poisson commute at $\pa M,$ so we conclude that
$\{x \alpha_i, x \alpha_j\}=0$ on $\fcal_I.$

This implies that in fact $\{\alpha_i, \alpha_j\}=0$ on $\fcal_I,$ as $\{x \alpha_i, x
\alpha_j\}-x^2\{\alpha_i,\alpha_j\}$ also vanishes there.
\end{proof}

We now fix an arbitrary open set $U \subset \Testar M$ disjoint from rays
meeting $x=\xi=0.$ We omit the set $U$ from the notation below, by abuse
of notation.

\begin{definition}\label{definition:coisotropic}
\renewcommand{\theenumi}{\alph{enumi}}

\ 

\begin{enumerate}
\item
Let $\ePs{}(U)$ be the subset of $\ePs{}(M)$ consisting of operators $A$
with $\eWF' A\subset U$.
\item
Let $\module$ denote the module of pseudodifferential operators in
$\ePs{1}(M)$ given by
$$\module = \{A \in \ePs{1}(U): \sigma(A)|_{\fcaldot}=0\}.$$

\item
Let $\coiso$ be the algebra generated by $\module,$ with $\coiso^k = \coiso
\cap \ePs{k}(M).$
Let $\hilbert$ be a Hilbert space on which $\ePs{0}(M)$ acts, and let $K\subset
\Testar (M)$ be a conic
set.
\item\label{coisocond} We say that $u$ has coisotropic regularity of order $k$ relative to $\hilbert$ in
$K$ if there exists $A \in \ePs{0}(M),$ elliptic on $K,$ such that
$\coiso^k A u \subset \hilbert.$
\item\label{nfcond} We say that $u$ satisfies the nonfocusing condition of degree $k$
  relative to $\hilbert$ on $K$ if there exists $A \in \ePs{0}(M),$ elliptic on
  $K,$ such that $A u \subset \coiso^k \hilbert.$  We say that $u$ satisfies the
  nonfocusing condition relative to $\hilbert$ on $K$ if it satisfies the
  condition to some degree.
\end{enumerate}
\renewcommand{\theenumi}{\Roman{enumi}}
\end{definition}
\begin{remark}
The conditions of coisotropic regularity and nonfocusing are dual
to one another, at least in a microlocal sense: if $u$ satisfies
\eqref{coisocond} and $v$ satisfies \eqref{nfcond} then the $L^2$ pairing
of $Au$ and $Av$ makes sense.
\end{remark}

\begin{lemma}\label{lemma:test-module}
The module $\module$ is a {\em test module}
in a sense analogous to \cite[Definition~6.1]{MR2020655}.
That is, $\module$ is closed under commutators and is finitely generated
in the sense that there exist finitely many $A_i\in\ePs{1}(M)$,
$i=0,1,\ldots,N$, $A_0=\Id$, such that
\begin{equation*}
\module=\{A\in\ePs{1}(U):\ \exists Q_i\in\ePs{0}(U),
\ A=\sum_{j=0}^N Q_iA_i\}
\end{equation*}
Moreover, we may take $A_N$ to have symbol $|\tau|\,a_N=x^2|\tau|^{-1}p$,
$p=\sigma(\Box)$, and $A_i$, $i=1,\ldots,N-1$ to have symbol $|\tau|\,a_i$
with $da_i(q)\in T^{*,-}_q(\Sestar M)$ for $q\in\pa\fcaldot$, where
we used the notation of Remark~\ref{rem:linearization}.
\end{lemma}

\begin{proof}
For  $A,B\in\module$, $[A,B]\in\ePs{1}(M)$, $\eWF'([A,B])\subset U$ and
$\sigma([A,B])=\frac{1}{i}\{\sigma(A),\sigma(B)\}$, so
$\sigma([A,B])|_{\fcaldot}=0$ as $\fcaldot$ is coisotropic. Thus,
$\module$ is closed under commutators.

Let $a_i$, $i=1,\ldots,N$,
be smooth
homogeneous degree $0$ functions on $\Testar M\setminus 0$ such
that $a_i$ vanish on $\fcaldot$ and $da_i(q)$, $i=1,\ldots,N$, spans
$N^*_q\fcaldot$ for each $q\in \fcaldot$. In particular, we
may take $a_N=x^2|\tau|^{-1}p_0=\varphi p_0$ and choose $a_i$, $i=1,\ldots,N-1$ as
in the statement of the lemma. Let $A_i$ be
any elements of $\ePs{1}(M)$ with $\sigma(A_i)=|\tau|\,a_i$. In view of this
spanning property, every homogeneous degree $1$ function
$a\in\CI(\Testar M\setminus 0)$ with $\supp a\subset U$ and $a|_\fcaldot=0$
can be written as $\sum_{i=1}^N q_i\, |\tau|\,a_i$ with $\supp q_i\subset U$,
$q_i$ smooth homogeneous degree $0$, so if $Q_i\in\ePs{0}(M)$ are
chosen to satisfy $\eWF'(Q_i)\subset U$ and $\sigma(Q_i)=q_i$, then
$Q_0=A-\sum_{i=1}^N Q_i A_i$ satisfies $\eWF'(Q)\subset U$,
$Q_0\in\ePs{0}(M)$, proving the lemma.
\end{proof}

We recall from \cite{MR2020655}:

\begin{lemma}\label{lemma:coiso-gens}(See \cite[Lemma~6.3]{MR2020655}.)
If $A_i,$ $0\le i\le N,$ are generators for
$\module$ in the sense of Lemma~\ref{lemma:test-module} with $A_0=\Id,$ then 
\begin{equation}
\coiso^{k}=\left\{\sum\limits_{|\alpha|\le k} Q_\alpha
\prod\limits_{i=1}^{N}A_i^{\alpha _i},\ Q_\alpha \in\ePs{0}(U)\right\}
\label{HMV.102}\end{equation}
where $\alpha$ runs over multiindices $\alpha
:\{1,\dots,N\}\to\bbN_0$ and $|\alpha |=\alpha _1+\dots+\alpha _N.$ 
\end{lemma}

\begin{remark}
The notation here is that the empty product is $A_0=\Id$, and the product
is ordered by ascending indices $A_i$. The lemma is
an immediate consequence of $\module$ being both a Lie algebra and a
module; the point being that products may be freely rearranged, modulo
terms in $\coiso^{k-1}$.
\end{remark}

Now, $A_N$ being a multiple of $\Box$, modulo $\ePs{0}(M)$,
regularity under summands $Q_\alpha
\prod\limits_{i=1}^{N}A_i^{\alpha _i}$ with $\alpha_N\neq 0$ and $|\alpha|=k$
is automatic for solutions of $\Box u=0$,
once regularity with respect to $\coiso^{k-1}$ is known. The key additional
information we need is:

\begin{lemma}
For $l=1,\ldots,N-1$,
\begin{equation}\label{eq:A_i-comm}
x^2i[A_l,\Box]=\sum_{j=0}^N C_{lj}A_j,\ C_{lj}\in\ePs{1}(M),
\ \sigma(C_{lj})|_{\pa\fcaldot}=0\ \text{for}\ j\neq 0.
\end{equation}
\end{lemma}

\begin{remark}
In fact, it would suffice to have $\sigma(C_{lj})|_{\pa\fcaldot}$ diagonal
with 
\begin{equation*}
\Re\sigma(C_{ll})|_{\pa\fcaldot}\geq 0,
\label{30.8.2006.36}\end{equation*}
as in \cite[Equation~(6.14)]{MR2020655} or an even weaker lower (or
upper) triangular statement, see \cite[Section~6]{MR2020655}.
\end{remark}

\begin{proof}
Let $|\tau|\,a_l=\sigma(A_l)$, considered as a homogeneous degree $1$ function
on $\Testar M$. Note that $a_N=\sigma(A_N)=x^2|\tau|^{-1}{p_0}$, so
$\sigma(\Box)=x^{-2}|\tau|^2
a_N$. As $x^2i[A_l,\Box]\in\ePs{2}(M)$ and $A_0=\Id$,
it suffices to prove that for suitable $c_{lj}$, homogeneous
degree $1$ with $c_{lj}|_\fcaldot=0$,
\begin{equation}\label{eq:symbolic-comm}
\sigma(x^2i[A_l,\Box])=\sum_{j=1}^N c_{lj}|\tau|\,a_j.
\end{equation}
But $\sigma(x^2i[A_l,\Box])=-x^2 H_{p_0} (|\tau|\,a_l)=-H_{p_0}(|\tau|\,a_l)$
is homogeneous of degree $2$,
vanishing at $\fcaldot$, so \eqref{eq:symbolic-comm} follows, though
not necessarily with $c_{lj}|_\fcaldot=0$.

By Lemma~\ref{lemma:linearization} and the part of the statement of
Lemma~\ref{lemma:test-module} regarding $da_l(q)$, $d a_l$ is an
eigenvector of the linearization of $\sigma H_{p_0}$ with eigenvalue
$-\xis$, and $d\sigma$ is also an eigenvector with eigenvalue
$-\xis$. Correspondingly, using the homogeneity and working on $\Testar M$,
$\sigma H_{p_0} (\sigma^{-1}a_l)=(\sigma H_{p_0}\sigma^{-1})a_l
+\sigma^{-1}(\sigma H_{p_0})(a_l)$ vanishes at $\pa\fcaldot$, so we may
indeed take $c_{lj}|_{\pa\fcaldot}=0$.
\end{proof}

The final piece of information we needed to analyze propagation
of regularity with respect to $\module$ is to recall from
Lemma~\ref{lemma:linearization} that $x$, resp.\ $\sigma$ are eigenvalues
of the linearization of $\sigma H_p$ of eigenvalue $\xis$, resp.\ $-\xis$,
and correspondingly for $W_{m,l}\in\ePs{m,l}(M)$
with symbol $w_{m,l}=|\tau|^m x^l
=\sigma^{-m} x^l$,
\begin{equation}\label{HMV.103}
i[W_{m,l},\Box]=W_{m-1,l-1}C_0, \\
\text{ where } C_0\in\ePs{0,0}(M),\ \sigma(C_0)|_{\pa\fcaldot}=-(m+l)\xis.
\end{equation}

\section{Domains}

An essential ingredient in our use of the edge calculus will be the
identification of domains of powers of the Laplacian with weighted edge
Sobolev spaces.

\begin{definition}
Let $\dom$ denote the Friedrichs form domain of $\Lap,$
i.e.\ the closure of $\dCI(X)$ with respect to the norm
$(\norm{du}_{L^2_g}^2 + \norm{u}_{L^2_g}^2)^{1/2}.$  Let $\dom_s$ denote the
corresponding domain of $\Lap^{s/2}$ (hence $\dom=\dom_1$),
and $\dom'=\dom_{-1}$.
\end{definition}
Recall that $L^2_g(X)= x^{-(f+1)/2} L^2_b(X),$ and that we use the b weight
in our definition of edge Sobolev spaces.

We remark that multiplication by elements of $\CI_Y(X)$ (the subspace
of $\CI(X)$ consisting of fiber constant functions at $\pa X$) preserves
$\dom$, as for $$\phi\in\CI_Y(X),\ V\in x^{-1}\Ve(X),$$ we have $$V\phi u=\phi Vu
+[V,\phi]u,$$ and  $$[V,\phi]\in\CI(X)$$ is bounded on $L^2_g$. Thus,
$\dom$ can be characterized locally away from $\pa X$, plus locally
in $Y$ near $\pa X$ (i.e.\ near $\pa X$ the domain does not have a local
characterization, but it is local in the base $Y$, so the non-locality is in
the fiber $Z$).

The basic lemma is the following.

\begin{lemma}\label{lemma:poincare}
We have, for $v \in \dCI(X),$ $f>1$,
\begin{equation}\label{normcontrol}
\norm{x^{-1} v}^2 + \norm{x^{-1}D_z v}^2+\norm{D_x v}^2+\norm{D_y v}^2
\leq C \norm{v}_{\dom}^2.
\end{equation}
\end{lemma}

\begin{remark}
This is the only result we need for the diffractive theorem, which uses the
b-calculus. If $f=1$, we use more carefully crafted ps.d.o's to make sure
that the control of $\|v\|$, rather than $\|x^{-1}v\|$, suffices---indeed,
$\|x^{-1}v\|$ is {\em not} controlled by $\|v\|_{\dom}$ if $f=1$.  A
fortiori, in the extreme case $f=0,$ we do not control $\norm{x^{-1} v}$ by
the domain norm.

That a reverse inequality to \eqref{normcontrol} holds too is immediate
from the definition of $\dom$ and the form of the metric.
\end{remark}

\begin{proof}
As $\norm{v}_{\dom}^2=\norm{dv}_{L^2_g}^2 + \norm{v}_{L^2_g}^2$, the
form of the metric $g$ implies
\begin{equation*}
\norm{x^{-1}D_z v}^2+\norm{D_x v}^2+\norm{D_y v}^2
\leq C \norm{v}_{\dom}^2.
\end{equation*}
So it remains to show that $\norm{x^{-1} v}^2
\leq C \norm{v}_{\dom}^2$ for $v$ supported near $\pa X$,
which in turn will follow from
$\norm{x^{-1} v}^2\leq C(\|D_x v\|^2+\|v\|^2)$. Treating $y,z$ as parameters,
this in turn is an immediate consequence of the standard analogous
one-dimensional result (with respect to the measure $x^f\,dx$).
\end{proof}

As an immediate consequence, we deduce the following:

\begin{proposition}\label{prop:domains}
Suppose that either $s \in [0,1),$ or $s=1$ and $f>1$. Then
$\dom_s = x^{-(f+1)/2+s} \eH{s}(X).$
\end{proposition}

\begin{proof}
For $s=1$, $f>1$,
this follows immediately from Lemma~\ref{lemma:poincare} and the
definition of the Friedrichs quadratic form domain $\dom$ as the closure
of $\dCI(X)$ with respect to
$\|\cdot\|_{\dom}.$ Indeed, the lemma gives $\|u\|_{x^{-(f+1)/2+1}\eH{1}(X)}
\leq C\|u\|_{\dom}$ for $u\in\dCI(X)$;
the analogous reverse inequality is immediate, so the density of
$\dCI(X)$ in both spaces proves the proposition for $s=1$. As $\dom_\theta$,
$\theta\in[0,1]$, resp.\ $x^{\theta-(f+1)/2}\eH{\theta}(X)$, are the
complex interpolation spaces for $(L^2,\dom)$, resp.\ $(L^2,x^{1-(f+1)/2}
\eH{1}(X))$, the proposition follows for $s\in[0,1]$, $f>1$.

A different way of identifying $\dom,$ which indeed works for arbitrary $f,$ is
to note that working in a local coordinate chart $O'$ in $Y$ over which
the fibration of $\pa X$ is trivial, and extending the fibration
to a neighborhood of $\pa X$, one has an open set $O$ in $X$
which one can identify with $[0,\ep)_x\times O'\times Z$. As $\dom$ is
a $\CI_Y(X)$ module, we only need to characterize the elements of
$\dom$ supported in compact subsets of $O$ which have a product decomposition.
But directly from the definition of $\dom$, such functions are exactly the
correspondingly supported elements
of
\begin{equation*}
H^1(Y;L^2([0,\ep)_x\times Z;x^f\,dx\,dz))\cap L^2(Y;D([0,\ep)\times Z)),
\end{equation*}
where $D$ stands for the Friedrichs form domain of any conic Laplacian on
$[0,\ep)\times Z$. These domains $D$ are well-known, and are described in
\cite[Equation~(3.11)]{Melrose-Wunsch1}. In particular, for $f>1$,
$D([0,\ep)\times Z)=x^{1-(f+1)/2}\bH{1}([0,\ep)\times Z)$, providing an
alternative method of characterization of the domain $\dom$ for $f>1$.

Now, $x^{1-(f+1)/2}\eH{1}(X)\subset\dom$ always (with a continuous inclusion)
as both are completions of $\dCI(X)$, with the first norm being stronger
than the second. Complex interpolation with $L^2$ (of which both are
a subspace) gives the inclusion $x^{\theta-(f+1)/2}\eH{\theta}(X)\subset
\dom_\theta$ for $\theta\in[0,1]$.

To see the reverse inclusion, we proceed as follows. In order to
simplify the notation, we localize to a product neighborhood $O$
of $\pa X$, but do not denote this explicitly.
Since complex interpolation of the spaces $L^2=\dom_0$ and $\dom=\dom_1$
yields $\dom_\theta$ for $0\leq\theta\leq 1$, with an analogous
statement for $D([0,\ep)\times Z)$, while complex interpolation of $L^2$ and
$H^1$ yields $H^\theta$ for $0\leq\theta\leq 1$, we deduce that
\begin{equation*}
\dom_\theta\subset H^\theta(Y;L^2([0,\ep)_x\times Z))\cap L^2(Y;D_\theta
([0,\ep)\times Z)),
\end{equation*}
meaning elements of $\dom_\theta$ multiplied by a cutoff supported in $O$
identically $1$ near $\pa X$ are in the right hand side.
For $\theta\in[0,1)$, $D_\theta([0,\ep)\times Z)=x^{\theta-(f+1)/2}\bH{\theta}
([0,\ep)\times Z)$, even if $f=1$, so we deduce that for $\theta\in[0,1)$,
\begin{equation*}
\dom_\theta\subset H^\theta(Y;L^2([0,\ep)_x\times Z))\cap L^2(Y;x^{\theta-(f+1)/2}
\bH{\theta}([0,\ep)\times Z)).
\end{equation*}

The right side is exactly $x^{\theta-(f+1)/2}\eH{\theta}(O).$ In fact,
this follows from
\begin{multline*}
H^1(Y;L^2([0,\ep)_x\times Z)),\ L^2(Y;x^{1-(f+1)/2}
\bH{1}([0,\ep)\times Z))\\ \Mand\ x^{1-(f+1)/2}\eH{1}(O)
\end{multline*}
being the form domains of commuting self-adjoint
operators $A_1,$ $A_2$, resp.\ $A=A_1+A_2,$
with $A_j\geq \Id,$ as this implies the non-trivial (and relevant)
direction, i\@.e\@. that
\begin{equation*}
H^\theta(Y;L^2([0,\ep)_x\times Z))\cap L^2(Y;x^{\theta-(f+1)/2}
\bH{\theta}([0,\ep)\times Z))\subset x^{\theta-(f+1)/2}\eH{\theta}(O)
\end{equation*}
by using $A^{(\theta-z)/2}$ to produce the holomorphic family used in
the definition of complex interpolation. For instance, one may take
$A_1=\Delta_Y$ for some metric on $Y$, and $A_2=\Delta_{[0,\ep)\times Z}
+cx^{-2}$ for a conic metric and $c>0$ sufficiently large.

Thus, the proposition follows for $s\in [0,1)$ and $f$ arbitrary (as well as
$s=1$, $f>1$). 
\end{proof}

We write $\domt$, etc.\ for the analogous
spaces on $M$:

\begin{definition}
$\|u\|_{\domt}^2=\|D_t
  u\|_{L^2(M)}^2+\|d_Xu\|^2_{L^2(M)}+\|u\|^2_{L^2(M)}$.
We also write $\domt([a,b])$ for the space with the same norm on
  $[a,b]\times X.$
\end{definition}

We will require the following result about the interaction of the test
module $\module$ from Definition~\ref{definition:coisotropic} with the norm
on $\domt:$

\begin{lemma}\label{lemma:adjoints}
Let $A \in \module.$  Let $A^*_\domt$ denote the adjoint of $A$ respect to
the norm on $\domt.$  Then
$$
A^*_\domt = A^*+B+R
$$
with $B \in \ePs{0,0}(M)$ and $R$ a smoothing operator mapping $\domt(M)
\to \CI(M).$
\end{lemma}
\begin{proof}
Let $\tBox = D_t^2+\Lap+1.$  If $u,v\in \dCI(M)$ have compact
support, we have
\begin{align*}
  \ang{A u,v}_\domt &= \ang{\tBox Au,v} \\
  & = \ang{u,A^*v}_\domt +\ang{[\tBox,A]u,v}.
\end{align*}
Using the usual elliptic parametrix construction, we may write
$$
[\tBox,A] = B^*\tBox+R'
$$
where $B^* \in \ePs{0,0}(M)$ and $R' \in \ePs{-\infty,-2}(M).$  By the
Riesz Lemma, $\ang{R'u,v} = \ang{Ru,v}_\domt$ for some smoothing operator
$R$, so we obtain the desired form of the adjoint.
\end{proof}

The domain $\dom$ localizes with respect to fiber-constant functions.
Thus, we make the definition:

\begin{definition}
For $u\in\CmI(X)$, we say that $u\in\dom_{\loc}$ if $\phi u\in\dom$
for all fiber constant $\phi\in\CI_c(X)$. Similarly,
for $u\in\CmI(X)$, we say that $u\in\dom'_{\loc}$ if $\phi u\in\dom'$
for all fiber constant $\phi\in\CI_c(X)$.

We define the localized domains on $M$ analogously.
\end{definition}

Note that if $u\in\dom$ then certainly $u\in\dom_{\loc}$.

\section{Wavefront sets ($f>1$)}
\emph{Throughout this section, we assume $f>1.$} For the $f=1$ case, see
Section~\ref{sec:f1}.\footnote{As always, for the case of a manifold with
boundary, i.e.\ $f=0,$ we refer the reader to \cite{MS1,MR83h:35120}.}

By convention all ps.d.o's in this section have Schwartz kernels
supported in a fixed compact set $\tilde{K} \times \tilde{K},$ and we
define associated local norms by fixing a function $\tilde \phi \in
\CIc(M)$ equal to $1$ on neighborhood of $\tilde K;$ we let
$$
\norm{u}_{\bullet, \loc} = \norm{\tilde \phi u}_{\bullet}
$$
for various choices of Sobolev space.

\begin{lemma}
If $V\in x^{-1}\Ve(M)$ and $A\in\Psibc^m(M)$ then 
\begin{equation*}
[A,V]\in x^{-1}\Psibc^m(M).
\label{30.8.2006.35}\end{equation*}
\end{lemma}
(Recall that the calculus with conormal estimates, $\Psibc^*(M),$ is
described at the end of \S\ref{section:calculi}.)

\begin{proof}
The result follows from the fact that $V\in x^{-1}\Diffb{1}(M)$.
\end{proof}

\begin{lemma}\label{lemma:Psibc-map}
Any $A\in\Psibc^0(M)$ with compact support defines a continuous linear map
$A:\domt\to\domt$ with a norm bounded by a seminorm of $A$ in
$\Psibc^0(M)$.

Moreover, for any $K\subset M$ compact, any $A\in\Psibc^0(M)$
with proper support
defines a continuous map from the subspace of $\domt$
consisting of distributions supported in $K$ to $\domt_c$.
\end{lemma}

\begin{proof}
The norm on $\domt$ is $(\|d_M u\|^2+\|u\|^2)^{1/2}$, and $\CI_c(M^\circ)$
is dense on $\domt$, so it suffices to show that for
$u\in\CI_c(M^\circ)$, $\|Au\|_{L^2}\leq C\|u\|_{L^2}$, and
$\|VAu\|_{L^2}\leq C(\|d_M u\|+\|u\|)$ for $V\in x^{-1}\Ve(M)$.
As $\Psibc^0(M)$ is bounded on $L^2$, the first claim is clear.
Moreover,
\begin{equation*}
\|VAu\|_{L^2}\leq\|AVu\|_{L^2}+\|[V,A]u\|_{L^2}\leq C\|Vu\|_{L^2}
+C\|x^{-1}u\|_{L^2}
\end{equation*}
where we used the previous lemma for the second term, and also that
$\Psibc^0(M)$ is bounded on $L^2$. Finally, as $f>1$,
$\|x^{-1}u\|_{L^2}\leq C(\|D_x u\|_{L^2}+\|u\|_{L^2}),$ finishing the
proof.
\end{proof}
We note that this lemma is in fact false without the assumption $f>1$.

We are now able to follow the treatment of \cite{Vasy5} almost verbatim,
using the preceding Lemma in place of \cite[Lemma~3.2]{Vasy5}.
\begin{definition}
For $m\geq 0$, we define $H^m_{\domt,b,c}(M)$
as the subspace of $\domt$
consisting of $u\in \domt$ with $\supp u$ compact and
$Au\in \domt$ for some (hence any, see \cite{Vasy5})
$A\in\Psib^m(M)$ (with compact support) which is elliptic over $\supp u$,
i.e.\ $A$ such that such that
$\sigma_{{\bo},m}(A)(q)\neq 0$ for any
$q\in\Tbstar_{\supp u}M\setminus 0$.

We let $H^m_{\domt,b,\loc}(M)$ be the subspace of $\domt_\loc$
consisting of $u\in \domt_{\loc}$ such that for any $\phi\in\CI_c(M)$,
$\phi u\in H^{m}_{\domt,b,c}(M)$.
\end{definition}

Thus, $H^m_{\domt,b,\loc}(M)$ consists of distributions conormal to finite
order relative to $\domt$.
We also define the spaces with a negative order of conormal regularity
relative to $\domt$:

\begin{definition}\label{Def:Hm-neg}
Let $m<0$, and $A\in\Psib^{-m}(M)$ be elliptic on $\Sbstar M$
with proper support.
We let $H^{m}_{\domt,\bo,c}(M)$ be the space of all $u\in\CmI(M)$ of the
form $u=u_1+Au_2$ with $u_1,u_2\in \domt_c$.
We let
\begin{equation*}
\|u\|_{H^{m}_{\domt,\bo,c}(M)}
=\inf\{\|u_1\|_{\domt}+\|u_2\|_{\domt}:\ u=u_1+Au_2\}.
\end{equation*}

We also let $H^{m}_{\domt,\bo,\loc}(M)$ be the space of all $u\in\CmI(M)$
such that $\phi u\in H^{m}_{\domt,\bo,c}(M)$ for all $\phi\in\CI_c(M)$.
\end{definition}

The spaces $H^m_{\domt,\bo,c}$ and $H^m_{\domt,\bo,\loc}$ do not
depend on the choice of $A$; see \cite{Vasy5}.

With this definition, $A\in\Psibc^k(M)$
defines a map
\begin{equation*}
A:H^m_{\domt,b,c}(M)\to H^{m-k}_{\domt,b,c}(M);
\end{equation*}
if the supports are kept in a fixed compact set, this is a continuous
linear map between Banach spaces. We now proceed to microlocalize these
spaces.

\begin{definition}\label{Def:dWF}
Let $u\in H^s_{\domt,\bo,\loc}$ for some $s$ and let $m\geq 0$.
We say that $q\in\Tbstar M\setminus 0$ is not in $\dWF^{m}(u)$ if
there exists $A\in\Psib^m(M)$ such that $\sigma_{{\bo},m}(A)(q)\neq 0$
and $Au\in \domt$.

For $m=\infty$, we say that $q\in\Tbstar M\setminus 0$ is not in $\dWF^{m}(u)$
if there exists $A\in\Psib^0(M)$ such that $\sigma_{{\bo},0}(A)(q)\neq 0$
and $LAu\in \domt$ for all $L\in\Diffb{}(M)$, i.e.\ if
$Au\in H^{\infty}_{\domt,b,\loc}(M)$.

We also use the notation $\dWF(u)=\dWF^\infty(u).$
\end{definition}

\begin{remark}
Note that $\dWF(u)$ is only defined for $u\in\bigcup_s
H^s_{\domt,\bo,\loc}(M)$.  It is reasonable thus to call $\bigcup_s
H^s_{\domt,\bo,\loc}(M)$ the space of $\domt$-admissible distributions, and
to refer to solutions of $\Box u=0$ with $u\in\bigcup_s
H^s_{\domt,\bo,\loc}(M)$ as admissible solutions---even, in a more
general setting, if $\Box$ is the `Laplacian' of a general
pseudo-Riemannian metric.
\end{remark}

We recall the definition of the operator wave front set of a family
of ps.d.o's.

\begin{definition}
Suppose that $\cB$ is a bounded subset of $\Psibc^k(M)$, and $q\in\Sbstar M$.
We say that $q\notin\bWF'(\cB)$ if there is some $A\in\Psib(M)$
which is elliptic at $q$ such that $\{AB:\ B\in\cB\}$ is a bounded
subset of $\Psib^{-\infty}(M)$.
\end{definition}

Note that the wave front set of a family $\cB$
is only defined for bounded families. It can
be described directly in terms of quantization of (full) symbols,
much like the operator wave front set of a single operator. All
standard properties of the operator wave front set also hold for a family;
e.g.\ if $E\in\Psib(M)$ with $\bWF'(E)\cap\bWF'(\cB)=\emptyset$ then
$\{BE:\ B\in\cB\}$ is bounded in $\Psib^{-\infty}(M)$.

The key lemma here concerns families $\cB$ with wave front set disjoint
from a wave front set $\dWF^k(u)$ of a distribution $u$; it is the
quantitative version of the microlocality of $\Psibc(M)$:
\begin{equation*}
\dWF^{m-k}(Bu)\subset\bWF'(B)\cap\dWF^m(u),\ B\in\Psibc^k(M),
\end{equation*}
see \cite[Lemma~3.9]{Vasy5}.

\begin{lemma}\label{lemma:WFb-mic-q}
Suppose that $K\subset\Sbstar M$ is compact, and $U$ a neighborhood of $K$
in $\Sbstar M$. 
Let $Q\in\Psib^k(M)$ be elliptic on $K$ with $\bWF'(Q)
\subset U$.
Let $\cB$ be a bounded subset of $\Psibc^k(M)$
with $\bWF'(\cB)\subset K$.
Then for any $s\in\RR$ there is a constant $C>0$ such that
\begin{equation*}
\|Bu\|_{\domt}\leq C(\|u\|_{H^s_{\domt,\bo,\loc}}+\|Qu\|_{\domt}).
\end{equation*}
 for $B\in\cB$, $u\in H^s_{\domt,\bo,\loc}$ with $\dWF^{k}(u)\cap
U=\emptyset$,
\end{lemma}
For the proof, see \cite[Lemma~3.13 and
Lemma~3.18]{Vasy5}.

\begin{remark}
Instead of working with $\domt$, one can also work with $\domt'$, as
$\Psibc^0(M)$ acts on it by duality. All the above results have
their analogues with $\dWF$ replaced by $\mWF$.
\end{remark}

We will also employ a wave front set on a quotient space of $\Sbstar M.$  In
\S\ref{section:edgeb} we will define a space $\Sbdotstar M$ with a proper,
surjective map $p:\Sbstar M\to\Sbdotstar M;$ $\Sbdotstar M$ will be
equipped with the quotient topology. In
this case, by abusing the notation slightly, we use the following terminology:
\begin{definition}\label{def:WF-Sbdotstar}
$\dWF(u)$ and $\mWF(u)$, considered as subsets of $\Sbdotstar M$,
are the images under the quotient map $p$ of $\dWF(u)$ and $\mWF(u)$,
considered as subsets of $\Sbstar M$.
\end{definition}

Thus, for $q\in\Sbdotstar M$, $q\notin\dWF(u)$ if and only if
$p^{-1}(q)\cap \dWF(u)=\emptyset,$ i.e.\ $u$ is well-behaved at all
pre-images of $q.$

\section{Edge/b relationship}
\label{section:edgeb}

Let $N$ denote an edge manifold with fibration
$$
Z \to \pa N \overset{\pi_0}{\to} W.
$$
Canonical coordinates on $\Testar N$ induced by coordinates $(x,w,z)$
(with $w$ pull-backs of coordinates on $W$)
on $N$ are $(x,w,z,\xi,\mu,\zeta)$ corresponding to writing covectors
in $\Testar N$ as
\begin{equation*}
\xi \frac{dx}x + \mu \cdot \frac{dw}x + \zeta \cdot dz.
\end{equation*}
Analogously, canonical coordinates on $\Tbstar N$ are
$(x,w,z,\xib,\mub,\zetab)$
corresponding to writing covectors
in $\Tbstar N$ as
\begin{equation}\label{bcoords}
\xib \frac{dx}x + \mub \cdot dw + \zetab \cdot dz.
\end{equation}

Let $\pi$ denote the bundle map $\Testar N \to \Tbstar N$ given in
canonical coordinates by
$$
\pi(x,w,z,\xi,\mu,\zeta)=(x,w,z,x\xi,\mu,x\zeta).
$$ This map can be obtained more invariantly by identifying $\Testar N$
with the bundle $\Tmstar N:=x \Testar N,$ i.e.\ the bundle whose sections
can be written as $x$ times an edge one-form; the identification is via
multiplication by $x$, naturally.  The map $\pi$ is then the inclusion map
of this bundle into $\Tbstar M.$ We define the compressed cotangent bundle
by setting
\begin{align*}
\Tbdotstar N &= \pi(\Testar N)/Z,\\ \dot\pi:\Testar N &\to\Tbdotstar N
\end{align*}
the projection, where, here and henceforth, we take the quotient by $Z$ to
act only over the boundary (where it is defined) and the topology to be
given by the quotient topology.  We remark that $\Tbdotstar_{\pa N} N$ can
be identified with $T^* W.$

In our geometric setting, we have $N=M:=\RR\times X,$ $W=\RR\times Y,$
and the $w$ variables
above are replaced by $w=(t,y),$ and correspondingly $\mu=(\tau,\eta).$

Below it is much more convenient to work with the cosphere bundles
rather than the cotangent bundles. These are defined
as
\begin{equation*}
\Sbstar M=(\Tbstar M\setminus 0)/\RR^+,\ \Sbdotstar M=\Sbstar M/Z,
\ \Sestar M=(\Testar M\setminus 0)/\RR^+,
\end{equation*}
where $0$ denotes the zero sections of the respective bundles and $\RR^+$
acts by scalar multiplication in the fibers. The maps $\pi$ and
$\dot\pi$
commute with the $\RR^+$-action on the cotangent bundles, but they
do map some non-zero covectors to the zero section of $\Tbstar M$
(at $x=0$),
so they do not descend to maps $\Sestar M\to\Sbdotstar M$.
(We remark that $\Sestar M$ and $\Smstar M$ are
naturally identified, which is one of the reasons for working with
cosphere bundles.) Nonetheless, we introduce the notation
\begin{equation*}\begin{split}
&\pi(\Sestar M)=(\pi(\Testar M)\setminus 0)/\RR^+\subset\Sbstar M,\\
&\dot\pi(\Sestar M)=(\dot\pi(\Testar M)\setminus 0)/\RR^+\subset\Sbdotstar M
\end{split}\end{equation*}
for convenience.

However, letting $p=\sigma_e(x^2\Box)\in\CI(\Testar M\setminus 0),$ we have
$$
p|_{x=0} = \tau^2-\left({\xi^2} + \abs{\eta}^2_h +
 \abs{\zeta}^2_k\right).
$$
Since
$$
\pi(0,t,y,z,\xi,\tau,\eta,\zeta)=(0,w,z,0,\tau,\eta,0),
$$
on the characteristic set
$p^{-1}(\{0\})\subset\Testar M$ of $\Box$,
where $\tau\neq 0$, non-zero covectors are mapped to non-zero covectors
by $\pi$, hence also by $\dot\pi$. (In other words, $\pa M$
is {\em non-characteristic}.)
Thus $\pi,\dot\pi$ also define maps, denoted with the same letter:
$$
\pi:\Sigma\to\Sbstar M,\ \dot\pi:\Sigma\to\Sbdotstar M,
\ \Sigma=p^{-1}(\{0\})/\RR^+.
$$ We set $$\dot\Sigma = \dot\pi(\Sigma);$$ this is the {\em compressed
characteristic set}.

In order to obtain coordinates on a subset of the cosphere bundle from the
canonical coordinates on a cotangent bundle we need to fix a homogeneous
degree $1$ function that does not vanish on the cone corresponding to this
subset.  In our case, on $\Sigma$, $\abs\tau$ can be taken as a canonical
choice.  We thus let
$$
\xis=\frac{\xi}{\abs\tau}, \etas=\frac{\eta}{\abs\tau}, \zetas=\frac{\zeta}{\abs\tau}
$$
(so $\taus\equiv 1$), and likewise, in the `b-coordinates' of \eqref{bcoords},
$$
\xibs=\frac{\xib}{\abs\taub}, \etabs=\frac{\etab}{\abs\taub}, \zetabs=\frac{\zetab}{\abs\taub}.
$$
In these coordinates, at $x=0$, $\Sigma$ is given by $\xis^2+|\etas|_h^2+
|\zetas|^2_k=1$, while at $x=0$, $\dot\Sigma$ is given by
$\xibs=0$, $\zetabs=0$, $|\etabs|_h^2\leq 1$ (with the subscripts $h,k$ as
usual denoting lengths with respect to the indicated metrics).

\begin{lemma}\label{lemma:dot-Sigma-nbhd}
Every neighborhood of a point $q_0=(t_0,y_0,\tau=1,\etabs_0)\in\dot\Sigma$
in $\dot\Sigma$
contains an open set of the form
\begin{equation*}
\{q:\ |x(q)|^2+|y(q)-y_0|^2+|t(q)-t_0|^2+|\etabs(q)
-\etabs_0|^2<\delta\}.
\end{equation*}
\end{lemma}

\begin{proof}
Let $\rho:\Sbstar M\to\Sbdotstar M$ be the quotient map.  A neighborhood
$U$ of $q_0$ in $\dot\Sigma$ is, by definition, of the form
$U'\cap\dot\Sigma$ where $U'\subset\Sbdotstar M$ is open. But, as $U'$
being open, $Z$ is compact, and $q_0\in U'$, $U'$ contains an open set of
the form
\begin{equation}\begin{split}\label{eq:Sbdotstar-nbhd}
\{q\in\Sbdotstar M:\ |x(q)|^2&+|y(q)-y_0|^2+|t(q)-t_0|^2\\
&+|\xibs(q)|^2+|\etabs(q)
-\etabs_0|^2+|\zetabs(q)|^2<\delta'\}
\end{split}\end{equation}
for some $\delta'>0$.
But $\dot\Sigma=\pi(\Sigma)$, so on $\dot\Sigma$,
$\xibs(q)=x(q')\xis(q')$ with $q'\in\Sigma$, $q=\dot\pi(q')$,
so $|\xibs(q)|\leq x(q) |\xis(q')|\leq x(q)$, hence on a set in $\dot\Sigma$
with $x(q)^2<\delta''$, we have in fact
$|\xibs(q)|^2<\delta''$, and similarly,
$|\zetabs(q)|^2<C\delta''$. Thus, the set
\begin{equation*}
\{q\in\dot\Sigma:
\ |x(q)|^2+|y(q)-y_0|^2+|t(q)-t_0|^2+|\etabs(q)
-\etabs_0|^2<(C+2)^{-1}\delta'\}
\end{equation*}
is indeed contained in \eqref{eq:Sbdotstar-nbhd}, proving the lemma.
\end{proof}

\begin{lemma}
$\dot\Sigma$ is a metrizable space.
\end{lemma}

\begin{proof}
This follows from $\Sbdotstar M$ being metrizable, which can be shown
easily directly. We provide an alternative argument
\mbox{}From Lemma~\ref{lemma:dot-Sigma-nbhd}
it is immediate that $\dot\Sigma$ is regular,
i.e.\ if $q_0\in\dot\Sigma$, $F\subset\dot\Sigma$ closed, $q_0\notin F$,
then there
are disjoint neighborhoods of $U$ of $q_0$ and $U'$ of $F$: indeed, we
can use sub- and superlevel sets of the continuous function
$f_{q_0}(q)=|x(q)|^2+|y(q)-y_0|^2+|t(q)-t_0|^2+|\etabs(q)
-\etabs_0|^2$. If now $K$ is compact, $F$ is closed, for each $q_0\in K$
there is a function $\chi_{q_0}\in\CI_c(\RR)$, identically $1$ near $0$ such
that $\supp(\chi_{q_0}\circ f_{q_0})\cap F=\emptyset$.
Let $V_{q_0}=\{q:\ \chi_{q_0}\circ f_{q_0}>1/2\}$, then
$\{V_{q_0}:\ q_0\in K\}$ is an open cover of $K$. Let $\{V_{q_1},\ldots,
V_{q_k}\}$ be a finite subcover. Then $F(q)=\sum_{j=1}^k
(\chi_{q_j}\circ f_{q_j})(q)$ satisfies $F>1/2$ on $K$, $F=0$ on $F$,
so again sub- and superlevel sets of $F$ can be used as disjoint
neighborhoods of $F$ and $K$. If $K$ is closed but not compact,
the sum can still be made locally finite, proving that $\dot\Sigma$
is normal. This proves the lemma.
\end{proof}

We may define the following three subsets of $\pi(\Sestar[{\pa M}]M)$
(`elliptic,' `glancing' and `hyperbolic'):
\begin{equation*}\begin{split}
&\cE = \pi(\Sestar[{\pa M}]M)\setminus \pi(\Sigma),\\
&\cG = \{q \in \pi(\Sestar[{\pa M}]M): \Card(\pi^{-1} (q)\cap\Sigma) =1\},\\
&\hyp = \{q \in \pi(\Sestar[{\pa M}]M): \Card(\pi^{-1} (q)\cap \Sigma) \geq 2\}.
\end{split}\end{equation*}
In coordinates, we have
$$
\pi(0,t,y,z,\xi,\taus,\etas,\zeta)=(0,t,y,z,0,
\taus,\etas,0),
$$ hence the three sets are defined by $\{\taubs^2<\abs{\etabs}^2\},$
$\{\taubs^2=\abs{\etabs}^2\},$ $\{\taubs^2>\abs{\etabs}^2\}$ respectively
inside $\pi(\Sestar[{\pa M}]M)$, which is given by $x=0$, $\zetabs=0$,
$\xibs=0$.

We may also define the corresponding sets in $\Sbdotstar_{\pa M} M$
 (hence quotiented by $Z$ and denoted with a dot):
\begin{equation*}\begin{split}
&{\dcE} = \cE/Z,\\
&{\dcG} = \cG/Z,\\
&{\dhyp} = \hyp/Z.
\end{split}\end{equation*}

\begin{remark}
The hyperbolic and glancing sets can be visualized as (co-) vectors at the
boundary that are respectively transverse to the geometric boundary $Y$ and
tangent to it.  A major difference from the situation of manifolds with
boundary is that given a point in $\pa X$ and tangential momenta
$\taubs,\etabs$ with $\taubs^2>\abs{\etabs}^2$ we have specified not just a
single point in the hyperbolic set but a whole manifold of them
parametrized by $Z.$ If $X=[\RR^3;\gamma]$ is given by the blowup of a
smooth curve $\gamma,$ then we may visualize these points as the tangents to the
cone of rays emanating from a single point in $\gamma$ making a fixed angle
to $\gamma$ (hence striking the blown-up space at all possible points in a
fixed fiber).  It is thus the whole cone that is being considered at once when
we specify a point in the quotient space $\dcE.$
\end{remark}

Bicharacteristics are usually defined as curves in the cotangent bundle.
It is in fact more natural to define them in the cosphere bundle.
However, the Hamilton vector field of a homogeneous degree $m$ function
is homogenous of degree $m-1$, so the curves in the cosphere bundle are
only defined up to reparametrization unless $m=1$. To fix a
parameterization, we renormalize $\sigma_e(\Box)$
to make it homogeneous of degree $1,$
by considering
$$
\hat p=x\abs{\tau}^{-1}\sigma_e(\Box)\in x^{-1}\CI(\Te^*M\setminus o).
$$
Note that $\hat p$ is actually independent of the choice of $x$; $x^{-1}\tau
=\sigma_e(D_t)$, which is constant along bicharacteristics in $T^*M^\circ$,
so the bicharacteristics do not `slow down' as they approach $\pa M$,
unlike the integral curves of $x^2H_p$. In particular, their projection to $X$
are constant speed geodesics.
Then $H_{\hat p}$ is homogeneous of
degree $0$, hence defines a vector field on $\Sestar[M^\circ] M$.
Equivalently, on $\Sigma$, we can renormalize $H_{\sigma_e(\Box)}$:
$$
H_{\hat p}=\abs{\tau}^{-1}xH_{\sigma_e(\Box)}=
-\frac{2}{\abs{\tau}x}H,\ H=-2^{-1}x^2H_{\sigma_e(\Box)},
$$
$H$ as in \eqref{eq:H-form}.
Explicitly,
\begin{equation}\begin{split}\label{eq:H_p-m}
&-\frac{1}{2}
H_{\hat p}=\\
&-\taus\,\pa_t+\xis\,\pa_x+x^{-1}(|\zetas|^2+O(x^2)\zetas^2
+O(x^2)\etas\zetas+O(x)\etas^2)\pa_{\xis}\\
&+ (x^{-1}\zetas_i \kbar^{ij} + O(1)\etas) \pa_{z_j} +
(-\h x^{-1} \zetas_i\zetas_j \frac{\pa \kbar^{ij}}{\pa z_k} + O(1)\etas \zetas+
O(x) \etas^2)\pa_{\zetas_k} \\
&+(\etas_j H^{ij} + O(x^2) \etas +
O(x)\zetas) \pa_{y_i} + (O(1)\etas^2 + O (x) \etas
\zetas + O(1)\zetas^2)\pa_{\etas_i}.
\end{split}\end{equation}
While $H_{\hat p}$ is not a $\CI$ vector field (it is $xH_{\hat p}$ that is
$\CI$) on $\Sestar M$, it acts on $\CI$ functions independent of
$z,\xis,\zetas$ at $x=0$, hence on pull-backs of functions from $\Sbdotstar
M$.

Note that by definition
$$
\dot\Sigma\cap{{}^b\overset{\!\! .}{T_{\pa M}^*}}M=\dcG\cup\dhyp.
$$

\begin{definition}\label{def:gbb}
Let $I$ be an interval.  We say that a continuous map
$\gamma:I\to\dot\Sigma$ is a generalized broken bicharacteristic if for all
$f\in\cC(\Sbdotstar M)$ real valued with $\dot\pi^*f\in\CI (\Sestar M)$,
\begin{equation*}\begin{split}
&\liminf_{s\to s_0}\frac{(f\circ\gamma)(s)-(f\circ\gamma)(s_0)}{s-s_0}\\
&\ \geq \inf\{H_{\hat p}(\dot\pi^* f)(q):
\ q\in\dot\pi^{-1}(\gamma(s_0))\cap\Sigma\}
\end{split}\end{equation*}
holds.
\end{definition}

\begin{remark}
This definition may appear slightly complicated, especially as compared to
the more explicit statement of the next lemma, which could be used as an
alternative definition. However, it {\em is} the natural definition, for
certainly any putative notion of generalized broken bicharacteristic should
have the stated property. While it may be possible to strengthen the
definition so as to rule out certain tangential rays in some cases (as
happens in the smooth boundary setting), this would be a delicate matter.
In $N$-body scattering one can also define generalized broken
bicharacteristics as in Definition~\ref{def:gbb}, but these do {\em not}
have a simple characterization analogous to that given below.
\end{remark}

Since $\Sestar M\ni
q\mapsto H_{\hat p}(\dot\pi^* f)(q)$ is continuous for $f$ as
in the definition, and as $\dot\pi^{-1}(\gamma(s_0))\cap\Sigma$ is
compact, the infimum on the right hand side is finite. Applying
the same estimate to $-f$, we deduce that $f\circ\gamma$ is locally Lipschitz.

\begin{lemma}
Suppose $\gamma$ is a generalized broken bicharacteristic.

\begin{enumerate}
\item
If $\gamma(s_0)\in\dcG$ then for all $f\in\cC(\Sbdotstar M)$
real valued with $\dot\pi^*f\in\CI
(\Sestar M)$, $f\circ\gamma$ is differentiable at $s_0$ with
$$
(f\circ\gamma)'(s_0)=H_p(\dot\pi^*f)(q),
$$
where $q$ is the unique point in $\dot\pi^{-1}(\gamma(s_0))\cap\Sigma$.
\item
If $\gamma(s_0)\in\dhyp$ then there exists $\ep>0$ such that
$|s-s_0|<\ep$, $s\neq s_0$ implies $x(\gamma(s))\neq 0$, i.e.\ $\gamma(s)$
does not lie over the boundary.
\end{enumerate}
\end{lemma}

\begin{proof}
Part (1) follows from the definition by applying it to both $f$ and $-f$,
and noting that the set over which the $\inf$ is taken has a single
element, so in fact
$\lim_{s\to s_0}({s-s_0})^{-1}\left({(f\circ\gamma)(s)-(f\circ\gamma)(s_0)}\right)$
exists.

For part (2), consider $f=\xibs\in\cC(\Sbdotstar M)$, so $\dot\pi^*f
\in\CI(\Sestar M)$, and indeed $\dot\pi^*f=x\xis$. For $q\in\Sigma$
with $x(q)=0$, \eqref{eq:H_p-m} gives
$$
H_{\hat p} \dot\pi^*f(q)=\xis^2(q)+|\zetas|_k^2(q)=\taus^2(q)-|\etas|_h^2(q)>0,
$$
so $\dot\pi^*f$ is strictly increasing,
and consequently non-zero in a punctured
neighborhood of $s_0$, so the same holds for $f$. But on $\dot\Sigma$,
at $x=0$, $f=0$, so we deduce that $x\neq 0$ in a punctured neighborhood
of $s_0$ as claimed.
\end{proof}

As an immediate corollary, we deduce that near points in $\dhyp$,
generalized broken bicharacteristics consist of two bicharacteristic
segments of $x^2H_p$ in $\Sestar M$ (projected to $\Sbdotstar M$), one
incoming and one outgoing.

The following are useful general facts about generalized broken
bicharacteristics.

\begin{corollary}(cf.\ Lebeau, \cite[Corollaire~2]{Lebeau5})
\label{cor:Lebeau-Lipschitz}
Suppose that $K$ is a compact subset of $\dot\Sigma$. Then there is
a constant $C>0$ such that for all generalized broken bicharacteristics
$\gamma:I\to K$, and for all functions $f$
on a neighborhood of $K$ in $\dot\Sigma$ with
$\dot\pi^*f\in\CI(\Sestar M)$, one has the uniform
Lipschitz estimate
\begin{equation*}
|f\circ\gamma(s_1)-f\circ\gamma(s_2)|\leq C\|\dot\pi^*f\|_{C^2}\,|s_1-s_2|,
\ s_1,s_2\in I.
\end{equation*}
In particular, (locally) the functions $x$, $y$, $t$, $\taub$ and $\etab$
are Lipschitz on generalized broken bicharacteristics.
\end{corollary}

We also need to analyze the uniform behavior of generalized broken
bicharacteristics. Here we quote Lebeau's results.

\begin{proposition}(cf.\ Lebeau, \cite[Proposition~5]{Lebeau5})
\label{prop:Lebeau-unif-limits}
Suppose that $K$ is a compact subset of $\dot\Sigma$, $\gamma_n:[a,b]\to K$ is
a sequence of generalized broken bicharacteristics which converge uniformly
to $\gamma$. Then $\gamma$ is a generalized broken bicharacteristic.
\end{proposition}

\begin{proposition}(cf.\ Lebeau, \cite[Proposition~6]{Lebeau5})
\label{prop:Lebeau-compactness}
Suppose that $K$ is a compact subset of $\dot\Sigma$, $[a,b]\subset\RR$
and
\begin{equation}
\rcal=\{\text{generalized broken bicharacteristics}\ \gamma:[a,b]\to K\}.
\end{equation}
If $\rcal$ is not empty then it is compact in the topology of uniform
convergence.
\end{proposition}

\begin{proof}
$\rcal$ is equicontinuous, as in Lebeau's proof,
so the proposition follows from
the theorem of Ascoli-Arzel\`a and Proposition~\ref{prop:Lebeau-unif-limits}.
\end{proof}

\begin{corollary}(Lebeau, \cite[Corollaire~7]{Lebeau5})
\label{cor:Lebeau-bichar-ext}
If $\gamma:(a,b)\to\dot\Sigma$ is a generalized broken bicharacteristic then
$\gamma$ extends to $[a,b]$.
\end{corollary}

We now make the connection with the $\fcal,$ $\fcaldot$ notation introduced
above.  Given $q \in \hyp,$ we recall from \S\ref{section:edge} that there
exist unique maximally extended incoming/outgoing bicharacteristics
$\gamma_{I/O},$ where $\sgn\xi=\pm\sgn \tau,$ such that
$q=\pa(\overline{\gamma_{\bullet}});$ as above, we denote these curves
$$
\fcal_{\bullet, q}.
$$
Likewise, for $p \in \dhyp$ we let
$$
\fcaldot_{\bullet,p}=\bigcup_{q \in \pi_0^{-1} p} \fcal_{\bullet, q}.
$$
We will abuse notation slightly to write
\begin{equation}
\dhyp_{I/O}=\pa\fcaldot_{I/O}
\end{equation}
for the endpoints of incoming/outgoing hyperbfolic bicharacteristics at the
boundary.

As opposed to generalized broken bicharacteristics, we may define a
different flow relation, at least away from glancing rays, as follows:
\begin{definition}\label{def:geomflow}
Let $q,q' \in \hyp,$ with $\pi_0(q) =\pi_0(q').$  We say that
$$
\fcal_{I,q},\ \fcal_{O,q'}
$$ are related under the forward geometric flow (and vice-versa under the
backward flow) if there exists a geodesic of length $\pi$ in $Z_{y(q)}$
connecting $z(q)$ and $z(q').$ If $a \in \fcal_{I,q}$ with $q \in \hyp,$ we
let the forward flowout of $a$ be the union of the forward geodesic segment
through $a$ and all the $\fcal_{O,q'}$ that are related to $\fcal_{I,q}$
under the forward geometric flow (and vice-versa for backward flow).  (If
the forward flow through $a$ stays in $M^\circ,$ we simply let its forward
flowout be the ordinary flowout under geodesic flow.)
\end{definition}

We note that the flow relation generated by generalized broken
bicharacteristics differs from that of Definition~\ref{def:geomflow} as
follows: a ray in $\fcal_{I,q}$ can be continued as a generalized broken
bicharacteristic by \emph{any} $\fcal_{O,q'}$ with $\pi_0(q')=\pi_0(q);$
there is no requirement on the relative locations of $q,q'$ in the fiber.

\section{Ellipticity}\label{sec:elliptic}

First note that products of elements of $x^{-1}\Diffe{1}(M)$ and
$\Psib^m(M)$ can be written with the products taken in either order:

\begin{lemma}\label{lemma:rearrange}
For any $Q\in x^{-1}\Diffe{1}(M)$ and $A\in\Psib^m(X)$ there exist
$A'\in\Psib^{m}(X)$, $Q'\in x^{-1}\Diffe{1}(M)$ such that 
\begin{equation}
QA=AQ+A'Q'.
\label{30.8.2006.17}\end{equation}
\end{lemma}

\begin{proof}
Since $x^{-1}\Diffe{1}(M)\subset x^{-1}\Diffb{1}(M),$ $[Q,A]\in
x^{-1}\Psib^m(M),$ so \eqref{30.8.2006.17} holds with
$A'=[Q,A]x\in\Psib^m(M)$ and $Q'=x^{-1}\in x^{-1}\Diffe{1}(M).$
\end{proof}

In fact the span of products on the left in \eqref{30.8.2006.17}, and hence
on the right, may be characterized as the subspace of
$x^{-1}\Psib^{m+1}(M)$ fixed by the vanishing of the principal symbol at
$x=0$, $\sigma=0,$ $\zeta=0.$ Similarly one may consider higher order
products. 

\begin{definition} Let 
\begin{equation}
x^{-k}\Diffe{k}\Psib^{m}(M)\subset x^{-k}\Psib^{k+m}(M)
\label{30.8.2006.18}\end{equation}
be the span of the products $QA$ with $Q\in x^{-k}\Diffe{k}(M)$ and
$A\in\Psib^{m}(M).$
\end{definition}

Induction based on Lemma~\ref{lemma:rearrange} gives

\begin{corollary}\label{cor:rearrange}
The space $x^{-k}\Diffe{k}\Psib^{m}(M)$ is also the span of the products
$AQ$ with $Q\in x^{-k}\Diffe{k}(M)$ and $A\in\Psib^{m}(M)$ and hence  
\begin{equation*}
\bigcup_{k,m}x^{-k}\Diffe{k}\Psib^m(M)
\label{30.8.2006.19}\end{equation*}
forms a bigraded ring which is closed under adjoints with respect to any
b-density.
\end{corollary}

Next consider the Hamilton vector field of the principal symbol of an element
of $\Psib^m(M).$ 

\begin{lemma}\label{lemma:b-Ham-vf} If $A\in\Psib^m(M)$ then the Hamilton
vector field $H_a$ of $a=\sigma_{{\bo},m}(A),$ defined initially on
$T^*M^\circ,$ extends to an element of $\Vb(\Tbstar M)$ and in canonical
local coordinates $(x,w,z,\xib,\mub,\zetab)$ on $\Tbstar M,$ induced by
local coordinates $(x,w,z)$ on $M,$
\begin{equation*}
H_a=(\pa_\xib a) x\pa_x+(\pa_{\mub} a) \pa_{w}+(\pa_{\zetab} a) \pa_{z}
-(x\pa_x a)\pa_\xib-(\pa_{w} a)\pa_{\mub}-(\pa_{z} a)\pa_{\zetab}.
\end{equation*}
\end{lemma}

\begin{proof}
Let $(x,w,z,\sigma,\mub,\zetab)$ be the canonical coordinates on $T^*M$
induced by local coordinates $(x,w,z)$ on $M.$ In $M^\circ,$
$T^*_{M^\circ}M$ and $\Tbstar_{M^\circ}M$ are naturally identified via the
map 
\begin{equation*}
\bpi:T^*M\to\Tbstar M,\ \bpi(x,y,z,\sigma,\mub,\zetab)
=(x,w,z,x\sigma,\mub,\zetab). 
\label{30.8.2006.34}\end{equation*}
Moreover,
\begin{equation*}
H_{\bpi^*a}=\pa_\sigma \bpi^*a \pa_x+\pa_{\mub}\bpi^*a \pa_{w}
+\pa_{\zetab}\bpi^*a \pa_{z}
-\pa_x \bpi^*a \pa_\sigma-\pa_{w} \bpi^*a\pa_{\mub}
-\pa_{z} \bpi^*a\pa_{\zetab}.
\end{equation*}
Since $\pa_\sigma \bpi^*a=x\bpi^*\pa_\xib a$ and $\pa_x \bpi^*a=\bpi^*(\pa_x
+x^{-1}\xib\pa_\xib)a,$
\begin{equation*}\begin{split}
H_{\bpi^*a}\bpi^*c
=\bpi^*(&\pa_\xib a (x\pa_x+\xib\pa_\xib)c
+\pa_{\mub} a \pa_{w}c+\pa_{\zetab} a \pa_{z}c\\
&-(x\pa_x a+\xib\pa_\xib a)\pa_\xib c
-\pa_{w} a\pa_{\mub}c-\pa_{z} a\pa_{\zetab}c),\ \Mif c\in\CI(\Tbstar M).
\end{split}\end{equation*}
The terms in which both $a$ and $c$ are differentiated with respect to
$\xib$ cancel, proving the lemma.
\end{proof}

The test operator used below has the special property that its symbol is
fiber constant at the compressed cotangent bundle:

\begin{definition}\label{def:basicop}
A symbol $a\in \CI(\Tbstar M)$ of order $m$ is said to be
\emph{basic} if it is constant on the fibers above $\Tbdotstar(M),$
i\@.e\@. in terms of local coordinates $\pa_z(a)=0$ at
$\{x=\xib=\zetab=0\}.$  An operator having such a principal symbol is also
said to be \emph{basic}.
\end{definition} 

\begin{lemma}\label{lemma:commutators}
If $A\in \Psib^{m}(M)$ there exist $B\in \Psib^{m}(M),$ $C \in
\Psib^{m-1}(M)$, depending continuously on $A$, such that
$$
[D_x, A] = B + C D_x
$$
with $\sigma(B) = -i\pa_x(\sigma (A)),$ $\sigma(C) = -i\pa_\xib(\sigma(A)).$
If in addition $A$ is a basic operator,
\begin{equation}\label{eq:D_z-comm}
[x^{-1} D_{z_i}, A] = B_i+ C_i D_x + \sum_j E_{ij} x^{-1} D_{z_j}+ x^{-1} F_i
\end{equation}
with $B_i \in \Psib^{m}(M),$ $C_i,$ $E_{ij},$ $F_i \in \Psib^{m-1}(M)$
depending continuously on $A$ and
\begin{equation}
-i\pa_{z_i}(\sigma(A))-i\zetab_j\pa_\xib(\sigma(A)=x\sigma(B_i)+\xib\sigma(C_i)+
\sum_{j}\zetab_j\sigma (E_{ij}).
\label{30.8.2006.20}\end{equation}
\end{lemma}

\begin{proof}
Writing
\begin{equation*}
[D_x,A]=[x^{-1}(xD_x),A]=x^{-1}[xD_x,A]+[x^{-1},A]xD_x,
\end{equation*}
we recall from the end of \S\ref{section:calculi} that $[xD_x, A]$ is one
power of $x$ more regular than the commutator of two generic b-operators,
i.e.\ $[xD_x,A]\in x\Psib^{m}(M)$, so setting
$B=x^{-1}[xD_x,A]\in\Psib^{m}(M)$ Lemma~\ref{lemma:b-Ham-vf} shows
$\sigma_m(B)=-ix^{-1}H_{\xib}a=-i\pa_x a$ with $a=\sigma (A).$ Moreover,
$C=[x^{-1},A]x\in\Psib^{m-1}(M)$ has
$\sigma_{m-1}(C)=-iH_{x^{-1}}a=-i\pa_\xib a.$

On the other hand, $x^{-1}D_{z_j}\in x^{-1}\Psib^1(M)$, so
$[x^{-1}D_{z_j},A]\in x^{-1}\Psib^{m}(M)$, and
$\sigma([x^{-1}D_{z_j},A])=-i(x^{-1}\zetab_j\pa_\xib a+x^{-1}\pa_{z_j}a)$.
By assumption, $\pa_{z_j}a$ vanishes at $x=0$, $\zetab=0$, $\xib=0$, so
we can write $\pa_{z_j} a=xb_j+\xib c_j+\sum_k\zetab_k e_{jk}$
with $b_j$ homogeneous of degree $m$, $c_j$, $e_{jk}$ homogeneous of
degree $m-1$.
Now let $B_j\in\Psib^m(M)$ with $\sigma(B_j)=-ib_j$,
$C_j, E_{jk}\in\Psib^{m-1}(M)$
with $\sigma(C_j)=-ic_j$,
$\sigma(E_{jk})=-i(e_{jk}+\delta_{jk}\pa_\xib a)$. Then
$R=[x^{-1}D_{z_j},A]-B_j-C_jD_x-\sum E_{jk}D_{z_k}$ satisfies
$R\in x^{-1}\Psib^m(M)$ and $\sigma(R)=0$, so $R\in x^{-1}\Psib^{m-1}(M)$,
i.e.\ $R=x^{-1}F$ with $F\in\Psib^{m-1}(M)$, finishing the proof of
\eqref{eq:D_z-comm}.
\end{proof}

\begin{remark}\label{rem:localize}
Since $M$ is non-compact and the results here are microlocal, we shall fix
a compact set $\tilde K\subset M$ and assume that all pseudodifferential
operators under consideration have their Schwartz kernels supported in
$\tilde K\times\tilde K.$ Choose $\phi\in\CIc(M)$ which is identically
equal to $1$ in a neighborhood of $\tilde K$ and fiber constant. Below we
use the notation $\|.\|_{\domt_\loc}$ for $\|\phi u\|_{\domt}$ to avoid
having to specify $\tilde U.$ We also write $\|\phi v\|_{\domt'}$ as
$\|v\|_{\domt'_{\loc}}.$
\end{remark}

Next comes the crucial estimate. 

\begin{lemma}\label{lemma:Dirichlet-form}
Suppose that $K\subset U\subset\Sbstar X$ with $K$ compact and $U$ open and
that $\cA=\{A_r\in\Psib^{s-1}(X):\ r\in(0,1]\}$ a basic family with
$\bWF'(\cA)\subset K$ which is bounded in $\Psibc^s(X).$
Then there exist $G\in\Psib^{s-1/2}(X),$ $\tilde G\in\Psib^{s+1/2}(X)$
with $\bWF'(G),$ $\bWF'(\tilde G)\subset U$
and $C_0>0$ such that
\begin{multline*}
\left\lvert\int_M \left(|d_X A_r u|^2-|D_t A_r u|^2\right)\right\rvert\\
\le C_0(\|u\|^2_{\domt_{\loc}}+\|Gu\|^2_{\domt}+\|\Box u\|^2_{\domt'_{\loc}}
+\|\tilde G \Box u\|^2_{\domt'}),\\
\forall\ r\in(0,1],\ u\in H^1_{0,\loc}(X)\\
\Mwith \dWF^{s-1/2}(u)\cap U=\emptyset,\
\mWF^{s+1/2}(\Box u)\cap U=\emptyset
\end{multline*}
where the meaning of $\|u\|^2_{\domt_\loc}$ and $\|\Box u\|^2_{\domt'_\loc}$
is stated above in Remark~\ref{rem:localize}.
\end{lemma}

\noindent It follows that if $\Box u=0$
\begin{equation*}
\left\lvert\int_M \left(|d_X A_r u|^2-|D_t A_r u|^2\right)\right\rvert\leq
C_0(\|u\|^2_{\domt_\loc}+\|Gu\|^2_{\domt}).
\end{equation*}
The main point of this lemma is that $G$ is $1/2$ order lower than the
family $\cA.$ In the limit, $r\to 0,$ this gives control of the Dirichlet
form evaluated on $A_0u,$ $A_0\in\Psibc^s(M),$ in terms of lower order
information. The role of $A_r$ for $r>0$ is to regularize such an argument,
i.e.\ to ensure that the terms in a formal computation, in which one uses
$A_0$ directly, actually make sense.

\begin{proof}
The assumption on the wavefront set of $u$ implies that $A_r u\in\domt$ for
$r\in(0,1]$, so, writing $\langle\cdot,\cdot\rangle$ for the pairing in $L^2(M),$
\begin{equation*}
\|d_X A_r u\|^2_{L^2(M)}-\|D_t A_r u\|^2_{L^2(M)}=
-\langle \Box A_r u,A_r u\rangle.
\end{equation*}
Here the right hand side is the pairing of $\domt'$ with $\domt.$
Writing $\Box A_r=A_r\Box +[\Box ,A_r],$ the right hand side can be estimated by
\begin{equation}\label{eq:Dirichlet-form-8}
|\langle A_r\Box u,A_r u\rangle|+|\langle [\Box ,A_r] u,A_r u\rangle|.
\end{equation}
The lemma is therefore proved once it is shown that the first term of
\eqref{eq:Dirichlet-form-8} is bounded by
\begin{equation}\label{eq:Dirichlet-bd}
C_0'(\|u\|^2_{\domt_\loc}+\|Gu\|^2_{\domt}+\|\Box u\|^2_{\domt'_\loc}
+\|\tilde G \Box u\|^2_{\domt'}),
\end{equation}
and the second term is bounded by
$C_0''(\|u\|^2_{\domt_\loc}+\|Gu\|^2_{\domt}).$

The first estimate is straightforward. Let
$\Lambda_{-1/2}\in\Psib^{-1/2}(M)$ be elliptic with
$\Lambda_{1/2}\in\Psib^{1/2}(M)$ a parametrix (hence also elliptic), so
\begin{equation*}
E=\Lambda_{-1/2}\Lambda_{1/2}-\Id,\quad E'=\Lambda_{1/2}\Lambda_{-1/2}-\Id\in\Psib^{-\infty}(M).
\end{equation*}
Then
\begin{equation*}\begin{split}
\langle A_r\Box u,A_r u\rangle&=\langle(\Lambda_{-1/2}\Lambda_{1/2} -E)A_r \Box u,A_r u\rangle\\
&=\langle\Lambda_{1/2} A_r \Box u,\Lambda_{-1/2}^* A_r u\rangle
-\langle A_r \Box u,E^*A_r u\rangle.
\end{split}\end{equation*}
Since $\Lambda_{1/2}  A_r$ is uniformly bounded in $\Psibc^{s+1/2}(M)$,
and $\Lambda_{-1/2}^* A_r$ is uniformly bounded in $\Psibc^{s-1/2}(M)$,
$\langle \Lambda_{1/2} A_r \Box u,\Lambda_{-1/2}^* A_r u\rangle$
is uniformly bounded, with a bound like \eqref{eq:Dirichlet-bd}
using Cauchy-Schwarz and Lemma~\ref{lemma:WFb-mic-q}.
Indeed, using Lemma~\ref{lemma:WFb-mic-q},
choosing any $G\in\Psib^{s-1/2}(M)$ which is elliptic on
$K$, there is a constant $C_1>0$ such that
\begin{equation*}
\|\Lambda_{-1/2}^*A_r u\|^2_{\domt}\leq
C_1(\|u\|^2_{\domt_\loc}+\|Gu\|^2_{\domt}).
\end{equation*}
Similarly, by Lemma~\ref{lemma:WFb-mic-q} (or more precisely, its
version with $\mWF$) choosing any $\tilde G\in\Psib^{s+1/2}(M)$ which is
elliptic on $K,$ there is a constant $C_1'>0$ such that
$\|\Lambda_{1/2}  A_r\Box u\|^2_{\domt'}\leq
C_1'(\|\Box u\|^2_{\domt'_\loc}+\|\tilde G\Box u\|^2_{\domt'})$.
Combining these gives, with $C_0'=C_1+C_1',$
\begin{multline*}
|\langle\Lambda_{1/2} A_r\Box u,\Lambda_{-1/2}^* A_r u\rangle|\leq
\|\Lambda_{1/2} A_r \Box u\|
\,\|\Lambda_{-1/2}^*A_r u\|\\
\leq \|\Lambda_{1/2} A_r \Box u\|^2+\|\Lambda_{-1/2}^*A_r u\|^2\\
\leq C_0'(\|u\|^2_{\domt_\loc}+\|Gu\|^2_{\domt}
+\|\Box u\|^2_{\domt'_\loc}+\|\tilde G\Box u\|^2_{\domt'}),
\end{multline*}
as desired.

A similar argument using the assumption that
$A_r$ is uniformly bounded in $\Psibc^{s}(M),$ and the 
uniform boundedness of $E^* A_r$ in $\Psibc^{s-1/2}(M)$
(in fact it is bounded in $\Psibc^{-\infty}(M),)$ shows that $\langle A_r
\Box u,E^* A_r u\rangle$ is uniformly bounded.

Now we turn to the second term in \eqref{eq:Dirichlet-form-8}.
Let $Q_i\in x^{-1}\Diffe{1}(M)$ be a local basis of $x^{-1}\Diffe{1}(M)$
as a $\CI(M)$-module.
Using Lemma~\ref{lemma:commutators},
\begin{equation*}
[\Box ,A_r]=\sum_{i,j}Q_i Q_j B'_{ij,r}+\sum Q_j B'_{j,r}+B'_r,
\end{equation*}
$B'_r\in\Psib^{s}(M)$,
$B'_{j,r}\in\Psib^{s-1}(M)$, $B'_{ij,r}\in\Psib^{s-2}(M)$,
uniformly bounded in $\Psibc^{s+1}(M)$, resp.\ $\Psibc^{s}(M)$,
resp.\ $\Psibc^{s-1}(M)$.
With $\Lambda_{-1/2}\in\Psib^{-1/2}(M)$ as above, using
Lemma~\ref{lemma:rearrange},
we can write further
\begin{equation*}
\Lambda_{1/2} [\Box ,A_r]=\sum_{i,j}Q_i Q_j B_{ij,r}
+\sum Q_j B_{j,r}+B_r,
\end{equation*}
with $B_r$, etc, as their primed analogues, but of order $1/2$ greater:
$B_r\in\Psib^{s+1/2}(M)$,
$B_{j,r}\in\Psib^{s-1/2}(M)$, $B_{ij,r}\in\Psib^{s-3/2}(M)$,
uniformly bounded in $\Psibc^{s+3/2}(M)$, resp.\ $\Psibc^{s+1/2}(M)$,
resp.\ $\Psibc^{s-1/2}(M)$.
Thus,
\begin{equation}\begin{split}\label{eq:ell-16}
\langle [\Box ,A_r]u,A_ru\rangle
&=\sum_{ij}\langle Q_iQ_j B_{ij,r}u,\Lambda_{-1/2}^*A_r u\rangle
-\sum_{ij}\langle Q_iQ_j E B_{ij,r}u,A_r u\rangle\\
&\quad+\sum_j \langle Q_j B_{j,r} u,\Lambda_{-1/2}^*A_r u\rangle
-\sum_j \langle Q_jE B_{j,r} u,A_r u\rangle\\
&\quad+\langle B_r u,\Lambda_{-1/2}^*A_r u\rangle
-\langle E B_r u,A_r u\rangle\\
&=\sum_{ij}\langle Q_jB_{ij,r}u,Q_i^*\Lambda_{-1/2}^*A_r u\rangle
-\sum_{ij}\langle Q_j E B_{ij,r}u,Q_i^* A_r u\rangle\\
&\quad+\sum_j \langle B_{j,r} u,Q_j^*\Lambda_{-1/2}^*A_r u\rangle
-\sum_j \langle E B_{j,r} u,Q_j^*A_r u\rangle\\
&\quad
+\langle B_r u,\Lambda_{-1/2}^*A_r u\rangle
-\langle E B_r u,A_r u\rangle,
\end{split}\end{equation}
where $Q_i^*\in x^{-1}\Diffe{1}(M)$
is the formal adjoint of $Q_i$ with respect to $dg$, and
where in the last step we used
\begin{equation*}
B_{ij,r}u,\Lambda_{-1/2}^*A_r u,\ EB_{ij,r}u,\ A_ru\in \domt.
\end{equation*}
We estimate the term
$$
|\langle Q_j B_{ij,r}u,Q_i^*\Lambda_{-1/2}^*A_r u\rangle|
$$ by Cauchy-Schwarz: both factors are uniformly bounded for $r\in(0,1]$
since $\Lambda_{-1/2}^* A_r$, $B_{ij,r}$ are uniformly bounded in
$\Psibc^{s-1/2}(M)$ with a uniform wave front bound disjoint from
$\dWF^{s-1/2}(u)$.  Indeed, as noted above, by Lemma~\ref{lemma:WFb-mic-q},
choosing any $G\in\Psib^{s-1/2}(M)$ which is elliptic on $K$, there is a
constant $C_1>0$ such that this term is bounded by
$C_1(\|u\|^2_{\domt_\loc}+\|Gu\|^2_{\domt})$.  Similar estimates apply to
the other terms on the right hand side of \eqref{eq:ell-16} (with the
slight technical point that for the penultimate term one uses the pairing
between $\domt^{-1}$ and $\domt$), showing that $\langle [\Box ,A_r] u,A_r
u\rangle$ is uniformly bounded for $r\in(0,1]$, indeed is bounded by
$C_0(\|u\|^2_{\domt_\loc}+\|Gu\|^2_{\domt})$, proving the lemma.
\end{proof}

Next we refine Lemma~\ref{lemma:Dirichlet-form} by including a second parameter
and arranging that $\tilde G$ have order $s$ as opposed to $s+1/2.$

\begin{lemma}\label{lemma:Dirichlet-form-2} Under the same hypotheses as
  Lemma~\ref{lemma:Dirichlet-form} there exist $G\in\Psib^{s-1/2}(X)$ and
  $\tilde G\in\Psib^{s}(X)$ with $\bWF'(G),$ $\bWF'(\tilde G)\subset U$
and $C_0>0$ such that for $\ep>0,$ $r\in(0,1]$ and $u\in\domt_\loc$
with $\dWF^{s-1/2}(u)\cap U=\emptyset,$ $\mWF^{s}(\Box u)\cap U=\emptyset$
\begin{equation*}\begin{split}
|\int_M \left(|d_X A_r u|^2-|D_t A_r u|^2\right)|\leq
\ep&\|d_M A_r u\|^2_{L^2(M)}
+C_0(\|u\|^2_{\domt_\loc}+\|Gu\|^2_{\domt}\\
&+\ep^{-1}\|\Box u\|^2_{\domt'_\loc}
+\ep^{-1}\|\tilde G \Box u\|^2_{\domt'}).
\end{split}\end{equation*}
\end{lemma}

\begin{proof}It is only necessary to treat the term $|\langle A_r\Box u,A_r
 u\rangle|$ slightly differently, using Cauchy-Schwarz:
\begin{equation*}
|\langle A_r\Box u,A_r u\rangle|\leq \|A_r \Box u\|_{\domt'}
\|A_r u\|_{\domt}\leq \ep\|A_r u\|^2_{\domt}+\ep^{-1}
\|A_r \Box u\|^2_{\domt'}.
\end{equation*}
Now the lemma follows by using Lemma~\ref{lemma:WFb-mic-q}.
Definition~\ref{Def:dWF}, namely
choosing any $\tilde G\in\Psib^{s}(M)$ which is elliptic on
$K$, there is a constant $C_1'>0$ such that
$\|A_r\Box u\|^2_{\domt'}\leq
C_1'(\|\Box u\|^2_{\domt'_\loc}+\|\tilde G\Box u\|^2_{\domt'});$
the end of the proof follows that of Lemma~\ref{lemma:Dirichlet-form}.
\end{proof}

\begin{proposition}\label{prop:elliptic}(Microlocal elliptic regularity.)
If $u\in \domt_\loc$ then
\begin{equation*}
\dWF^m(u)\subset \mWF^m(\Box u)\cup\dot\pi(\Sestar M),\Mand
\dWF^m(u)\cap\dcE\subset \mWF^m(\Box u).
\end{equation*}
\end{proposition}

\begin{proof}
We first prove a slightly weaker result in which
$\mWF^m(\Box u)$ is replaced by $\mWF^{m+1/2}(\Box u)$---we rely on
Lemma~\ref{lemma:Dirichlet-form}. We then prove the
original statement using Lemma~\ref{lemma:Dirichlet-form-2}.

Suppose that either $q\in \Sbdotstar M\setminus\dot\pi(\Sestar M)$ or
$q\in\dcE$. We may assume iteratively that $q\notin\dWF^{s-1/2}(u)$;
we need to prove then that $q\notin\dWF^{s}(u)$ (note that the inductive
hypothesis holds for $s=1/2$ since $u\in \domt_\loc$).
Let $A\in\Psib^{s}(M)$ be basic and such that
$\bWF'(A)\cap \dWF^{s-1/2}(u)=\emptyset$, $\bWF'(A)\cap\mWF^{s+1/2}(\Box u)
=\emptyset$, with $\bWF'(A)$ in a small
neighborhood $U$ of $q$ so that for a suitable $C>0$ or $\ep>0$, in $U$
\begin{enumerate}
\item
$1<C (\xibs^2+|\zetabs|^2)$ if $q\in \Sbdotstar M
\setminus\dot\pi(\Sestar M)$,\ \text{\emph{or}}
\item
$|\xibs|+|\zetabs|<\ep(1+|\etabs|^2)^{1/2}$, and
$\etabs>1+\ep$, if $q\in\dcE$.
\end{enumerate}
Let $\Lambda_r\in\Psib^{-2}(M)$ for $r>0$, such that $\cL=\{\Lambda_r: \
r\in(0,1]\}$ is a bounded family in $\Psib^0(M)$, and $\Lambda_r\to\Id$ as
$r\to 0$ in $\Psib^{\tilde\ep}(M)$, $\tilde\ep>0;$ we let the symbol of
$\Lambda_r$ be $(1+r(\taub^2+|\etab|^2+|\zetab|^2+|\xib|^2))^{-1}$. Let
$A_r=\Lambda_r A$.  Let $a$ be the symbol of $A$, hence for $r>0,$
$$\sigma(A_r)=(1+r(\taub^2+|\etab|^2+|\zetab|^2+|\xib|^2))^{-1}a.$$ We now have
$A_r\in \Psib^{s-2}(M)$ for $r>0$, and $A_r$ is uniformly bounded in
$\Psibc^{s}(M)$, $A_r\to A$ in $\Psibc^{s+\tilde\ep}(M)$.

By Lemma~\ref{lemma:Dirichlet-form},
\begin{equation*}
\|d_X A_r u\|^2-\|D_t A_r u\|^2
\end{equation*}
is uniformly bounded for $r\in(0,1]$. Write the dual metric as
\begin{equation*}\begin{split}
g^{-1}=\pa_x^2+&2\sum_i B_i\pa_x\pa_{y_i}
+2\sum_i C_i \pa_x\, (x^{-1}\pa_{z_i})
+\sum_{i,j}B_{ij}\pa_{y_i}\pa_{y_j}\\
&+\sum_{i,j} C_{ij} (x^{-1}\pa_{z_i})
(x^{-1}\pa_{z_j})
+2\sum_{i,j}E_{ij}\pa_{y_i}\,(x^{-1}\pa_{z_j}).
\end{split}\end{equation*}
Then
\begin{multline*}
\int_M |d_X A_r u|^2\, dg
=\int_M \abs{D_x A_r u}^2\, dg+\int_M \sum B_i D_x
A_ru\,\overline{D_{y_i}A_ru} \, dg\\
+\int_M \sum C_i D_x A_ru\,\overline{x^{-1}D_{z_i}A_ru} \, dg
+\int_M \sum B_{ij} D_{y_i}A_r u \,\overline{D_{y_j} A_r u} \, dg\\
+\int_M \sum C_{ij} x^{-1}D_{z_i}A_r u \,\overline{x^{-1} D_{z_j} A_r u}\,
dg
+\int_M \sum E_{ij} D_{y_i}A_r u \,\overline{x^{-1} D_{z_j} A_r u} \, dg.
\end{multline*}
As $B_{ij}(x,y,z)=B_{ij}(0,y,z)+x B'_{ij}(x,y,z)$
(in fact, $B(0,y,z)$ is independent of $z),$ we see that
if $A_r$ is supported in $x<\delta$,
\begin{equation}\label{eq:elliptic-28}
|\int_M xB'_{ij} D_{y_i}A_r u \,\overline{D_{y_j} A_r u}\, dg|
\leq C\delta\sum_{i',j'}\|D_{y_{i'}}A_r u\|\,\|D_{y_{j'}}A_r u\|,
\end{equation}
with analogous estimates for $C_{ij}(x,y,z)-C_{ij}(0,y,z)$ and
for $B_i(x,y,z)$, $C_i(x,y,z)$ and $E_{ij}(x,y,z)$.
Moreover, as the matrix $B_{ij}(0,y,z)$ is positive definite,
for some $c>0$,
\begin{equation*}
c\int_M \sum_j |x^{-1}D_{z_j}A_r u|^2\, dg \leq
\frac{1}{2}\int_M \sum_{ij} B_{ij}(0,y,z)
x^{-1}D_{z_i}A_r u\, \overline{x^{-1} D_{z_j} A_r u}\, dg;
\end{equation*}
we also make $c<1/2$.
Thus,
there exists $\tilde C>0$ and $\delta_0>0$ such that if $\delta<\delta_0$
and $A$ is supported
in $x<\delta$ then
\begin{equation}\begin{split}\label{eq:ell-32}
&c\int_M (|D_x A_r u|^2+\sum_j |x^{-1} D_{z_j} A_r u|^2)\, dg\\
&\qquad
+\int_M((1-\tilde C\delta)\sum_j |D_{y_j} A_r u|_h^2-|D_t A_r u|^2)\, dg\\
&\qquad\qquad\leq \int_M (|d_X A_r u|^2-|D_t A_r u|^2)\, dg,
\end{split}\end{equation}
where we used the notation
\begin{equation*}
\sum_j |D_{y_j} A_r u|_h^2=\sum_{ij} B_{ij}(0,y,z) D_{y_i}A_r u\,\overline{D_{y_j}A_r u},\end{equation*}
i.e.\ $h$ is the dual metric $g$ restricted to the span of the $dy_j$,
$j=1,\ldots,l$.

Now we distinguish the cases $q\in\dcE$ and
$q\in \Sbdotstar M\setminus\dot\pi(\Sestar M)$.
If $q\in\dcE$, $A$ is supported near $\dcE$,
we choose $\delta\in(0,(2\tilde C)^{-1})$ so that
$(1-\tilde C\delta)\etabs^2>1+\delta$ on
a neighborhood of $\bWF'(A)$, which is possible in view of (2) at the
beginning of the proof.
Then the second integral on the left hand side of \eqref{eq:ell-32}
can be written
as $\| BA_r u\|^2$, with the symbol of $B$ given by
$((1-\tilde C\delta)|\etab|^2
-\taub^2)^{1/2}$ (which is $\geq\delta\taub$), modulo a term
\begin{equation*}
\int_M F A_r u\,\overline{A_r u},\ F\in\Psib^{1}(M).
\end{equation*}
But this expression is uniformly bounded as $r\to 0$ by the
argument above.
We thus deduce that
\begin{equation*}
c\left( \norm{D_x A_r u}^2+\sum_j\norm{x^{-1}D_{z_j} A_r u}^2\right)+\| BA_r u\|^2
\end{equation*}
is uniformly bounded as $r\to 0$.

If $q\in \Sbdotstar M\setminus\dot\pi(\Sestar M)$,
and $A$ is supported in $x<\delta$,
\begin{equation*}
\int_M \delta^{-2}|D_{z_j} A_r u|^2\leq \int_M |x^{-1}D_{z_j} A_r u|^2,
\ \int_M \delta^{-2}|xD_x A_r u|^2\leq \int_M |D_x A_r u|^2.
\end{equation*}
On the other hand, near $q\in \Sbdotstar M\setminus\dot\pi(\Sestar M)$,
for $\delta>0$ sufficiently small,
\begin{equation}\label{foobar}\begin{split}
&I =\int_M\left(\frac{c}{2\delta^{2}}
(|xD_x A_r u|^2+\sum_j |D_{z_j}A_r u|^2)-|D_t A_r u|^2\right)\\
&\qquad\qquad=\|BA_r u\|^2
+\int_M F A_r u\,\overline{A_r u},
\end{split}\end{equation}
with the symbol of $B$ given by
$(\frac{c}{2\delta^{2}}(\xib^2+\sum \zetab_j^2)-\taub^2)^{1/2}$
(which does not vanish on $U$ for
$\delta>0$ small), while
$F\in\Psib^{1}(X)$, so the second term on the right hand side
is uniformly bounded as $r\to 0$.
Now the LHS of \eqref{eq:ell-32} is bounded below by
$$
I+ \frac c2 \left( \norm{D_x A_r u}^2+ \sum_j \norm{x^{-1} D_{z_j} A_r
  u}^2 \right).
$$
We thus deduce in this case that
\begin{equation*}
\frac{c}{2}\left(\norm{D_x A_r u}^2+\sum_j\norm{x^{-1}D_{z_j} A_r
u}^2\right)+\| BA_r u\|^2
\end{equation*}
is uniformly bounded as $r\to 0$.

We thus conclude that
$D_x A_r u$, $x^{-1}D_{z_j} A_ru, BA_ru$ are uniformly bounded $L^2(M)$.
Correspondingly
there are sequences $D_xA_{r_k}u$, $x^{-1}D_{z_j}A_{r_k}u$, $BA_{r_k}u$,
weakly convergent in $L^2(M)$, and such that $r_k\to 0$, as $k\to\infty$.
Since they converge to $D_x Au$,
$x^{-1}D_{z_j}Au$, $BAu$, respectively, in $\CmI(M)$,
we deduce that the weak limits are $D_x Au$,
$x^{-1}D_{z_j}Au$, $BAu$, which therefore
lie in $L^2(M)$. Consequently,
$d Au\in L^2(M)$ proving the
proposition with $\mWF^{m}(\Box u)$ replaced by $\mWF^{m+1/2}(\Box u)$.

To obtain the full result observe from Lemma~\ref{lemma:Dirichlet-form-2}
that for any $\ep>0$
\begin{equation*}\begin{split}
&\int_M \left(|d_X A_r u|^2-|D_t A_r u|^2-\ep|d_M A_r u|^2\right)\\
&\qquad=\int_M \left((1-\ep)|d_X A_r u|^2-(1+\ep)|D_t A_r u|^2\right)
\end{split}\end{equation*}
is uniformly bounded above for $r\in(0,1].$ By repeating the argument for
sufficiently small $\ep>0$, the right hand side gives an upper bound for
\begin{equation*}
\frac{c}{2}\int_X (|D_x A_ru|^2+\sum_j| x^{-1}D_{z_j} A_r u|^2)+\| BA_r u\|^2,
\end{equation*}
which is thus uniformly bounded as $r\to 0$ and the result follows.
\end{proof}

\begin{theorem}(Elliptic regularity)
Suppose $\Box u=0$. Then $\dWF(u)\subset \dot\Sigma$.
\end{theorem}
\begin{proof}
Over $M^\circ,$ this is ordinary elliptic regularity.  Over $\pa M,$ we
note that it follows from the first part of Proposition~\ref{prop:elliptic}
that the wavefront set is contained in $\dot\pi(\Sestar M).$  On the other
hand, $\dot\pi(\Sestar M) \backslash \dot\Sigma = \dcE,$ so the result then
follows from the second part of the proposition.
\end{proof}

\section{Law of reflection}\label{sec:reflection}

In this section we show that singularities interacting with the boundary
may be microlocalized to give a propagation theorem in the slow variables,
i\@.e\@.\ constraining the values of $t,y$ and their dual variables at which
singularities come off the boundary. The argument is global in the fast
variables in the fiber (and in a certain sense the dual to the normal
variable is also `fast'). Recall that we may freely regard $\dWF(u)$ as a
subset of $\Sbstar M$ or of $\Sbdotstar M$ (see
Definition~\ref{def:WF-Sbdotstar}).

As is already clear from \cite{Melrose-Wunsch1}, subprincipal terms matter
in the commutator argument, so our arguments are somewhat delicate. To
carry out commutator computations we trivialize the fibration locally near
a point of $Y,$ and hence extend the trivialization to a neighborhood in
$M$ of the fiber above $y.$ Thus an open set $O$ in $M$ is identified with
$[0,\ep)_x\times O'\times Z,$ where $O'\ni y$ is open in $Y$ and we take
$O'$ to be a coordinate patch. Then $\Tbstar_O M$ is identified with
$\Tbstar([0,\ep)\times O')\times T^*Z$, and this product decomposition
allows us to pull back functions from $\Tbstar([0,\ep)\times O')$ to
$\Tbstar_O M.$

In particular this gives a connection for the fibration near the chosen
fiber. Any b-vector field $V$ on $[0,\ep)\times O'$ lifts to $O$ via the
product identification and then for any vector field $W$ tangent to the
fibers, the commutator $[V,W]$ is also tangent to the fibers. Indeed, if
$W$ is the lift of a vector field from $Z,$ the conclusion is clear, and in
general $W$ is a finite linear combination of such vector fields with
coefficients on $M.$ If $W$ is only tangent to the fibers over the boundary
then the same is true of $[V,W].$ Also, if $f$ is the pull-back of a
function on $[0,\ep)\times O',$ its commutator with $W$ vanishes. Using
this local trivialization of the fibration we also choose an explicit
quantization map for b-pseudodifferential operators.

\begin{definition}\label{Op-Def}
Fix a partition of unity $\phi_i$ on $Z$ supported in coordinate charts,
cut-offs $\chi_i$ identically equal to $1$ on a neighborhood of
$\supp\phi_i$ and still supported in coordinate charts and a function
$\rho\in\CI_c((-1/2,1/2))$ identically $1$ near $0.$ Then for $a\in
S^m(\Tbstar O)$ set
\begin{multline}
\Op_\trv(a)=\sum_i \chi_i A_i\chi_i,\ \text{where}\\
A_i=(2\pi)^{-n}\int e^{i(\xib\cdot \frac{x-x'}{x'}+\etab\cdot(y-y')
+\zetab\cdot (z-z'))}\\
\rho(\frac{x-x'}{x'})
a(x,y,\xib,\etab,\zetab)\phi_i(z)\,d\xib\,d\zetab\,d\etab\nu. 
\label{30.8.2006.22}\end{multline}
\end{definition}

This is an explicit semi-global quantization, with $\nu$ a choice of right
density, for the algebra of (compactly supported) b-pseudodifferential
operators, as can be seen directly by lifting the kernels to the
b-stretched product on which the cut-off $\rho(x/x'-1)$ is smooth. In
particular it is surjective, modulo terms of order $-\infty,$ for operators
essentially supported near the given fiber over the boundary.
 
To construct a `test operator' we choose a symbol $\tilde a$ on
$\Tbstar([0,\ep)\times O'),$ lift it to $\Tbstar O$ (still denoting it by
$\tilde a$) and multiply it by a cut-off that is identically $1$ near
$\dot\Sigma.$ Thus set
\begin{equation}\label{eq:a-tilde-a}
a=\tilde a \psi(\xibs^2+|\zetabs|^2)
\end{equation}
with $\psi\in\CI_c(\RR)$ supported in $[-2c_1,2c_1]$ and identically equal
to $1$ on $[-c_1,c_1].$ Thus, $d\psi(\xibs^2+|\zetabs|^2)$ is supported in
$\xibs^2+|\zetabs|^2\in [c_1,2c_1]$ and hence outside $\dot\Sigma.$ This
may be thought of as a factor microlocalizing near the characteristic set
but effectively commuting with $\Box.$ Note that $a$ is indeed a symbol on
$\Tbstar O.$

\begin{lemma}\label{lemma:D_z-comm-fine}
With $\tilde a\in S^m(\Tbstar([0,\ep)\times O'))$ given
by \eqref{eq:a-tilde-a}, $A=\Op_\trv(a)$ and $W$ a vector field lifted from
$Z$ using the local trivialization of the fibration,
\begin{equation}\label{eq:D_z-comm-fine}
[W,A] = B\in\Psib^{m}(M)
\end{equation}
has wave front set disjoint from $\dot\Sigma$ and depends continuously on
$a.$
\end{lemma}

\begin{proof} 
Changing $A$ by an operator of order $-\infty$ operator does not affect the
conclusion, so in particular the commutator terms involving the cut-offs
$\chi_i$ in the definition of $\Op_\trv(a)$ may be ignored.

Working modulo $\Psib^{-\infty}(M)$, so localizing near the diagonal as
necessary, the Schwartz kernel of $A$ can be computed using product
coordinate charts $\tilde O\times \tilde O.$ In these,  where we allow
coordinate changes in $z,$ $A$ takes the form
\begin{multline*}
(2\pi)^{-n}\int e^{i(\xib\cdot \frac{x-x'}{x'}+\etab\cdot(y-y')
+\zetab\cdot (z-z'))}\\
\rho(\frac{x-x'}{x'}) \tilde a(x,y,\xib,\etab)
\tilde\psi(x,y,z,\xib,\etab,\zetab)
\,d\xib\,d\zetab\,d\etab,
\end{multline*}
where $d(\tilde\psi)$ is supported away from $\dot\Sigma.$  The commutator
$[A,W]$ thus arises from the derivative of $\tilde\psi$ along $W,$ hence is
microsupported away from $\dot\Sigma.$
\end{proof}

The adjoint of $\Op_\trv(a)$ with respect to $dg$ is not quite of the form
$\Op_\trv(b)$ because we have not used Weyl quantization in
\eqref{30.8.2006.22}, however it is of this form modulo terms of order
$-\infty$ since $\Op_\trv$ is a full quantization of the b-calculus. If
$\nu=\nu_B\nu_Z$ is a product type edge density relative to the
trivialization, with $B= [0,\ep)\times O'$, and $a$ is given by
\eqref{eq:a-tilde-a} then, modulo terms with $\bWF'$ disjoint from
$\dot\Sigma$ (arising from the cut-off $\psi$) the adjoint
$\Op_\trv(a)^\dagger$ of $A$ with respect to $\nu$ reduces to the adjoint
with respect to $\nu_B.$ Since the quantization is of product type, the
full symbol of the adjoint is of the same form \eqref{eq:a-tilde-a}. Thus
$A$ may be arranged to be self-adjoint with respect to $\nu.$ Then $dg=J\nu,$
for some smooth nonvanishing $J,$ and the adjoint $A^*$ of $A$ with
respect to $dg$ is $A^*=J^{-1}A^\dagger J,$ i\@.e\@. if $A^\dagger=A$ then
$A^*=J^{-1}A J.$

\begin{lemma}\label{lemma:D_z-comm-fine-star}
If $A=A^\dagger$ with symbol given by \eqref{eq:a-tilde-a} and $W$ is a
vector field lifted from $Z$ using the trivialization of the fibration,
then $E\in\Psib^{2m}(M),$ and $F\in\Psib^{2m-1}(M)$ may be chosen to
depend continuously on $A$ so that
\begin{equation}
[W, A^*A] = F+E,\ \sigma(F)=a\{a,W(\log J)\},\ \bWF'(E)\cap\dot\Sigma=\emptyset,
\label{eq:D_z-comm-fine-star}\end{equation}
where $W(\log J)\in\CI(M)$ is lifted to $\Tbstar M$ by the bundle projection.
\end{lemma}

\begin{proof}
Working modulo operators which have wavefront set disjoint from $\dot\Sigma,$
$[A,W]\equiv0$ so
\begin{equation*}\begin{split}
[W, A^*A] &\equiv[W,J^{-1}]AJA+J^{-1}A[W,J]A\\
&\equiv -J^{-2}[W,J]AJA+J^{-1}A[W,J]A
\equiv J^{-1}[A,[W,J]J^{-1}]JA.
\end{split}\end{equation*}
Then $F$ and $E$ may be defined by inserting appropriate cut-offs.
\end{proof}

Using the chosen local trivialization of the fibration, set 
\begin{equation}
Q_0=x^{-1},\ Q_1=D_x\Mand Q_j=x^{-1}W_j,\ j\geq 2
\label{23.10.2006.1}\end{equation}
where the $W_j$ are the lifts of vector fields from $Z$ spanning $\Vf(Z)$
over $\CI(Z).$ 

\begin{lemma}\label{lemma:edge-b-comp} The $Q_i$ span the $\CI(M)$-module
$x^{\ell-1}\Diffe{1}(M)/x^{\ell}\Diffb{1}(M)$ over $x^{\ell}\CI(M);$ in
fact for every $\ell,$ 
\begin{multline}
P_1\in x^{\ell-1}\Diffe{1}(M)\Longrightarrow \ \exists\ B_i\in
  x^{\ell}\CI(M)\Mand C\in x^{\ell}\Vb(M)\\
\text{such that }P_1=\sum B_i Q_i+C.
\label{23.10.2006.2}\end{multline}
Similarly, each element $P\in x^{\ell-2}\Diffe{2}(M)$ may be decomposed as
\begin{multline*}
P=\sum Q_i^* B_{ij} Q_j+\sum C_i Q_i+L,\\
\Mwith B_{ij}\in x^{\ell}\CI(M),\ C_i\in x^{\ell}\Vb(M)\Mand L\in x^{\ell}
  \Diffb{2}(M).
\end{multline*}
\end{lemma}

\begin{proof}
If \eqref{23.10.2006.2} holds locally over a covering of $M$ by coordinate
charts then summing, on the left, over a partition of unity subordinate to
this cover gives the semi-global result. Thus, it suffices to work in
coordinates compatible with the chosen trivialization. The vector fields $xD_x,$
$xD_{y_i}$ and $W_j$ span $\Ve(U)$ so if $P_1\in x^{\ell-1}\Diffe{1}(M)$
then over $U,$ $xP_1=a (xD_x)+\sum b_j W_j+\sum c_i (xD_{y_i})+c$ with $a,$
$c,$ $b_j,$ $c_i\in x^{\ell-1}\CI(U).$ This is the local form of
\eqref{23.10.2006.2}:
\begin{equation*}
P_1=cQ_0+aQ_1+\sum b_jQ_j+\sum (c_i x)D_{y_i}.
\label{23.10.2006.3}\end{equation*}

For the second claim, write $P=\sum P_i^* P_i'$ with
$P_i\in x^{-1}\Diffe{1}(M)$, $P_i'\in x^{\ell-1}\Diffe{1}(M)$, and apply
the first part to $P_i$ and $P'_i$.
\end{proof}

\begin{lemma}\label{lemma:Box-form} With the $Q_i$ as in \eqref{23.10.2006.1}
there exist $M_i,$ $N_i\in x\Diffb{1}(M),$ $\Kf_{ij}\in\CI(M)$ and
$H\in\Diffb{2}(M)$ such that
\begin{multline}
-\Box=\sum Q_i^*\Kf_{ij} Q_j+\sum Q_i^*M_i+\sum N_iQ_i+H\\
\Mwith\sigma(H)|_{x=0}=h-\taub^2,
\label{eq:Box-form}\end{multline}
$h$ being the metric on the base $Y.$
\end{lemma}

\begin{proof} If the requirements are weakened to $M_i,$ $N_i\in\Diffb{1}(M)$
and the conclusion regarding $\sigma(H)|_{x=0}$ is dropped then
\eqref{eq:Box-form} follows from Lemma~\ref{lemma:edge-b-comp} applied to
$\Box\in x^{-2}\Diffe{2}(M).$ Moreover, every element of $x^{-1}\Diffe{2}(M)$
has the form \eqref{eq:Box-form} with $M_i,$ $N_i$ as stated, and with
$\sigma(H)|_{x=0}=0.$ So the conclusions of the lemma regarding $M_i,$ $N_i$
and $\sigma(H)|_{x=0}$ only depend on the normal operator of $\Box$ in
$x^{-2}\Diffe{2}(M).$

Now, the normal operator is determined by the restriction of $x^{-2} g$ to
$\Testar_{\pa X}X,$ which, by \eqref{edgemetric1}, is
$x^{-2}(dx^2+h(0,y,dy))+k(0,y,z,0,0,dz).$ The dual of $x^{-2} g$, as a
symmetric edge 2-tensor, i\@.e\@.\ as a symmetric section of $\Testar
X\otimes\Testar X,$ is $x^2(\pa_x^2+H(0,y,\pa_y))+K(0,y,z,0,0,\pa_z)$
modulo sections vanishing at $\pa X$, with $H$, resp\@. $K,$ the duals of
$h$ and $k.$ The standard expression for the Laplacian yields
\begin{equation*}\begin{split}
\Delta=&D_x^*D_x
+\sum_{i,j}D_{y_i}^*H^{ij}(y)D_{y_j}
+\sum_{i,j}D_{z_i}^*x^{-2}K^{ij}(y,z) D_{z_j}+P,
\end{split}\end{equation*}
$P\in x^{-1}\Diffe{2}(M)$. As $\det g=x^{2f}a$, $a\in\CI(X)$, $a>0$,
\begin{equation*}
D_{y_i}^*-D_{y_i}=(\det g)^{-1/2}D_{y_i} (\det g)^{1/2}-D_{y_i}\in\CI(X),
\end{equation*}
and the lemma follows.
\end{proof}

The $Q_i$ in \eqref{23.10.2006.1} are homogeneous of degree $-1$ with
respect to dilations in the $x$ factor in the trivialization. This leads to
the term $4a(\pa_{\xib}a)\Kf_{ij}$ in $\sigma(L_{ij})$ in the following
commutator calculation; ultimately this is the dominant term.

\begin{lemma}\label{cor:Box-comm} With the notation of
Lemma~\ref{lemma:Box-form}, if $A=A^\dagger=\Op_\trv(a)\in\Psib^{m}(M)$ then  
\begin{equation}
\begin{gathered}
i[A^*A,\Box]=\sum Q_i^* L_{ij}Q_j+\sum Q_i^* L'_i+\sum L''_i Q_i+L_0+E\\
\Mwith L_{ij}\in\Psib^{2m-1}(M),\  L'_i,\ L''_i\in\Psib^{2m}(M),\
L_0\in\Psib^{2m+1}(M),\\
E\in x^{-2}\Diffe{2}\Psib^{2m-1}(M)\Mhaving\bWF'(E)\cap\dot\Sigma=\emptyset,\\
\sigma(L_{ij})=4a(\pa_{\xib}a)\Kf_{ij}+aV_{ij}a,\ \sigma(L'_i)=aV'_ia,\\
\sigma(L''_i)=aV''_ia\Mand\sigma(L_0)=aV_0 a
\end{gathered}
\label{23.10.2006.4}\end{equation}
where the $V_{ij},$ $V'_i,$ $V''_i$ and $V_0$ are in the span, over
  homogeneous functions on $\Tbstar M$ (independent of $a$), of
  vector fields on $\Tbstar ([0,\ep)\times Y)$ that are tangent to
  $\{\xib=c\}$ at $x=0$ for each $c,$ and
\begin{equation*}
V_0|_{x=0}=2H_{h-\taub^2}.
\end{equation*}
\end{lemma}
\noindent
We remark that while the factors in
$x^{-k}\Diffe{k}\Psib^m(M)$ may be arbitrarily rearranged, according to
Lemma~\ref{lemma:rearrange} and Corollary~\ref{cor:rearrange}, the
principal symbols of the factors in the rearranged version may involve {\em
derivatives} of the principal symbol of the original operators.  This is
due to the appearance of $Q_0=x^{-1}$ terms which are lower order than the
$Q_j$ in the standard sense, but not in the sense of the grading of
$x^{-k}\Diffe{k}\Psib^m(M).$ This is of no consequence if the terms being
rearranged are {\em a priori} bounded. However, for positive commutator
estimates such as those below, the principal symbols needs to be
controlled, and the basic problem we face is that the differential of a
non-zero function is never bounded by a multiple of the function near the
boundary of its support, so in microlocalizing it is necessary to exclude
such undesirable derivatives, hence the careful structuring of
\eqref{23.10.2006.4}.

\begin{proof} Using the decomposition of $\Box$ in
Lemma~\ref{lemma:Box-form}
\begin{multline}
-i[A^*A,\Box]=i\sum [A^*A,Q_i^*]\Kf_{ij} Q_j+i\sum Q_i^*\Kf_{ij}[A^*A,Q_j]\\
+i\sum [A^*A,Q_i^*]M_i+i\sum N_i[A^*A,Q_i]\\
+i\sum Q_i^*[A^*A,\Kf_{ij}]Q_j
+i\sum Q_i^*[A^*A,M_i]+i\sum [A^*A,N_i]Q_i
+i[A^*A,H].
\label{23.10.2006.5}\end{multline}
The last four terms involve commutators of the form $i[A^*A,R],$ with
$R\in\Diffb{l}(M)$ and Lemma~\ref{lemma:b-Ham-vf} shows that such terms are
in $\Psib^{2m+l-1}(M)$, with principal symbol of the form $aVa$ with $V$ a
homogeneous vector field on $\Tbstar M,$ tangent to $\xib=c$ at $x=0.$
These terms
make obvious contributions as required by \eqref{23.10.2006.4}, so it only
remains to consider the terms involving commutators with the $Q_i$ and $Q_i^*.$

Consider first $Q_1=D_x.$ Lemma~\ref{lemma:commutators} applies and shows
that $i[A^*A,D_x]$ is of the form $B+CD_x$ with $B\in\Psib^{2m}(M),$ and
$C\in\Psib^{2m-1}(M)$ where $\sigma(B)=-\pa_x a^2$ and $\sigma(C)=-\pa_{\xib}
a^2.$ The term involving $B$ can then be absorbed in the $L'_i,$ $L''_i$ or
$L_0$ terms, while $C$ enters into $L_{ij}$ (or $L'_i,$ $L''_i,$ depending
on which factor of $Q_i$ we were considering). Notice though that if $C$
enters into $L'_i,$ i\@.e\@.\ if we were considering $i[A^*A,Q_i^*]M'_i,$ the
vanishing of the principal symbol of $M'_i$ at $x=0$ shows that the
principal symbol of $L'_i$ is indeed given by $aVa$ with $V$ a vector field
tangent to $\xib=c$ at $x=0,$ namely $x\pa_{\xib}.$

Next consider the terms $[A^*A,Q_i]$ with $Q_i=x^{-1}W_i$ for $i\neq 0,$ to
which Lemma~\ref{lemma:D_z-comm-fine-star} applies for all $W_i.$ Since $F$
in Lemma~\ref{lemma:D_z-comm-fine-star} is in $\Psib^{2m-1}(M),$ for all
$W_i$ but $Q_0$ (for which this commutator vanishes), we may regard
$x^{-1}F$ as a term of the form $F'Q_0$, $F'\in\Psib^{2m-1}(M)$, thus
giving terms as in the statement of the lemma.

So it remains to consider $i[A^*A,x^{-1}].$ This lies in $\Psib^{2m-1}(M),$
and by Lemma~\ref{lemma:b-Ham-vf}, it has principal symbol
$-2x^{-1}\,a\pa_{\xib} a.$ This completes the proof.
\end{proof}

In $O_{\delta}=O\cap\{x<\delta\},$ with $|\zetabs|$ arising from the metric
on the fibers
$$
p_0=\sigma_b(\Box)
 = \taub^2-\left(\frac{\xib^2}{x^2} + (1+O(x)) \abs{\etab}^2_h + O(1)
\etab \zetab + (\frac{1}{x^2} + O(1)) \abs{\zetab}^2_k\right).
$$
and for $\delta>0$ sufficiently small,
$\xibs^2+|\zetabs|^2<c_1/2$ in $\dot\Sigma\cap U_0$, $U_0
=\Tbstar O_{\delta}.$
The Hamilton vector field of $p_0$ in $\Tbstar (M)$ is then
\begin{equation}\label{b-hamvf}
\begin{aligned}
H_{p_0} &= \taub \pa_t -\frac \xib x \pa_x -\frac{2(\xib^2+\abs{\zetab}^2)}{x^2}
\pa_\xib -
H_{y,\etab} -\frac 1{x^2} H_{z,\zetab} \\ &+ (O(x \etab) + O(\zetab))
\pa_y+(O(x\etab^2) +O(\etab\zetab) 
+O(\zetab^2/x^2))\pa_\etab \\ &+ (O(\etab)+O(\zetab))\pa_z + (O(\etab\zetab) +
O(\zetab^2))\pa_\zetab\\ &+(O(x \etab^2)+O(x\etab\zetab) +
O(\zetab^2/x))\pa_\xib.
\end{aligned}
\end{equation}
where $H_{y,\etab}$ and $H_{z,\zetab}$ denote respectively the Hamilton
vector fields of $h$ in $(y,\etab)$ (with $x$ as a parameter) and of $k$ in
$(z,\zetab)$ (with $(x,y)$ as parameters). So
\begin{multline}\label{Hpalpha}
-\xibs = -\frac{\xib}{\abs{\taub}}\Msatisfies\\
H_{p_0}(-\xibs) = \frac{1}{\abs\taub}\left( \frac{\xib^2}{x^2} +
\frac{\abs{\zetab}^2}{x^2} + O(x \etab^2) + O (x\etab\zetab) +
O(\zetab^2/x)\right).
\end{multline}
The function $-\xibs$ is therefore a
`propagating variable' decreasing under the flow. The commutator argument
used below makes essential use of this decrease.

Now we are ready to state and prove the diffractive propagation theorem in
the hyperbolic region; the proof is similar to that of the corresponding
result in \cite{Vasy5}.

\begin{theorem}\label{theorem:diffractive}
Let $u$ be an admissible solution of $\Box u=0$ and suppose $U\subset\dot\Sigma$
is a neighborhood of $q_0=(t_0,\taub_0,y_0,\etab_0)\in \dhyp$ then 
\begin{equation}
U\cap \{-\xibs<0\}\cap\WF u = \emptyset\Longrightarrow q_0\notin\dWF u.
\label{23.10.2006.6}\end{equation}
\end{theorem}

\noindent
Since $\WF u$ and $\dWF u$ are the same away from $x=0,$ and $\dWF$ is a
closed set, this result suffices to get propagation of wavefront set into
and back out of the edge boundary along rays striking it transversely
(i.e.\ arriving at points in $\dhyp.$)

\begin{proof} Fix a small coordinate neighborhood $U_0$ of $q_0.$ By
Lemma~\ref{lemma:dot-Sigma-nbhd} and the hypotheses, there exists
$\delta>0$ such that
\begin{equation}\label{eq:prop-rem-9b} \{q:\
|x(q)|^2+|y(q)-y_0|^2+|t(q)-t_0|^2+|\etabs(q)
-\etabs_0|^2<\delta\}\cap\{-\xibs<0\} \cap \WF u=\emptyset
\end{equation}
where $\etabs=\frac{\etab}{\abs{\taub}}.$ Set
\begin{equation}\label{eq:prop-omega-def}
\begin{gathered}
\omega = x^2 + (y-y_0)^2 + (\etabs-\etabs_0)^2 + (t-t_0)^2\Mand\\
\phi=-\xibs +\frac{1}{\beta^2\delta}\omega
\end{gathered}
\end{equation}
with the idea that at the edge, the propagation which we cannot
control is in the `fast variables' $(\xib,\zetab)$ but we may
localize in the remaining (`slow') variables by cutting off in $\omega$
and using a multiple of $-\xibs$ to construct a positive commutator. Note that
\begin{equation}\label{eq:omega-est-a0}
|\taub|^{-1}H_{p_0}\omega=O\left(\sqrt\omega  (\xibs^2/x^2 +\zetabs^2/x^2+1)\right).
\end{equation}

Next we select some cut-off functions. First choose $\chi_0\in\CI(\RR)$ with
support in $[0,\infty)$ and $\chi_0(s)=\exp(-1/s)$ for $s>0.$ Thus,
$\chi_0'(s)=s^{-2}\chi_0(s).$ Take $\chi_1\in\CI(\RR)$ to have support in
$[0,\infty),$ to be equal to $1$ on $[1,\infty)$ and to have
$\chi_1'\geq0$ with $\chi_1'\in\CIc((0,1)).$ Finally, let $\chi_2\in\CIc(\RR)$  
be supported in $[-2c_1,2c_1]$ and be identically equal to $1$ on $[-c_1,c_1].$
Here we will take $c_1$ to be such that 
\begin{equation}
\xibs^2+|\zetabs|^2<c_1/2\Min\dot\Sigma\cap U_0
\label{30.8.2006.22a}\end{equation}
where $|\zetabs|$ is the metric length of fiber covectors with respect
to some fiber metric. Thus, $\chi_2(\xibs^2+|\zetabs|^2)$ is a cut-off such that 
$d\chi_2(\xibs^2+|\zetabs|^2)$ is supported in 
$\xibs^2+|\zetabs|^2\in [c_1,2c_1]$ hence outside $\dot\Sigma.$ Such a
factor microlocalizes near the characteristic set but effectively commutes
with $\Box.$  We shall further insist that all cut-offs and their
derivatives have (up to sign) smooth square roots.

Now consider the test symbol
\begin{equation}
a =\chi_0(1-\phi/\delta)\chi_1(-\xibs/
\delta+1)\chi_2(\xibs^2+|\zetabs|^2).
\label{eq:a_m-def}\end{equation}
We proceed to show that
\begin{equation}
\begin{gathered}
\text{
for given $\beta>0,$ if $\delta>0$ is sufficiently small then }
a\in\CI(\Tbstar M\setminus0)\\
\text{is a basic symbol with support in any preassigned}\\
\text{conic neighborhood of }q_0\in\Sbdotstar M.
\end{gathered}
\label{30.8.2006.23}\end{equation}

To see this, observe that from the choices made for the cut-offs,
$\phi\leq \delta$ and $\xibs\leq\delta$ on the supports. Since
$\omega\geq 0$ in \eqref{eq:prop-omega-def}, $-\xibs\leq \delta$ so 
\begin{equation}
|\xibs|\leq \delta\Longrightarrow 
\omega\leq \beta^2 \delta (\phi+\xibs) \leq \beta^2\delta(\delta+\xibs)\leq
2\delta^2\beta^2.
\label{eq:omega-delta-est}\end{equation}
In view of \eqref{eq:prop-omega-def} and \eqref{eq:prop-rem-9b},
this shows that for any $\beta>0,$ $a\in\CI(\Tbstar M\setminus0)$ is
supported in a given conic neighborhood of $q_0,$ provided $\delta>0$ is
sufficiently small. Since the only dependence on the fiber in \eqref{eq:a_m-def}
comes through the metric, as $|\zetabs|^2,$ and $\chi_2$ is constant near
$0,$ $a$ is constant on the fibers near $\xibs^2+|\zetabs|^2=0,$ which is
to say it is basic; this proves \eqref{30.8.2006.23}.

Now, we may take $A=\Op_\trv(a_0)$ using Definition~\ref{Op-Def} with $a_0$
from \eqref{eq:a_m-def}. To extract some positivity from the commutator of
$A^*A$ with $\Box$ using Lemma~\ref{cor:Box-comm} we need to show that the
term $4a\pa_{\xib} a \Kf_{ij}$ in the symbol of $\sigma (L_{ij})$ dominates
the others. To do this it is enough to show that for a homogeneous vector
field $V,$ tangent to $\xib=0,$ $Va$ can be estimated by
$\pa_{\xib}a,$ if $\beta,$ $\delta>0$ are chosen appropriately.

We shall rewrite the identity \eqref{23.10.2006.4} for this particular
choice of $A,$ using \eqref{eq:Box-form}. Choose a basic operator $\tilde
B\in\Psib^{1/2}(M)$ with 
\begin{equation}\label{eq:tilde-b}
\tilde b=\sigma_{{\bo},1/2}(\tilde B)=|\taub|^{1/2}
\delta^{-1/2}(\chi_0\chi_0')^{1/2}\chi_1\chi_2\in
\CI(\Tbstar M\setminus 0)
\end{equation}
where the regularity uses the properties of the cut-offs, in particular that
they and their derivatives are squares. Thus $\tilde b^2$ occurs as a factor in
any first derivative of $a$ in which the derivatives strike the $\chi_0$
term.  (Note that other derivatives are acceptable either because they are
supported away from $\dot\Sigma$ or because they are supported in the
region where we have assume regularity.)
 Also, choose $C\in\Psib^0(M)$
with symbol 
\begin{equation}
\sigma_{{\bo},0}(C)=|\taub|^{-1}(2\taub^2-2|\etab|^2_y)^{1/2}\psi
\label{30.8.2006.24}\end{equation}
where $\psi\in S^0(\Tbstar M)$ is identically $1$ on $U$ considered as a subset
of $\Tbstar M.$ Then, recalling that the $Q_i$ are given in
\eqref{23.10.2006.1} we proceed to show that
\begin{equation}
\begin{gathered}
i[A^*A,\Box]
=\sum Q_i^* L_{ij}Q_j+\sum Q_i^* L'_i+\sum L''_i Q_i+L_0+E\\
=R' \Box+\tilde B^*(C^*C+R_0+\sum_i Q_i R_i
+\sum_{ij} Q_iR_{ij}Q_j)\tilde B+R''+E'+E''
\end{gathered}
\label{eq:P-comm}\end{equation}
where all factors have wavefront set near $q_0,$ 
\begin{equation}\begin{gathered}
R_0\in\Psib^0(M),\ R_i\in\Psib^{-1}(M),\ R_{ij}\in\Psib^{-2}(M),\\
R'\in\Psib^{-1}(M),\ R''\in x^{-2}\Diffe{2}\Psib^{-2}(M),\\
E',\ E''\in x^{-2}\Diffe{2}\Psib^{-1}(M),
\end{gathered}\label{comm-orders}\end{equation}
have symbols
\begin{equation}
\begin{gathered}
r_0=\sigma_{{\bo},0}(R_0),\ r_i=\sigma_{{\bo},-1}(R_i),\ r_{ij}\in\sigma_{{\bo},-2}(R_{ij})\\
\Msatisfying\supp(r_j)\subset\{\omega\leq 9\delta^2\beta^2\},\
|r_0|\leq C_2(1+\frac{1}{\beta^2\delta})\omega^{1/2},
\\
|\taub r_i|\leq C_2(1+\frac{1}{\beta^2\delta})\omega^{1/2},\
|\taub^2 r_{ij}|\leq C_2(1+\frac{1}{\beta^2\delta})\omega^{1/2},
\\
|r_0|\leq 3C_2(\delta\beta+\beta^{-1}),\
|\taub r_i|\leq 3C_2(\delta\beta+\beta^{-1}),\\
|\taub^2 r_{ij}|\leq 3C_2(\delta\beta+\beta^{-1})
\end{gathered}\label{comm-est}
\end{equation}
and 
\begin{equation}
\bWF'(E')\subset\xibs^{-1}((0,\infty))
\cap U,\ \bWF'(E'')\cap\dot\Sigma=\emptyset.
\label{30.8.2006.25}\end{equation}

The first equality in \eqref{eq:P-comm} is just \eqref{23.10.2006.4}. The
operators $L'_i$ have symbols as described in \eqref{23.10.2006.4} and in
particular these are smooth multiples of $\tilde b^2$ plus
noncharacteristic terms and terms supported in $\xibs^{-1}((0,\infty))$
so we may write $L'_i=\tilde B^*R_i\tilde B+R'_i+E'+E''$ where
$R'_i\in\Psib^{-2}(M)$ and the $E',$ $E''$ terms satisfy
\eqref{30.8.2006.25}; henceforth, the meaning of the terms $E'$ and $E''$
(and corresponding symbols $e',e''$) will change from line to line, but all
will satisfy \eqref{30.8.2006.25}.  Thus
\begin{equation}
\sum\limits_{i}Q^*_iL_i'=\tilde B^*\sum\limits_{i}Q^*_iR_i
\tilde B^*+\sum\limits_{i}[Q^*_i,\tilde B^*]R_i\tilde B +\sum\limits_{i}Q^*_iR'_i.
\label{30.8.2006.30}\end{equation}
Since $\tilde B$ is basic, the commutator term may be absorbed into $E''.$
The final term may be absorbed into $R''$ in \eqref{eq:P-comm}. The
$L''_iQ_i$ terms are similar, with a final commutation of $Q_i$ to the
left. Similarly, $L_0$ gives rise to a term involving $R_0$ with additional
contributions to $R''$ and $E''.$

Now consider the first term in which the $L_{ij}$ appear. The definitions
above give
\begin{multline}
4a\pa_{\xib}a =-\frac{4}{\abs{\taub}^2}\tilde b^2+r_{ij}+e'+e''\\
\Longrightarrow 
L_{ij}=B^*(T^* \Kf_{ij}T+R_{ij})B+L'_{ij}+E'+E'',\ L'_{ij}\in\Psib^{-2}(M)
\label{30.8.2006.26}\end{multline}
where $T \in \Psib^{-1}(M)$ is elliptic.  The sum involving the $L_{ij}$ in \eqref{eq:P-comm} can therefore be
written, using \eqref{eq:Box-form}
\begin{multline*}
\sum\limits_{i} Q_i^* L_{ij}Q_j=\tilde R' \Box \\
+ B^*(-H-\sum_i Q^*_i M_i-\sum_i N_iQ_i+\sum_{ij} Q_iR_{ij}Q_j)
\tilde B+\sum Q_i^* L'_{ij}Q_j+E'+E''
\label{30.8.2006.27}\end{multline*}
where the $E',$ $E''$ terms now involve commutators of the $Q_i$ or $Q_i^*$
and $\tilde B$ or $\tilde B^*,$ as well as commutators involving $T, T^*.$
This leads to the final line in \eqref{eq:P-comm}.

By Lemma~\ref{cor:Box-comm}, the symbol $r_{ij}$ can be written as $\chi_0
V \chi_0$ times $\chi_1$ and $\chi_2$ factors, with $V$ tangent to $\xib=c$
at $x=0.$ (We have absorbed $\chi_1',$ $\chi_2'$ terms in $E',E''.$) As we
can estimate $V \omega$ by a multiple of $\sqrt{\omega}$ and $V\xib$ by a
multiple of $x$ and hence of $\sqrt{\omega}$, we obtain the desired symbol
estimates.  The estimates on $r_i$ and $r_0$ proceed analogously.

Having calculated the commutator in \eqref{eq:P-comm}, we proceed to
estimate the `error terms' $R_0$, $R_i$, $R_{ij}$ as \emph{operators}. We
start with $R_0.$ As follows from the standard square root construction to
prove the $L^2$ boundedness of pseudodifferential operators, there exists
$R_0'\in\Psib^{-1}(M)$ such that
\begin{equation*}
\|R_0v\|\leq2\sup|r_0|\,\|v\|+\|R_0'v\|
\end{equation*}
for all $v\in L^2(M,dg)$. Here $\|\cdot\|$ is the $L^2(M,dg)$-norm, as usual.
Thus, we can estimate, for any $\gamma>0$,
\begin{multline}
|\langle R_0 v,v\rangle|\leq \|R_0 v\|\,\|v\|
\leq 2\sup |r_0|\,\|v\|^2+\|R_0' v\|\,\|v\|\\
\leq 6C_2(\delta\beta+\beta^{-1})\|v\|^2
+\gamma^{-1}\|R_0' v\|^2+\gamma \|v\|^2. 
\label{30.8.2006.33}\end{multline}

Now we turn to $R_i.$ Let $T_{-1}\in\Psib^{-1}(M)$ be elliptic (we will use
this to keep track of orders), and $T_1\in\Psib^1(M)$ a parametrix, so
$T_1T_{-1}=\Id+F$, $F\in\Psib^{-\infty}(M)$.  Then there exist
$R'_i\in\Psib^{-1}(M)$ such that
\begin{multline*}
\|R_i w\|=\|R_i (T_1 T_{-1} -F)w\|\leq\|(R_i T_1)(T_{-1}w)\|+\|R_iFw\|\\
\leq 6C_2(\delta\beta+\beta^{-1})\|T_{-1}w\|+
\|R_i' T_{-1}w\|+\|R_iFw\|
\end{multline*}
for all $w$ with $T_{-1}w\in L^2(M,dg)$.
Similarly, there exist $R'_{ij}\in\Psib^{-1}(M)$ such that
\begin{equation}\label{tminusrij}
\|(T_1)^*R_{ij} w\|\leq 6C_2(\delta\beta+\beta^{-1})\|T_{-1}w\|+
\|R_{ij}' T_{-1}w\|+\|(T_1)^*R_{ij}Fw\|
\end{equation}
for all $w$ with $T_{-1}w\in L^2(M,dg)$.
Thus,
\begin{multline}
|\langle R_i Q_i v,v\rangle|\leq
6C_2(\delta\beta+\beta^{-1})\|T_{-1}Q_i v\|\,\|v\|\\
\qquad+2\gamma\|v\|^2+\gamma^{-1}\|R'_i T_{-1}Q_i v\|^2+\gamma^{-1}\|F_i
Q_iv\|^2, 
\label{30.8.2006.31}\end{multline}
and, writing $Q_jv=T_1T_{-1}Q_jv-FQ_jv$ in the right factor, and taking the
adjoint of $T_1$, \eqref{tminusrij} gives
\begin{multline}
|\langle R_{ij} Q_i v,Q_jv\rangle|\leq
6C_2(\delta\beta+\beta^{-1})\|T_{-1}Q_iv\|\,\|T_{-1}Q_jv\|\\
+2\gamma\|T_{-1}Q_jv\|^2+\gamma^{-1}\|R'_{ij}T_{-1} Q_iv\|^2
+\gamma^{-1}\|F_{ij} Q_iv\|^2+\|R_{ij}Q_iv\|\,\|FQ_jv\|, 
\label{30.8.2006.32}\end{multline}
with $F_i,$ $F_{ij}\in\Psib^{-\infty}(M).$

Finally we turn to the positive commutator argument itself.  Let
$s<\sup\{s': q_0 \notin \dWF u\};$ hence, without loss of generality, we
have $\dWF^s u \cap U =\emptyset$ where we shrink $U$ as necessary. We will
prove that $q_0 \notin \dWF^{s+1/2} u,$ a contradiction unless $s=+\infty.$

Let $\Lambda_r$ be a quantization of
\begin{equation}
|\taub|^{s}(1+r|\taub|^2)^{-s/2},\quad r\in[0,1),
\end{equation}
and set $A_r=A\Lambda_r\in\Psib^{0}(M)$ for $r>0.$ Then $A_r$ is uniformly
bounded in $\Psibc^{s}(M),$ and we may further arrange that $[\Box,
\Lambda_r]=0.$

By \eqref{eq:P-comm},
\begin{equation}\begin{split}\label{eq:pos-comm}
\langle i[A_r^*A_r,\Box]u,u\rangle
&=\|C\tilde B\Lambda_r u\|^2
+\langle R'\Box\Lambda_r u,\Lambda_r u\rangle
+\langle R_0 \tilde B\Lambda_r u,\tilde B\Lambda_r u\rangle\\
&\qquad+\sum \langle R_i Q_i \tilde B \Lambda_r u,\tilde B\Lambda_r u\rangle
+\sum \langle R_{ij} Q_i \tilde B\Lambda_r u,
Q_j\tilde B\Lambda_r u\rangle\\
&\qquad+\langle R''\Lambda_r u,\Lambda_r u\rangle
+\langle (E+E')\Lambda_r u,\Lambda_r u\rangle.
\end{split}\end{equation}
On the other hand, as $A_r\in\Psib^0(M)$ for $r>0$ and $u\in \domt,$
$A_r^*A_ru\in \domt,$ and the pairing in the following computation is
well-defined:
\begin{multline}\label{eq:comm-expansion}
\langle [A_r^*A_r,\Box]u,u\rangle
=\langle A_r^*A_r \Box u,u\rangle-\langle \Box A_r^*A_ru,u\rangle\\
=\langle A_r \Box u,A_ru\rangle-\langle A_ru,A_r \Box u\rangle
=0.
\end{multline}

The first term on the right in \eqref{eq:pos-comm} is thus bounded by the
sum of the absolute values of the others. The second term on the right in
\eqref{eq:pos-comm} vanishes, to the third we can apply
\eqref{30.8.2006.33}, to the fourth \eqref{30.8.2006.31} and
\eqref{30.8.2006.32} applies to the fifth. The penultimate term is bounded
uniformly as $r\downarrow0$ by the hypothesis on $s,$ i\@.e\@.\ the assumed
regularity of $u.$ The last term is also uniformly bounded, the $E$ part
for the same reason and the $E'$ term by elliptic regularity. Similarly, in
applying \eqref{30.8.2006.33}---\eqref{30.8.2006.32} all the terms but the
first two terms on the right in each are also uniformly bounded as
$r\downarrow0$ by the assumed property of $\dWF^s u.$ So, with a constant
$C'(\gamma)$ independent of $r$ and for some $C_3>0$ depending only on the
geometry
\begin{multline}
\|C\tilde B\Lambda_r u\|^2\leq C'(\gamma)\\
+\left(6C_2(\delta\beta+\beta^{-1})+C_3\gamma\right)
\left(\|\tilde B\Lambda_r u\|^2+
\sum_i\|T_{-1}Q_i\tilde B\Lambda_ru\|^2\right)\\
+6C_2(\delta\beta+\beta^{-1})\|\tilde B\Lambda_ru\|
\sum_i \|T_{-1}Q_i\tilde B\Lambda_ru\|
+2 \gamma\|\tilde B\Lambda_ru\|^2.
\label{eq:prop-64}\end{multline}
The remaining terms on the right may be estimated as follows: writing
$T_{-1}Q_i=Q_iT'_i+x^{-1}T_i''$ for some $T'_i,T''_i\in\Psib^{-1}(M)$
(recall that $x^{-1}=Q_0$), and using Lemma~\ref{lemma:Dirichlet-form}
where necessary, we may choose our constants so as to absorb the highest
order terms into the left-hand-side (using ellipticity of $C$ on $\WF'
(\tilde B)$ and microlocal elliptic regularity).  In particular, then, we
may pick $\beta$ sufficiently large, and $\gamma=\gamma_0>0,$ $\delta_0>0,$ so
small that for all $\delta<\delta_0$
\begin{equation*}
C_4\|\tilde B\Lambda_r u\|^2\leq C_5+C_6
\|d_M T_{-1}^2\tilde B\Lambda_r u\|^2,
\end{equation*}
with $C_4>0.$ Letting $r\to 0$ now keeps the right hand side bounded,
proving that $\|\tilde B\Lambda_r u\|$ is uniformly bounded as $r\to 0,$
hence $\tilde B\Lambda_0 u\in L^2(M,dg)$ (cf\@. the proof of
Proposition~\ref{prop:elliptic}). In view of
Lemma~\ref{lemma:Dirichlet-form} this proves that
$q_0\notin\dWF^{s+1/2}(u)$, and hence proves the first statement of the
proposition.

Finally we need check that the neighborhoods of $q_0$ which are disjoint
from $\dWF^{s}(u)$ do not shrink to $\{q_0\}$ as $s\to\infty.$ This
argument is parallel to the last paragraph of the proof of
\cite[Proposition~24.5.1]{Hormander3}. In each iterative step it is only
necessary to shrink the elliptic set of $\tilde B_s$ by an arbitrarily
small amount, which allows us to conclude that $q_0$ has a neighborhood
$U'$ such that $\dWF^{s}(u) \cap U'=\emptyset$ for all $s.$ This proves
that $q_0\notin \dWF^{\infty}(u).$
\end{proof}

We note that in the proof above the ability to control $Q_0 u$ terms in
terms of the $\norm{u}_{\domt}$ was crucial, and it is this aspect of the
proof that breaks down in fiber dimension $1.$  The necessary modifications
to the proof in this case are discussed in the following section.

To carry through similar estimates in the glancing region we need the
following technical lemma. In essence this shows that, applied to solutions
of $\Box u=0$ near $\dcG,$ $D_x$ and $x^{-1}D_{z_j}$ are not merely bounded
by $D_t$ but small compared to it. Such an estimate is natural since
\begin{equation*}
x^2 p_0|_{x=0}=x^2 \tau^2\hat p_0=\tau^2-\xi^2-|\eta|^2_h-|\zeta|^2_k
\end{equation*}
gives
\begin{equation*}
\xibs^2+|\zetabs|^2_k\leq C(x^2 |\hat p|+|x|+|1-|\etabs|_h^2|),
\end{equation*}
and $1-|\etabs|_h^2$ vanishes at $\dcG,$ so the right hand side is small
near $\dcG.$ In the remainder of this section, a $\delta$-neighborhood
will refer to a $\delta$-neighborhood with respect to the metric distance
associated to any Riemannian metric on the manifold $\Sbstar M.$ The
notation $\|u\|_{\domt_{\loc}}$ and $\|\Box u\|^2_{\domt'_\loc}$ is fixed in
Remark~\ref{rem:localize}.

\begin{lemma}\label{lemma:Dt-Dx} Suppose $u\in\domt_\loc$ and that
$K\subset\Sbstar M$ is compact with 
\begin{equation*}
K\subset\cG\setminus\mWF^{s+1/2}(\Box u)
\end{equation*}
then there exist $\delta_0>0$ and $C_0>0$ such that if
$\cA=\{A_r:\ r\in(0,1]\}$ is a bounded and basic
family in $\Psibc^s(M)$ with $\bWF'(\cA)\subset U \subset\Sbstar M,$ a
$\delta$-neighborhood of $K,$ $0<\delta<\delta_0,$ and with
$A_r\in\Psib^{s-1}(M)$ for $r\in(0,1]$ then for some
$$
G\in\Psib^{s-1/2}(M),\ \tilde G\in\Psib^{s+1/2}(M)\Mwith
\bWF'(G),\bWF'(\tilde G)\subset U
$$
and $\tilde C_0=\tilde C_0(\delta)>0$
\begin{multline*}
\|D_x A_r u\|^2+\sum_i\|x^{-1}d_Z A_r u\|_k^2\\
\leq C_0\delta\|D_t A_r u\|^2+\tilde C_0(\|u\|^2_{\domt_\loc}
+\|Gu\|^2_{\domt}+\|\Box u\|^2_{\domt'_\loc}
+\|\tilde G \Box u\|^2_{\domt'})\ \forall\ r>0.
\end{multline*}
\end{lemma}

\begin{remark}
As $K$ is compact, this is a local result in the extended base,
$[0,\ep)_x\times Y$. In particular, we may assume that $K$ is a subset of
$\Sbstar M$ over a suitable local coordinate patch in the extended base.
Moreover, we may assume that $\delta_0>0$ is so small that $D_t$ is
elliptic on $U.$
\end{remark}

\begin{proof}
The proof is again very similar to \cite{Vasy5}. By
Lemma~\ref{lemma:Dirichlet-form} we already know that
\begin{equation}\begin{split}\label{eq:mic-ell-gl-8}
\|d_X A_r u\|^2
\leq &\|D_t A_r u\|^2\\
&\ +C'_0(\|u\|^2_{\domt_\loc}+\|Gu\|^2_{\domt}
+\|\Box u\|^2_{\domt'_\loc}
+\|\tilde G \Box u\|^2_{\domt'}).
\end{split}\end{equation}
for some $C'_0>0$ and for some $G$, $\tilde G$ as in the statement of the
lemma. Thus, we only need to show that if we replace the left hand side
by $\|D_x A_ru\|^2+\|x^{-1}d_Z A_r u\|^2_k$
(i.e.\ we drop the tangential derivatives,
at least roughly speaking), the constant in front of $\|D_t A_r u\|^2$
can be made small.

A simple argument as in \cite{Vasy5},
amounting to moving the tangential derivatives to the right side
and freezing the coefficients at $x=0$ (cf.\ also the
proof of Proposition~\ref{prop:elliptic}),
reduces the problem to showing that
\begin{equation}\begin{split}\label{eq:Dt-Dx-16}
&\left|\int_M
\left((D_t^2-\sum h_{ij}(0,y)D_{y_i}D_{y_j}) A_r u\,\overline{A_r
  u}\right)\, dg \right|\\
&\qquad\leq C_2\delta \|D_t A_r u\|^2+\tilde C_2(\delta)
(\|u\|^2_{\domt_\loc}+\|Gu\|^2_{\domt}),
\end{split}\end{equation}
which we proceed to do.

The key point is that the symbol of $D_t^2-\sum h_{ij}(0,y)D_{y_i}D_{y_j}$
is small on $\bWF'(\cA)$; more precisely, it is of size $\delta$.  It is
convenient to microlocalize, i.e.\ replace this tangential part of the wave
operator by an operator $F$ microsupported near $\cG$.  So let
$\psi\in\CI(\Sbstar M)$ ($\psi$ can thus be identified with a homogeneous
degree zero function on $\Tbstar M\setminus 0$) with $\psi\equiv 1$ near
$\bWF'(\cA)$, $\supp \psi\subset U$, $|\psi|\leq 1$, and let
$F\in\Psib^0(M)$ be such that
\begin{multline}\label{eq:Dt-Dx-20}
\bWF'(F)\subset U,
\ \bWF'\left(D_t FD_t-(D_t^2-\sum h_{ij}D_{y_i}D_{y_j})\right)
\cap\bWF'(\cA)=\emptyset\\
f=\sigma_{{\bo},0}(F)=\psi(1-\sum h_{ij}\etabs_i\etabs_j),
\end{multline}
where such $\psi$ and $F$ exist, since $D_t$ is elliptic on $\bWF'(\cA).$
Now,
\begin{equation*}\begin{split}
&\left|\int_M
\left((D_t FD_t-(D_t^2-\sum h_{ij}(0,y)D_{y_i}D_{y_j})) A_r u\,\overline{A_r u}\right)\right|
\leq C_2'
\|u\|^2_{\domt_\loc}
\end{split}\end{equation*}
since $(D_t FD_t-(D_t^2-\sum h_{ij}D_{y_i}D_{y_j}))A_r$
is uniformly bounded in $\Psib^{-\infty}(M)$,
by the first line of
\eqref{eq:Dt-Dx-20}.
Moreover, $1-\sum h_{ij}\etabs_i\etabs_j$ is a $\CI$
function on
a neighborhood of $K$ in $\Sbstar M$ which
vanishes at $\cG$, so
$|1-\sum h_{ij}\zetabs_i\zetabs_j|<C_3\delta$ on
a $\delta$-neighborhood of $K$, and hence
\begin{equation*}
\sup|f|\leq C_3\delta.
\end{equation*}
Since there exists
$F'\in\Psib^{-1}(M)$ with $\bWF'(F')\subset U$
satisfying
\begin{equation*}
\|Fv\|\leq2\sup|f|\,\|v\|+\|F'v\|
\end{equation*}
for all $v\in L^2(M,dg)$, we deduce that $\|Fv\|\leq2 C_3 \delta\|v\|+\|F'v\|$
for all $v\in L^2(M,dg)$. Applying this with $v=D_t A_r u$, and estimating
$\|F'v\|$ using Lemma~\ref{lemma:WFb-mic-q},
\eqref{eq:Dt-Dx-16} follows,
which in turn completes the proof of the lemma.
\end{proof}

For propagation at the glancing points it is convenient to introduce
notation for the Hamilton vector field of the `tangential part' of the
wave operator. Thus, we let
\begin{equation}
W^\flat=\pa_t-H_{\hat h},\ \hat h(y,\etab)=\abs{\taub}^{-1}|\etab|_h^2,
\end{equation}
so $W^\flat$ is homogeneous of degree zero, hence can be regarded as
a vector field on the cosphere bundle. This vector field is
well-defined at $\cG$ as a vector field tangent to $\cG$.
Away from $\cG$ it depends on
choices, but as pointed out in a remark below, these choices do not
affect the statement of the following proposition. Moreover, this
vector field also makes sense on $\dcG$, which is a smooth manifold.

It is also useful to extend $\dot\pi$ to a neighborhood of $\pa M$
as projection to $\pa M$ followed by $\dot\pi$:
\begin{equation*}
\dot\pi^{\eo}(x,t,y,z,\xis,1,\etas,\zetas)=(t,y,1,\etas).
\end{equation*}

\begin{proposition}\label{prop:tgt-prop}
Let $u\in \domt_\loc$ and suppose
\begin{equation}
K\subset\Sbstar M\text{ is compact with }
K\subset\dcG\setminus\mWF^{\infty}(\Box u).
\end{equation}
Then there exist constants $C_0>0,$ $\delta_0>0$ such that if
$q_0=(t_0,y_0,1,\etabs_0)\in K$ (hence $|\etabs|_h=1$), to conclude that
$q_0\notin\dWF(u)$ it suffices to know that for some $0<\delta<\delta_0,$
with $C_0\delta\leq\beta<1,$ and for all
$q=(x,t,y,z,\xis,1,\etas,\zetas)\in \Sigma$
\begin{multline}\label{eq:tgt-prop-est}
|\dot\pi^{\eo}(q)-\exp(-\delta W^\flat)(q_0))|
\leq\beta\delta\Mand
|x(q)|\leq\beta\delta\\
\Longrightarrow
\pi(q)\notin\dWF(u).
\end{multline}
\end{proposition}

\begin{remark}\label{rem:tgt-est}
In the estimate \eqref{eq:tgt-prop-est}, $W^\flat$ can be replaced by any
$\CI$ vector field which reduces to $W^\flat$ at $q_0$ since flow to
distance $\delta$ along a vector field only depends on the vector field
evaluated at the initial point of the flow, up to an error $\cO(\delta^2).$
Similarly, changing the initial point of the flow by $\cO(\delta^2)$ does
not affect the endpoint up to an error $\cO(\delta^2).$
Thus, estimate \eqref{eq:tgt-prop-est} can be further
rewritten, at the cost of changing $C_0$ again, as
\begin{equation}\begin{split}\label{eq:tgt-prop-est-3}
|\exp(\delta W^\flat)(\dot\pi^{\eo}(q))-q_0|&\leq\beta\delta\Mand
|x(\exp(\delta W^\flat)(q))|\leq\beta\delta\\
&\Longrightarrow
\pi(q)\notin\dWF(u);
\end{split}\end{equation}
here we interchanged the roles of the initial and final points of the flow.
\end{remark}

\begin{proof} To prove this result it suffices to modify the proof of the
glancing propagation result of \cite{Vasy5}, with the modifications similar
to those leading to the normal propagation above. The key ingredient is
Lemma~\ref{lemma:Dt-Dx}, which allows the $Q_i$ factors, including $Q_0,$
in the commutator in Corollary~\ref{cor:Box-comm} to be shown to be small,
provided the symbol of the commutant is arranged to be supported
sufficiently close to the glancing set. Thus, the $L_0$ term of
Corollary~\ref{cor:Box-comm} dominates the commutator (there is no $\xib$
dependence of the commutant in this case near the characteristic set), with
principal symbol at $x=0$ given by $2aH_{h-\tau^2}a.$

First take a function $\omega_0\in\CI(S^*(\RR\times Y))$ which is a
sum of squares of $2l$ ($l=\dim Y$) homogeneous degree zero functions $\rho_j$:
\begin{equation*}
\omega_0=\sum_{j=1}^{2l} \rho_j^2,\ W^\flat\rho_j(q_0)=0,\ \rho_j(q_0)=0,
\end{equation*}
$d\rho_j(q_0)$, $j=1,\ldots,2l$ linearly independent at $q_0$. Since
$\dim S^*(\RR\times Y)=2l+1$,
$d\rho_j(q_0)$, $j=1,\ldots,2l$, together with $dt$
($t$ is also homogeneous of degree zero), span the cotangent space of
the $S^*(\RR\times Y)$, for dimensional reasons
(note that $W^\flat t(q_0)\neq 0$). In particular,
\begin{equation*}
|\taub^{-1}W^\flat\omega_0|\leq C_1'\omega_0^{1/2}(\omega_0^{1/2}+|t-t_0|)
\end{equation*}
Then extend $\omega_0$ to a function on $\Tbstar M$ (using the
trivialization) and set
\begin{equation}
\omega=\omega_0+x^2.
\end{equation}
Then the `naive' estimate, playing an analogous role to
\eqref{eq:omega-est-a0} in the hyperbolic region, is
\begin{equation}\begin{split}\label{eq:omega-est-a0-t}
|\taub^{-1}H_{p_0}\omega|&\leq \tilde C_1''\omega^{1/2}(\omega^{1/2}+|t-t_0|
+|\xibs|^2+|\zetabs|^2)\\
&\leq C_1''\omega^{1/2}(\omega^{1/2}+|t-t_0|+\taub^{-2}|p_0|);
\end{split}\end{equation}
here we used $p_0|_{x=0}=\taub^2-|\xib|^2-|\etab|^2_h-|\zetab|^2_k,$ which
lets us estimate
\begin{equation*}
\taub^{-2}|\xib|^2\leq C(\taub^{-2}|p_0|+|x|+\omega_0^{1/2}+|t-t_0|),
\end{equation*}
for
$1-|\zetabs|_y^2$ is homogeneous of degree zero and vanishes at $\dcG$
(recall that this last estimate
motivates Lemma~\ref{lemma:Dt-Dx}).
Note that \eqref{eq:omega-est-a0-t}
is much more precise than \eqref{eq:omega-est-a0}:
we have a factor of $\omega^{1/2}+|t-t_0|+\tau^{-2}|p_0|$ in addition to
$\omega^{1/2}$---this is crucial since we need to get the direction
of propagation right.

Finally put
\begin{equation}\label{eq:glancing-phi-def}
\phi=t-t_0+\frac{1}{\beta^2\delta}\omega,
\end{equation}
and define $a=a_0$ almost as in \eqref{eq:a_m-def}, with $-\xibs$
replaced by $t-t_0$, namely
\begin{equation}\label{eq:prop-22t}
a=\chi_0(2-\phi/\delta)\chi_1((t-t_0+\delta)/
\beta\delta+1)\chi_2(|\sigma|^2/\tau^2).
\end{equation}
The slight difference is in the argument of $\chi_1$, in order to
microlocalize more precisely in the `hypothesis region', i.e.\ where $u$ is
a priori assumed to have no wave front set. This is natural, since for the
hyperbolic points we only needed to prove that singularities cannot stay at
the boundary, while for glancing points we need to get the correct
direction of propagation.  We always assume for this argument that
$\beta<1$, so on $\supp a$ we have
\begin{equation*}
\phi\leq 2\delta\Mand t-t_0\geq-\beta\delta-\delta\geq-2\delta.
\end{equation*}
Since $\omega\geq 0$, the first of these inequalities implies that
$t-t_0\leq 2\delta$, so on $\supp a$
\begin{equation}
|t-t_0|\leq 2\delta.
\end{equation}
Hence,
\begin{equation}\label{eq:omega-delta-est-t}
\omega\leq \beta^2\delta(2\delta-(t-t_0))\leq4\delta^2\beta^2.
\end{equation}
Moreover, on $\supp d\chi_1$,
\begin{equation}\label{eq:dchi_1-supp}
t-t_0\in[-\delta-\beta\delta,-\delta],\ \omega^{1/2}\leq 2\beta\delta,
\end{equation}
so this region lies in the hypothesis region of \eqref{eq:tgt-prop-est-3}
after $\beta$ and $\delta$ are both replaced by appropriate constant
multiples.

Now, using \eqref{eq:omega-est-a0-t}, \eqref{eq:omega-delta-est-t}, and
$\taub^{-1}H_{p_0}(t-t_0)=2$, we deduce that at $p_0=0$,
\begin{equation*}\begin{split}
\taub^{-1} H_{p_0}\phi&=\taub^{-1} H_{p_0}(t-t_0)+\frac{1}{\beta^2\delta}\taub^{-1}
H_{p_0}\omega\\
&\geq 1
-\frac{1}{\beta^2\delta}C_1''\omega^{1/2}(\omega^{1/2}+|t-t_0|)\\
&\geq 1-2C_1''(\delta+\frac{\delta}{\beta})\geq
c_0/4>0
\end{split}\end{equation*}
provided that $\delta<\frac{2}{16 C_1''}$,
$\frac{\beta}{\delta}>\frac{16C_1''}{2}$, i.e.\ that $\delta$ is small, but
$\beta/\delta$ is not too small---roughly, $\beta$ can go to $0$ at most
as a multiple of $\delta$ (with an appropriate constant) as $\delta\to 0.$
Recall also that $\beta<1$, so there is an upper bound as well for $\beta$,
but this is of no significance as we let $\delta\to 0$. It is also
worth remembering that in the hyperbolic region, $\beta$ roughly played the
same role as here, but was bounded below by an absolute constant, rather
than by a suitable multiple of $\delta$, hence could not go to $0$ as
$\delta\to 0$. With this, we can proceed exactly as in the
hyperbolic region, so (recalling that $\taub>0$ on $\supp a$!)
\begin{equation*}
H_{p_0} a^2=-b_0^2+e,\ b_0=\tau^{1/2} (2\tau^{-1} H_{p_0}\phi)^{1/2}
\delta^{-1/2}(\chi_0\chi_0')^{1/2}\chi_1\chi_2,
\end{equation*}
with $e$ arising from the derivative of $\chi_1\chi_2.$ Again, $\chi_0$ stands for
$\chi_0(2-\frac{\phi}{\delta})$, etc.  In view of \eqref{eq:dchi_1-supp}
and \eqref{eq:tgt-prop-est-3} on the one hand, and the fact that $\supp
d\chi_2$ is disjoint from the characteristic set on the other, both $\supp
d\chi_1$ and $\supp d\chi_2$ are disjoint from $\dWF(u)$. Thus,
$i[A^*A,\Box]$ is positive modulo terms that are controlled a priori, so
the standard positive commutator argument gives an estimate for $Bu,$ where
$B$ has symbol $b_0.$ Replacing $a$ by $a\taub^{s+1/2},$ still gives a
positive commutator since $D_t$ commutes with $\Box,$ which now gives (with
the new $B$) $Bu\in L^2(M,dg).$ In particular $q_0\notin\dWF^{s}(u).$
\end{proof}

Now, applying arguments that go back to \cite{MR83h:35120} we find

\begin{theorem}\label{thm:prop-sing} If $u\in H^{-\infty}_{\domt,\loc}(M),$
then for all $s\in\RR\cup\{\infty\},$
\begin{equation*}
\dWF^{s}(u)\setminus\mWF^{s+1}(\Box u)\subset\dot\Sigma
\end{equation*}
is a union of maximally extended generalized broken bicharacteristics of
$\Box$ in $\dot\Sigma\setminus\mWF^{s+1}(\Box u).$
\end{theorem}

\section{Fiber dimension $1$}\label{sec:f1}

In this section we indicate the changes necessary in the previous sections
to accommodate fiber dimension $1$. Fortunately, these are quite minor, due
to the rather trivial character of $1$ dimensional Riemannian geometry. The
basic reason for treating fiber dimension $1$ separately is that $x^{-1}$
is not bounded from $\domt$ to $L^2$, and terms with $x^{-1}$ arise
throughout the previous section. Here we introduce a class of operators
that we call {\em very basic} in order to eliminate these terms---this is
possible as the metric on the fibers can be put in a rather simple
form.\footnote{Another, perhaps more natural way to proceed would be to
  consider operators that commute with the projection to fiber-constant
  functions at the boundary.}

We can assume that all fibers are circles
(disconnected fibers can be dealt with similarly). Moreover, one
may arrange that the fiber metric is $k(y)dz^2$, $dz^2$ denoting the
standard metric on the circle (corresponding to a circumference of
$2\pi$, say). Correspondingly, {\em locally in the base $Y$}
(and all our considerations are local in the base)
there is a circle action on $\pa M$,
with infinitesimal generator $\pa_z$.

\begin{definition}
Let $x^{-1}\Ve(M;Y)$ be the subspace of $x^{-1}\Ve(M)$ consisting
of vector fields $V$ such that $[V,\pa_z]\in\Ve(M)$.
\end{definition}

\begin{remark}
$x^{-1}\Ve(M;Y)$ is not a left $\CI(M)$-module, but it is a
left $\CI_Y(M)$-module, where
\begin{equation*}
\CI_Y(M)=\{f\in\CI(M):\ f|_{\pa M}
\ \text{is fiber constant}\}.
\end{equation*}
Note also that $x$ is determined by
the form of the metric up to multiplication by $a\in\CI_Y(M)$.
\end{remark}

\begin{definition}
We say that $A\in\Psibc^m(M)$ is {\em very basic} if it commutes with
$\pa_z$ at $x=0$, i.e.\ if $[A,\pa_z]\in x\Psibc^m(M)$.
\end{definition}

\begin{remark}
Very basic operators form a $\CI_Y(M)$-bimodule. Moreover, they form a ring
under composition as $[AB,\pa_z]=[A,\pa_z]B+A[B,\pa_z]$. Also, if
$A\in\Psibc^m(M)$ is very basic, $x^{-p}[A,x^p]\in\Psibc^{m-1}(M)$ is also
very basic, by the Jacobi identity.
\end{remark}

The point of this definition is:

\begin{lemma} If $A\in\Psibc^m(M)$ is very basic then
\begin{enumerate}
\item
(replaces Lemma~\ref{lemma:rearrange})
for any $Q\in x^{-1}\Ve(M;Y)$ there exist very basic $A_j\in\Psib^{m-1}(M)$ and
$Q_j\in x^{-1}\Ve(M;Y)$, $j=1,\ldots,l$, $A_0\in\Psib^m(M)$ such that
\begin{equation*}
[Q,A]=A_0+\sum A_j Q_j.
\end{equation*}
\item
(replaces Lemma~\ref{lemma:commutators})
there exist $B\in\Psibc^m(M)$, $E\in\Psibc^{m-1}(M)$ depending continuously
on $A$ such that
\begin{equation}
[x^{-1}D_z,A]=B+Ex^{-1}D_z,
\end{equation}
with $E$ very basic, and
there exist $B\in\Psibc^m(M)$, $C\in\Psibc^{m-1}(M)$ depending continuously
on $A$ such that
\begin{equation}
[D_x,A]=B+CD_x,
\end{equation}
and $C$ is very basic.
\item
(replaces Lemma~\ref{lemma:Psibc-map})
if $m=0$, $A:\domt\to\domt$ with norm bounded by a seminorm of $A$
in $\Psib^0(M)$.
\end{enumerate}
\end{lemma}

Thus, for very basic operators $x^{-1}D_z$ behaves much like $D_x$.

\begin{proof}
(1,2): As very basic ps.d.o's form a $\CI_Y(M)$ bimodule, it suffices
to check that the conclusion holds for $Q=D_x, D_{y_i}$ and $Q=x^{-1}D_z$
in some local coordinates on $Y$. For $D_{y_i}$ this is
immediate as $D_{y_i}$ is a very basic element of $\Psib^1(M)$, so
$[A,D_{y_i}]$ is a very basic element of $\Psib^m(M)$.

For $Q=x^{-1}D_z$, we compute
$[x^{-1}D_z,A]=[x^{-1},A]x(x^{-1}D_z)+x^{-1}[A,D_z]$. The second term
is in $\Psibc^m(M)$ as $A$ is very basic, while
the coefficient of $x^{-1}D_z$, $E=[x^{-1},A]x\in\Psibc^{m-1}(M)$ is
very basic in view of the previous remark, proving the first line of
(2). A similar argument
applies to $D_x=x^{-1}(xD_x)$, noting that $[A,xD_x]\in x\Psibc^m(M)$
just by virtue of $A\in\Psibc^m(M)$.

(3): As in the proof of Lemma~\ref{lemma:Psibc-map}, this reduces
to $[A,V]:\domt\to L^2(M)$ being bounded for $A\in\Psibc^0(M)$ and
$V\in x^{-1}\Ve(M;Y)$. By (1), this follows from $A_0,A_j$ being bounded
on $L^2$ while $Q_j:\domt\to L^2(M)$.
\end{proof}

Of course, we need to know that there is a plentiful supply of very basic
operators.

\begin{lemma}
Suppose that $a\in S^m(\Tbstar M)$ is such that $a|_{x=0}$ is invariant
under the circle action.  Then there exists $A\in\Psib^m(M)$ very basic
such that $\sigma(A)=a$.
\end{lemma}

\begin{proof}
Using a product decomposition of $M$ near $\pa M$, extend the circle action
to a neighborhood $O$ of $\pa M$.  Let $A'\in\Psib^m(M)$ be a quantization
of $a$, so $\sigma(A')=a$, and let $A$ be the average of $A'$ under the
circle action, i.e.\ $A''=(2\pi)^{-1}\int_{S^1}A_\theta\,d\theta$, where
$A_\theta$ is the conjugate of $A$ by pull-back by translation by $\theta$.
Due to the averaging, $A$ commutes with the infinitesimal generator of the
circle action, $\pa_z$, so $A''$ is indeed very basic. Moreover,
$\sigma(A'')|_{x=0}=a|_{x=0}$ as $a$ is invariant under the circle action,
so the symbol of $A_\theta$ at $x=0$ is also $a$. Now let $\tilde A\in
x\Psib^m(M)$ be such that $\sigma(\tilde A)=a-\sigma(A'')$; then
$A=A''+\tilde A$ has all the desired properties.
\end{proof}

\emph{With the help of this lemma, all proofs in Section~\ref{sec:elliptic} go
through, provided we use commutants that are very basic.}

\begin{corollary}
(Cf.\ Lemma~\ref{lemma:D_z-comm-fine}.)  With a trivialization in which the
fiber metric at $x=0$ is $k(y)dz^2$, suppose that $\tilde a$, $a$ are as in
Lemma~\ref{lemma:D_z-comm-fine}, where $|\zetab|^2$ stands for the lift of
the standard metric on the circle $Z$. Let $A$ be as in the previous lemma.
Then $[\pa_z,A]=B\in x\Psib^m(M)$ and $\bWF'(B)\cap\dot\Sigma =\emptyset$.
\end{corollary}

\begin{proof}
We merely need to observe that $a|_{x=0}$ is invariant under the circle
action, and apply the previous lemma to deduce $[\pa_z,A]=B\in x\Psib^m(M)$.
The wave front set statement follows from the construction of
$A$: it holds (in a uniform fashion) for each $A_\theta$.
\end{proof}

With the notation of the paragraph
preceding Lemma~\ref{lemma:D_z-comm-fine-star},
since $dg=J\nu$ with $J\in\CI_Y(M)$, as $\pa_z\log J|_{x=0}=0$,
Lemma~\ref{lemma:D_z-comm-fine-star}
can be strengthened as follows:

\begin{lemma}
With $A=A^\dagger$ as in \S\ref{sec:reflection}
\begin{equation}\begin{split}\label{eq:D_z-comm-fine-star-1}
&[W, A^*A] = B+F,\ B\in x\Psib^{2m}(M),\ F\in x\Psib^{2m-1}(M)\\
&\bWF'(B)\cap\dot\Sigma=\emptyset,\ \sigma(F)=a\{a,W\log J\},
\end{split}\end{equation}
with both $B$ and $F$
depending continuously on $A$. Here $W\log J\in x\CI(M)$ is lifted
to $\Tbstar M$ by the bundle projection.
\end{lemma}

Now Lemma~\ref{lemma:Box-form} can be strengthened to include a statement
that $\Kf_{ij}$ is very basic, i.e.\ $\Kf_{ij}\in\CI_Y(M)$.

Then Corollary~\ref{cor:Box-comm} remains valid even if the sums are
so that $Q_0=x^{-1}$ is excluded from them, and $L_{ij}$ are very basic.

Consequently, the proofs in the rest of Section~\ref{sec:reflection} go
through.

\section{Edge propagation}\label{sec:edge-propagation}

The flow along \eqref{hamvf.normform} in the boundary, i.e.\ at $x=0,$ is
explicitly solvable in $\Testar_{\pa M} M$: we have
$$
\xi=\abs{\zeta}\tan (\abs\zeta s+C)
$$ while $\abs{\zeta}$ is conserved, $(z,\zeta/\abs{\zeta})$ undergo geodesic flow
at speed $\abs\zeta$ with respect to the base metric on the boundary, and
$$
\tau=A \sec(\abs\zeta s+C),\quad \eta=B \sec(\abs\zeta s +C).
$$ Thus along maximally extended integral curves, $(z,\zeta/\abs{\zeta})$
undergoes time-$\pi$ geodesic flow and $(\tau,\eta)$ traverse a line
through the origin in $\RR^{b+1}.$ Since the Hamilton vector field along
these curves is nowhere vanishing, we conclude from the arguments of
\cite{MR95k:58168} that for any $k,l$, $\eWF^{k,l} u$ is a union of such
integral curves.  The Hamilton vector field vanishes at the endpoints of
these integral curves, at the codimension-two corner of the radial
compactification of $\Testar X,$ hence the question remaining is how the
interior wavefront set estimates on $u$ extend to the corner, and thence
into and out of $\Testar_{\pa X} X.$ In fact, we can do a (crucial) bit
more, obtaining propagation of \emph{coisotropic} regularity, i.e.\
regularity under $\coiso^k,$ along the flow discussed above.

\begin{theorem}\label{theorem:edge}
Let $u \in \eH{-\infty, l}(I \times X)$ be a distributional solution to the
wave equation with $\Box u=0,$ $\bar t \in I\subset \RR$ with $I$ open.

\begin{enumerate}
\item
Let $m>l+f/2.$ Given $p \in \dhyp_I,$ if $(\fcal_{I,p}\backslash{\pa M})
\cap \WF^m Au=\emptyset,$ for all $A \in \coiso^k$ then $p \notin
\eWF^{m,l'} B u$ for all $l'<l$ and all $B \in \coiso^k.$
\item
Let $m<l+f/2.$ Given $p \in \dhyp_O,$ if a neighborhood $U$ of $p$ in
$\Sestar|_{\pa M}M$ is such that $\eWF^{m,l}(Au)\cap U\subset\pa\fcal_O$
for all $A \in \coiso^k$ then
$p \notin \eWF^{m,l}(Bu)$ for all $B \in \coiso^k.$
\end{enumerate}
\end{theorem}

\begin{remark}
For any point $p\in \Testar_{\pa M} M\backslash \{\zeta=0\}$ there is an
element of $\coiso^k$ elliptic there, hence (1), with $k=\infty,$ shows
that solutions with (infinite order) coisotropic regularity have no
wavefront set in $\Testar_{\pa M} M\backslash\{\zeta=0\}.$ Indeed this
result holds microlocally in the edge cotangent bundle.  Note that the set
$\{x=0,\ \zeta=0\}$ is just the set of radial points for the Hamilton
vector field.
\end{remark}

As the proof follows quite closely that of Theorem~8.1 of \cite{Melrose-Wunsch1},
we shall be somewhat concise here. However, as the flowout cannot be put
in a simple model form as in \cite{Melrose-Wunsch1}, we need a lemma adapted
from \cite{Hassell-Melrose-Vasy2} to handle factors from $\cA^m,$ replacing
the powers of the fiber Laplacian used in \cite{Melrose-Wunsch1}. This
is stated in the Appendix, in Proposition~\ref{prop:HMV1-mod}.

\begin{proof}[Proof of Theorem~\ref{theorem:edge}.]
We only consider part (1) of the Theorem; the proof of (2) is completely
analogous, cf.\ \cite{Melrose-Wunsch1}.

Let $p=(\bar t,\bar y,\bar z,\bar \tau,\bar \xi,\bar \eta) \in \dhyp_I.$ We
begin by constructing a localizer in the fast variables.  Let
$$
\Upsilon: \Testar_{\pa M} (M)\backslash \{\xi=0\} \to Z
$$
be locally defined by
$$ \Upsilon (q)=z(\exp_{z_0,\zeta_0} s_\infty H_q), \quad s_\infty =
\frac{\sgn \xi}{\abs{\zeta}_{K_q}} \arctan\frac{\abs{\zeta}}{\xi}.
$$ where $q \in \Testar_{\pa M} (M)$ has coordinates
$(t_0,y_0,z_0,\tau_0,\xi_0,\eta_0,\zeta_0).$ This map simply takes a
point over the boundary to its limit point in the fiber variables along the
forward bicharacteristic flow, hence on the boundary, we certainly have
$\Upsilon_* (H) = 0.$ We now extend $\Upsilon$ smoothly to the interior of
$\Testar M,$ thus obtaining a map satisfying
$$
\Upsilon_*(H)=O(x).
$$

Let $a_i$ be homogeneous degree zero defining functions for $\fcalbar_I$ as
in Section~\ref{sec:coiso-reg}. (In fact, for the proof of (1), one could
use $\zeta_i$ in place of $a_i$, as one is only concerned about positivity
at $x=0$, where one can take $a_i=\zeta_i$,
but for the proof of (2), the use of $a_i$ is important.)

Now fix any sufficiently small neighborhood $U$ of $p$ in $\Testar M.$
There exists $\delta>0$ such that
$$
\{ x<\delta,\abs{a_i}<\delta, \abs{\eta/\tau-\bar
  \eta/\bar\tau}^2+ \abs{y-\bar y}^2+ d(\Upsilon, \bar
z)^2+ \abs{t-\bar t}^2<\delta\}\subset U.
$$
There exists $\beta$ such that $$H\left(\abs{\eta/\tau-\bar \eta/\bar\tau}^2+
  \abs{y-\bar y}^2+ d(\Upsilon, \bar z)^2+ \abs{t-\bar t}^2-\beta x\right)>0$$ on
  $U,$ since the derivative falling on each term in the cut-off is $O(x),$
  and $H(x) = \xi x +O(x^2),$ i.e.\ is bounded above on $U$ by a negative
  multiple of $x.$

We now need some cut-off functions.  Let $\psi$ be nonnegative,
nonincreasing, and supported in $(-\infty,\ep),$ and $\chi \in \CI(\RR)$ be
nonnegative, nondecreasing, and supported in $(\ep,\infty).$ We may further
arrange that $\psi$, $\chi$,
$-\psi'$ and $\chi'$ are squares of smooth functions.
We also let $\phi\in\CI_c([0,\ep))$ be nonnegative,
nonincreasing, identically $1$
near $0$, with $\phi$, $-\phi'$ the squares of a smooth function.

Set
\begin{multline}\label{symboldefn}
a_{m,l}^2 = \chi(\pm \tau)\chi(\pm\xi) \phi (x) 
\phi(\ep^{-1}\sum a_j^2) 
\psi(p^2(\cdot)/\tau^4) \\ \psi \big((\abs{\eta/\tau-\bar
  \eta/\bar\tau}^2+ \abs{y-\bar y}^2+ d(\Upsilon, \bar
z)^2+ \abs{t-\bar t}^2 -\beta x\big)
(\pm \tau)^m x^l,
\end{multline}
with $\pm = \sgn\bar\xi \bar\tau,$

Note that applying $H$ to $\tau^m x^l$ on the support of $a_{m,l}$ gives a
main term $(m+l)\xi^m x^l,$ hence this term is positive provided $m+l>0.$
As long as $\ep$ is sufficiently small that $(\beta+1)\ep<\delta,$ the
support of the symbol lies in $U,$ hence choosing $\ep$ sufficiently small
ensures the positivity of the term in $H (a_{m,l}^2)$ coming from the last
cut-off term. The cut-off factors involving $\chi$ are harmless as the
support of the Hamilton vector field applied to these intersects the
characteristic set inside $\Testar M$ in a compact set. Finally,
the cut-off terms involving $\phi$ have the following properties:
the Hamilton derivative of $\phi(x)$ supported
in the interior of $\Testar M$, while the Hamilton derivative of
$\phi(\ep^{-1}\sum a_j^2)$ is supported away from $\fcalbar_I$.

Using the above observations, provided $m+l>0,$ by choosing $\ep$
sufficiently small, we can ensure that
$$
H (a_{m,l}^2)= \pm (a')^2\pm\sum_j b_j^2+e+c+ k
$$ where $\supp e\subset T^* X^\circ,$ $\supp c$ is compact in
$\Testar(M)$ and $\Sigma\cap\supp(k)=\emptyset.$

For $k=0$, the proof is finished with a positive commutator argument as in
\cite[Theorem~8.1]{Melrose-Wunsch1}.  For arbitrary $k$, we can now use
Proposition~\ref{prop:HMV1-mod} to finish the proof inductively.
\end{proof}

\begin{remark}

The reason for the appearance of $f/2$ in the statement of the theorem is
that in order to prove the theorem, we need to use $a_{r,s}$ with
$r=m-\frac{1}{2}$, $s=-l-\frac{f-1}{2}$. Indeed, the commutator of
$A_{r,s}^*A_{r,s}$ (with $A_{r,s}\in\ePs{r,s}(M)$ having principal symbol
$a_{r,s}$) with $\Box$ is in $\ePs{2r-1,2s-2}(M)$. As in the commutator
argument we use the metric density $dg$ (to make $\Box$ formally
self-adjoint), which is $x^{f+1}$ times a non-degenerate smooth b-density
(with respect to which $\eH{m,l}(M)$ is weighted), a principal term in this
commutator of the form $(A')^*A'$ provides a bound for the $\eH{m,l}(M)$
norm of $u$ provided that $2m-(2r-1)=0$, $2l+(2s-2)+(f+1)=0$, leading to
the stated values of $r$ and $s$ in terms of $m$ and $l$.  As in order to
obtain a positive commutator we need $r+s>0$ (or $<0$ in the outgoing
region), we deduce that $m>l+\frac{f}{2}$ (resp. $m<l+\frac{f}{2}$) must be
satisfied.

Moreover, the reason for $l'<l$ in the statement of the theorem is the
possible need for an interpolation argument, in case the a priori regularity
of $u$ is weak, i.e.\ $u\in\eH{q,l}$ with $q$ too small (possibly negative).
This is because, due to the standard error terms, the regularity of $u$
can only be improved by $1/2$ order at a time, so already $m=q+1/2$ would
have to satisfy $m>l+\frac{f}{2}$, i.e.\ we would need $q>l+\frac{f-1}{2}$.
In particular, if  $q>l+\frac{f-1}{2}$ is satisfied, i.e.\ if we have this
much a priori regularity, we can take $l'=l$.
\end{remark}

\section{Propagation of coisotropic regularity}\label{section:coisotropic}

\begin{theorem}\label{theorem:coisotropic}
Let $p\in \dhyp$, $\ep>0$ and $k \in \NN.$ Then there is $k'$ (depending on $k$
and $\ep$) such that if $\Box u=0$ and $u$ has coisotropic regularity of
order $k'$ relative to $H^s$ (on the coisotropic $\fcaldot_I$)
near $\fcaldot_{I,p}$ strictly away from
$\pa M,$ then $u$ has coisotropic regularity of order $k$ relative to
$H^{s-\ep}$ (on the coisotropic $\fcaldot_O$)
near $\fcaldot_{O, p},$ strictly away from $\pa M.$

In particular, if $\Box u=0$ and $u$ is coisotropic (i.e.\ has infinite
regularity) relative to $H^s$ (on the coisotropic $\fcaldot_I$)
near $\fcaldot_{I,p}$ strictly away from
$\pa M,$ then for all $\ep>0$, $u$ is coisotropic relative to
$H^{s-\ep}$ (on the coisotropic $\fcaldot_O$)
near $\fcaldot_{O, p},$ strictly away from $\pa M.$
\end{theorem}

Explicitly, with the argument presented below, one can take $k'=k$ if
$\ep>1/2$, and any $k'>\frac{1}{2\ep}\,k$ if $\ep\leq 1/2$.

\begin{proof}
First we show that the
result holds with $H^{s-0}$ replaced by $H^{s-1/2-0}.$ Indeed, by
Theorem~\ref{theorem:diffractive}, $u$ is $\CI$ microlocally near
$\fcaldot_{O,p}$ strictly away from $\pa M$ provided it is $\CI$
microlocally near $\fcaldot_{I,p}$ away from $\pa M$, i.e.
\begin{equation*}
\WF(u)\cap(
\fcaldot_{I,p}\setminus\pa\fcaldot_{I,p})=\emptyset
\Longrightarrow
\WF(u)\cap (\fcaldot_{O,p}\setminus\pa\fcaldot_{O,p})=\emptyset.
\end{equation*}
In particular, $u$ is coisotropic as stated, provided that $\WF(u)\cap(
\fcaldot_{I,p}\setminus\pa\fcaldot_{I,p})=\emptyset$. Thus, we we may
assume that $u$ is microlocalized near $\fcaldot_{I,p}$ and away from $\pa
M$ for $t<t(p)$, hence that $u$ and $\coiso^k u$ lie in
$H^s_{\loc}(M^\circ)$ for $t<t(p)$ (with $\loc$ indicating locally in
time). It is convenient to normalize $s$ to make
Proposition~\ref{prop:domains} easier to use.  Let $\Theta_s$ have Schwartz
kernel
\begin{equation}
\kappa (\Theta_s)(t,t') = \psi(t-t') \kappa(\abs{D_t}^s)(t,t')
\label{Theta}
\end{equation}
where $\psi (t)$ is a smooth function of compact support, equal to one near
$t=0$ (cf.~\cite{Melrose-Wunsch1}).  Since applying $\Theta_{-s}$ to $u$
preserves the hypotheses, except replacing $s$ by $0$, and then the
conclusion is preserved (shifting $s$ back) upon the application of
$\Theta_{s}$, we may  apply Proposition~\ref{prop:domains},
to conclude that $u\in\eHloc{s,s-\frac{f+1}{2}}(M)$ for
$t<t(p)$, hence for all $t$.

Then Theorem~\ref{theorem:edge} part (1)
can be applied, first with $l=s-\frac{f+1}{2}$,
$m=s-\frac{1}{2}+0$ near incoming points (so $m>l+\frac{f}{2}$),
to conclude that
$u$ is coisotropic of order $k$ relative to
$\eH{s-\frac{1}{2}+0,s-\frac{f+1}{2}-0}(M)$ near $\pa\fcaldot_I$.
Note that if $U$ is a sufficiently small neighborhood of
$\dhyp_I=\pa\fcaldot_I$ in
$\Testar_{\pa M}M$, this implies that in fact
$u$ is microlocally in $\eH{s+k-\frac{1}{2}+0,s-\frac{f+1}{2}-0}(M)$
on $U\setminus\dhyp_I$,
since $\module$ has elliptic elements near each point
in $U\setminus\dhyp_I$. Since bicharacteristics $\gamma$ in $\dot\pi^{-1}(\dhyp)
\setminus (\dhyp_I\cup\dhyp_O)\subset\Testar_{\pa M}M$ tend to $\dhyp_I,$ resp.
$\dhyp_O$, as the parameter along $\gamma$ tends to $\pm\infty$,
the standard (non-radial) propagation
of singularities,
\cite[Theorem~8.1(ii)]{Melrose-Wunsch1},
yields that $u$ is coisotropic of order $k$ relative to
$\eH{s-\frac{1}{2}+0,s-\frac{f+1}{2}-0}(M)$ at $\dot\pi^{-1}(\dhyp)
\setminus\dhyp_O$---indeed, $u$ is simply {\em in}
$\eH{s+k-\frac{1}{2}+0,s-\frac{f+1}{2}-0}(M)$ microlocally
on $\dot\pi^{-1}(\dhyp)
\setminus (\dhyp_I\cup\dhyp_O)$.
Applying Theorem~\ref{theorem:edge} with $m=s-\frac{1}{2}-0$,
$l=s-\frac{f+1}{2}-0$ (the two small constants in the $-0$ are
chosen so that $m<l+\frac{f}{2}$), we conclude that 
$u$ is coisotropic of order $k$ relative to
$\eH{s-\frac{1}{2}-0,s-\frac{f+1}{2}-0}(M)$ near $\fcaldot_O$, as claimed.

On the other hand,
$u$ is in $H^s$ along $\fcaldot_{O,p}$ by
Theorem~\ref{theorem:diffractive}.  Hence the theorem follows by the
interpolation result of the following lemma.
\end{proof}

\begin{lemma}\label{lemma:coisotropic-interpolate}
Suppose that $u$ is in $H^s$ microlocally near some point $q$ away from
$\pa M$, and it is coisotropic of order $N$ relative to $H^m$ near $q$
with $s>m$. Then for $\ep>0$ and $k<\frac{\ep N}{s-m}$,
$u$ is coisotropic of order $k$ relative to
$H^{s-\ep}$ near $q$.

In particular, if $u$ is in $H^s$ microlocally near some point $q$ away from
$\pa M$ and $u$ is coisotropic (of order $\infty$, that is)
relative to $H^m$ near $q$ with $s>m$, then $u$ is coisotropic
relative to $H^{s-\ep}$ for all $\ep>0$.
\end{lemma}

\begin{proof}
If $Q\in\Psi^0(M)$ and $\WF'(Q)$ lies sufficiently close to $q$, then the
hypotheses are globally satisfied by $u'=Qu$. Moreover, being coisotropic,
locally $\fcal$ can be put in a model form $\zeta=0$ by a symplectomorphism
$\Phi$ in some canonical coordinates $(y,z,\eta,\zeta)$, see
\cite[Theorem~21.2.4]{Hormander3} (for coisotropic submanifolds one has
$k=n-l$, $\dim S=2n$, in the theorem). Further reducing $\WF'(Q)$ if
needed, and using an elliptic $0$th order Fourier integral operator $F$
with canonical relation given by $\Phi$ to consider the induced problem for
$v=Fu'=FQu$, we may thus assume that $v\in H^s$, and $D_z^\alpha v\in H^m$
for all $\alpha,$ i.e.\ $\langle D_z\rangle^N v\in H^m.$ Considering the
Fourier transform $\hat v$ of $v$, we then have $\langle\eta,\zeta\rangle^s
\hat v\in L^2$, $\langle\eta,\zeta\rangle^m\langle\zeta\rangle^N \hat v\in
L^2.$ But this implies
$\langle\eta,\zeta\rangle^{m\theta+s(1-\theta)}\langle\zeta\rangle^{N\theta}
\hat v\in L^2$ for all $\theta\in[0,1]$ by interpolation (indeed, in this
case by H\"older's inequality).  In particular, taking
$\theta=\frac{\ep}{s-m}$,
$\langle\eta,\zeta\rangle^{s-\ep}\langle\zeta\rangle^k \hat v\in L^2$ if
$k<\frac{N\ep}{s-m}$, and the lemma follows.
\end{proof}

\section{Geometric theorem}\label{section:geometric}

We now prove Theorem~\ref{theorem:2}, using as our main ingredients
Theorem~\ref{theorem:coisotropic} and a duality argument.

Let $U_1$ be a small open neighborhood of a single point $w \in
\fcal^\circ_{O,p}$ such that all points in $\fcal_O\cap U_1$ have a
distance from $\pa M$ between $0.9 d_0$ and $d_0$ for some small $d_0.$
Pick any time $T$ greater than $d_0,$ so that all interactions of $U_1$
with the boundary occur under backward generalized bicharacteristic flow
for time \emph{less than $T.$}

Let $U_0^G$ and $U_0^D$ now denote two open sets in $T^*M^\circ$ such that
$U_0^D$ contains the whole time-$T$ backward flowout of $U_1$ in the sense
of generalized broken bicharacteristics.  $U_0^G,$ by contrast, must
contain only the \emph{geometric} backward flow from $U_1$ (see
\S\ref{section:edgeb} for the definitions of these flows).  By hypothesis,
we may choose $U_0^D$ such that the nonfocusing condition holds relative to
$H^s$ on $U_0^D.$ We may split the initial data into a piece microsupported
in $U_0^G$ and supported away from a boundary, and a remainder.  The former
piece, by hypothesis, is globally in $\mathcal{C}(\RR; \dcal_s),$ hence
satisfies the conclusion of the theorem.  Thus it suffices to consider only
the latter piece, which has no wavefront set \emph{of any order} in
$U_0^G.$

Let $r<s.$ Let $A_i,$ $i=1,\dots,N,$ denote first-order pseudodifferential
operators, generating $\module$ as above, but now locally over a large set
in $M^\circ,$ and with kernels compactly supported in $M^\circ.$ Let $B_0$
and $B_1$ be pseudodifferential operators of order $0,$ compactly supported
in $M^\circ,$ and with microsupport in $U_0^D$ resp.\ $U_1,$ such that
$B_1$ is elliptic at $w$ and such that $B_0$ is elliptic on the
time-$T$ backward flowout of $\WF' B_1$ along generalized broken
bicharacteristics.  Let $U$ denote the forward time-evolution operator,
taking wave equation solutions on $[0,\delta]$ to wave equation solutions
on $[T-\delta, T].$ Theorem~\ref{theorem:coisotropic}, together with
time reversal symmetry, gives
\begin{align*}
&\sum_{\abs{\alpha}\leq k'} \norm{\Theta_{-r+1} A^\alpha B_1 U (u)}_{\domt([T-\delta,T])}^2 < \infty\\ \Longrightarrow
& \sum_{\abs{\alpha}\leq k} \norm{\Theta_{-r+1-\ep} A^\alpha B_0 u}_{\domt([0,\delta])}^2 < \infty,
\end{align*}
with $\Theta_{\cdot}$ defined by \eqref{Theta}.  Thus, if we define Hilbert
spaces $$ \hilbert_0 = \big\{ u \in \domt([0,\delta]): \Box u=0,\
\sum_{\abs{\alpha}\leq k} \norm{\Theta_{-r+1-\ep} A^\alpha B_0
u}_{\domt([0,\delta])}^2 < \infty \big\}
$$
and
$$
\hilbert_1 = \big\{ u \in \domt([T-\delta,T]): \Box u=0,
\ \sum_{\abs{\alpha}\leq k'} \norm{\Theta_{-r+1} A^\alpha B_1 u}_{\domt([T-\delta,T])}^2 < \infty \big\},
$$ we have $U^{-1}: \hilbert_1 \to \hilbert_0.$ Hence by unitarity of $U$ with
respect to the norm on $\domt,$ $U: \hilbert_0^* \to \hilbert_1^*,$ \emph{where
the dual spaces are taken with respect to energy norm, i.e.\ the norm on
$\domt$}.  Now
$$
\hilbert_0^* = \big\{ u \in \domt([0,\delta]): \Box u=0,\
u=\sum_{\abs{\alpha}\leq k} (A^\alpha)^*_\dcal (B_0)^*_\dcal \Theta_{-r+1-0} v_\alpha,\ v_\alpha \in \domt \big\}
$$
and
$$
\hilbert_1^* = \big\{ u \in \domt([T-\delta,T]): \Box u=0,\
u=\sum_{\abs{\alpha}\leq k'} (A^\alpha)^*_\dcal (B_1)^*_\dcal \Theta_{-r+1} v_\alpha,\ v_\alpha \in \domt \big\}
$$ Thus since $u$ satisfies the nonfocusing condition of degree $k$ w.r.t.\
$H^{s}$ for $t \in [0,\delta]$ on $U_0^D,$ it lies in $\hilbert_0^*;$ hence
$U(u) \in \hilbert_1^*,$ i.e.\ $u$ also satisfies the nonfocusing condition
of degree $k'$ w.r.t.\ $H^{r}$ for $t\in [T-\delta,T]$ on $U_1.$ (Here we
have used Lemma~\ref{lemma:adjoints} to see that lying in the range of
$(A^\alpha)^*_\dcal$ yields nonfocusing; recall that the $\dcal$ subscript
means adjoint with respect to the $\dcal$ inner product.)  So microlocally
near $w$ and for $t\in [T-\delta,T],$
\begin{equation}\label{nonfoc-reg}
u \in \sum_{\abs{\alpha}\leq k'} A^\alpha (H^r).
\end{equation}
On the other hand, since $\WF u \cap U_0^G=\emptyset,$
Theorem~\ref{theorem:edge} shows that there exists $\beta \in \RR$ such
that microlocally near $w$ and for
$t\in [T-\delta,T],$
\begin{equation}\label{coiso-reg}
A^{\alpha} u \in H^\beta\quad \forall \alpha
\end{equation}
(The particular choice of $\beta$ is dependent on the background regularity
of the solution.)  We can interpolate \eqref{nonfoc-reg} and
\eqref{coiso-reg} similarly to Lemma~\ref{lemma:coisotropic-interpolate};
the only difference is that one has
$\langle\eta,\zeta\rangle^r\langle\zeta\rangle^{-k'} \hat v\in L^2$ rather
than $\langle\eta,\zeta\rangle^s\hat v\in L^2$, with the notation of that
lemma. This shows that near microlocally $w,$ for $t\in [T-\delta, T],$ $u
\in H^{r-0}.$ Since $r<s$ is arbitrary, this proves
Theorem~\ref{theorem:2}.\qed

\section{Applications to Lagrangian data}

An important class of examples of solutions satisfying the nonfocusing
condition is given by choosing as Cauchy data Lagrangian distributions in
$X^\circ,$ with respect to Lagrangian manifolds transverse to the
coisotropic submanifold $\F$ obtained by flowout from the boundary.
For the proof of the proposition below it is useful to have the following
lemma which puts the coisotropic and Lagrangian manifolds into model form,
and which seems to be well-known although it is hard to find a published
reference.

\begin{lemma}
Suppose that $\F$ is coisotropic of codimension $k$, $\lag$ is Lagrangian,
and $\F$ and $\lag$ are transverse. Let $\omega$ denote the standard
symplectic form on $T^*\RR^n$, and write $\RR^n=\RR^{n-k}_y\times\RR^k_z$,
with dual coordinates $(\eta,\zeta)$.  Then there is a local
symplectomorphism mapping $\F$, resp.\ $\lag$ to $\{\zeta=0\}$, resp.\
$\{y=0,z=0\}$ in $(T^*\RR^n,\omega)$.

If in addition $\F$ and $\lag$ are \emph{conic} and the canonical one-form
does not vanish on $T_q \F,$ with $q \in \F \cap \lag$ then there is a
\emph{homogeneous} symplectomorphism mapping a neighborhood of $q$ into
$(T^*\RR^n,\omega)$, equipped with the standard (fiberwise) $\RR^+$-action,
such that $\F$, resp.\ $\lag$ are mapped into $\{\zeta=0\}$, resp.\
$\{y=0,z=0\}$.
\end{lemma}

\begin{proof}
By \cite[Theorem~21.2.4]{Hormander3} (and the subsequent remark),
we may assume that $\F$ is
given by $\{\zeta=0\}$ in $(T^*\RR^n,\omega)$. As $\lag$ is transverse
to $\F$, $d\zeta_j$, $j=1,\ldots,k$,
restrict to be linearly independent on $\lag$ at points
in $\lag\cap\F$. Moreover, for $q\in\lag\cap\F$, we may assume
that $d\eta_j(q)$, $j=1,\ldots,n-k$, together with $d\zeta_j(q)$, are
linearly independent on $\lag$. Indeed, a linear combination
$\sum c_j\pa_{z_j}$ cannot be tangent
to $\lag$ at $q$, as it is the image of $\sum c_j\,d\zeta_j$
under the Hamilton map;
for a Lagrangian submanifold the tangent space
is the image of the conormal bundle
under the Hamilton map, and we just established that $\sum c_j\,d\zeta_j(q)$ is
not an element of the conormal bundle of $\lag$. By a symplectic
linear transformation in $T^*\RR^{n-k}$ that switches the roles of
some components of $dy$ and $d\eta$ we can arrive at the stated situation.

Thus, $\lag$ is locally a graph over
$\{0\}\times(\RR^n)^*_{\eta,\zeta}$, i.e.\ on $\lag$, $y=Y(\eta,\zeta)$,
$z=Z(\eta,\zeta)$.
As $\lag$ is Lagrangian with respect
to $\omega=\sum d\eta_j\wedge dy_j+\sum d\zeta_j\wedge dz_j$, it
follows that $\lag$ is locally the graph of the differential of a function
$F:(\RR^n)^*_{\eta,\zeta}\to\RR$, i.e.\ $Y=d_\eta F$, $Z=d_\zeta F$.
The map
\begin{equation}\label{eq:move-lag}
(y,z,\eta,\zeta)\mapsto (y-d_\eta F,z-d_\zeta F,\eta,\zeta)
\end{equation}
is a local symplectomorphism, it preserves $\F=\{\zeta=0\}$, and
maps $\lag$ to $\{y=0, z=0\}$ as desired.

For the conic version, one may apply \cite[Theorem~21.2.4]{Hormander3},
so one can assume that $\F$ is given by $\zeta=0$. We may
assume that $y(q)=0$; $\eta(q)\neq 0$ since $q$ does not lie in the
zero section. Now,
the span of the $\pa_{z_j}$ intersects $T_q\lag$ trivially as above.
Let $V^*$ be the span of the $\sum c_jd\eta_j(q)$ which lie in $N^*_q\lag$,
i.e.\ for which the image $V$ under the Hamilton map
$\sum c_j\pa_{y_j}\in T_q\lag$.
Taking $W$ to be a complementary subspace to $V$ in the span of the $dy_j$,
we see that functions with differentials in $V$ plus functions with
differentials in $W^*$, together with the $\zeta_j$, give coordinates
on $\lag$ locally. Moreover, by the conic hypothesis, $W$ is at least
one-dimensional. There is a homogeneous symplectomorphism
switching the roles of the variables in $V$ and $V^*$ (if $V$ is not
already trivial); if $\eta''_1$ is away from $0$, as one may assume,
one can take the map on $V\oplus W\oplus V^*\oplus W^*$ given by
\begin{equation*}
(y',y_1'',\eta',\eta_1'')\mapsto(-\frac{\eta'}{\eta''_1},y_1''+\frac{y'\cdot\eta'}{\eta''_1},\eta''_1
y',\eta_1''),
\end{equation*}
with the rest of the $y''$ and $\eta''$ as well as the $z$ and $\zeta$ are
unchanged under the map. Note that the pullback of the differentials of
the new $\eta$ coordinates is $\eta''_1\,dy'_j+y'_j\,d\eta''_1$ and
$d\eta''_j$, which have the same span as $dy'_j$ and $d\eta''_j$, hence
pull-back to be linearly independent in $T^*_q\lag$.  Thus, one can arrange
that $\lag$ is locally a graph over $\{0\}\times(\RR^n)^*_{\eta,\zeta}$,
which in addition is conic, so $Y$ and $Z$ as above are homogeneous of
degree zero. Then $F$ may be arranged to be homogeneous of degree $1$ (we can
take $F=\sum \eta_j Y_j+ \sum \zeta_j Z_j$---see \cite[Proof of
Theorem~21.2.16]{Hormander3}), so \eqref{eq:move-lag} is homogeneous,
proving the lemma.

\end{proof}

\begin{proposition}
Let $p\in \fcal_I$ and let
$q\in \fcal^X_I$ be its projection.
Let $u$ be a solution to the wave equation with Cauchy
data that is in $\dcal_\infty$ near $\pa X$ and that in $X^\circ$ is given
by $(u, D_t u) \restrictedto_{t=0} =(u_0, u_1)$ with $(u_0, u_1) \in
(I^s(\lag), I^{s-1}(\lag))$ where $\lag$ is a Lagrangian intersecting
$\F^X$ transversely at $\fcaldot^X_{q}.$ Then $u$ satisfies the nonfocusing
condition of degree $(-s-n/4)+f/2-,$ microlocally along $\fcaldot_{I,p}.$
\end{proposition}
Note that $-s-n/4-0$ is the a priori Sobolev regularity of the initial data,
hence the gain in this result is of $f/2-0$ derivatives.

\begin{proof}
Without loss of generality, we may take $u_1=0$ (if instead, $u_0=0$ we
simply consider the solution $\pa_t u$). Decomposing by a partition
of unity, it suffices to assume that $u_0$ has support in a small open set
in $X^\circ.$

By the preceding lemma, there exists a local symplectomorphism $\Phi$ that
reduces this geometric configuration to a normal form: $\Phi$ maps $\F^X$
to the set $\{\zeta=0\}$ and $\lag$ to $N^*\{p\}$ with $p \in X^\circ.$
Quantize $\Phi$ to an FIO $T$ of order $0,$ with parametrix $S$, so that
$ST-I,\ TS-I \in \Psi^{-\infty}(M^\circ).$
Thus if $Y \in \coiso^m,$
$$
Tu = (T Y S) Tv \Longrightarrow u = Y v +r,\quad r \in \CI(M^\circ),
$$ and this would give the desired conclusion.  As we are working locally
in $M^\circ$, $T$ and $S$ preserve the scale of spaces $\dcal_r;$ thus it
now suffices to show that there exists $w\in \dcal_{-s-n/4+f/2-0}$ such
that $(T (1+ Y)^m S) w = Tu,$ with $Y \in \coiso^2.$ By Egorov's theorem,
$T Y S$ is a pseudodifferential operator of order $2$ whose symbol vanishes
quadratically on $\{\zeta=0\};$ it thus suffices to take $TYS=\Lap_z,$
where $\Lap_z$ is the (indeed, any) Laplacian in the fibers.

Meanwhile, $Tu\in I^s(N^*\{p\}).$ We thus simply drop $T$ and take $u\in
I^s(N^*\{p\}),$ $\F^X = \{\zeta=0\},$ $Y=\Lap_z.$  We lump together the base variable
$(x,y)$ to a single set of variables $\tilde y,$ and shift our variables so
that $\tilde y(p)=z(p)=0.$ Thus, in our transformed coordinates we have
$$
(1+Y)^{-N} u = \int e^{i (\tilde y\cdot \eta + z'\cdot \zeta'+
  (z-z')\cdot \zeta)}a(z,\zeta) b(x,\tilde y,z',\xi,\eta,\zeta')\, dz' \,
d\xi\, d\eta \, d\zeta\, d\zeta',
$$ where $a,$ the symbol of $(1+\Lap_z)^{-N},$ is a symbol of order
$-2N,$ and where $b$ is a symbol of order $s-n/4.$

For $N$ sufficiently large, the integral in $\zeta$ is absolutely
convergent, taking values in Lagrangian distributions of order $s+f/4$ in
the $\tilde y$ variables, hence
$$ (1+Y)^{-N} u \in \mathcal{C}(z; H_\loc^{-s-n/4+f/2-0})\subset
H_\loc^{-s-n/4+f/2-0}(X^\circ).
$$ Thus $u \in \dcal_{-s-n/4+f/2-0}+ Y (\dcal_{-s-n/4+f/2-0})+\dots +
Y^{N} (\dcal_{-s-n/4+f/2-0}).$
\end{proof}

\appendix
\renewcommand{\theequation}{A.\arabic{equation}}
\section*{Appendix: Iterative regularity under modules}
\renewcommand{\thesection}{A}

In this section we adapt the result of \cite{MR2020655} and
\cite{Hassell-Melrose-Vasy2} to the present setting. The key new ingredient
is that the set of radial points is not discrete, so we cannot be quite as
specific with microlocalization as in these papers; we use the construction
of Section~\ref{sec:edge-propagation}.  The notation below is that of
Definition~\ref{definition:coisotropic} and Lemma~\ref{lemma:test-module};
additionally, we use the `reduced' multiindex notation, dropping $A_N$,
which is an elliptic multiple of $\Box$ and thus treated separately, and
set
\begin{equation}
A_{\alpha}= \prod_{i=1}^{N-1}A_i^{\alpha_i},\ \alpha_i\in\bbN_0,\ 
1 \leq i \leq N-1.
\label{eq:red-multi}\end{equation}
Also let $W_{m,l}\in\ePs{m,l}(M)$ with $\sigma(W_{m,l})=|\tau|^m x^l$, and
let
\begin{equation*}
A_{\alpha,m,l}=W_{m,l}A_\alpha.
\end{equation*}
The justification for treating $A_N$ separately is:

\begin{lemma}(See \cite[Corollary~6.4]{MR2020655}.)
\label{lemma:reduced-basis}
Suppose $\module$ is a test module (see \cite[Definition~6.1]{MR2020655}),
and $A_N$ is a generator with principal symbol $a_N$ an elliptic multiple
of $\Box$. Suppose $u\in\CmI$ satisfies $\Box u=0$ and $u$ is coisotropic
of order $k-1$ on $O$ relative to $\eH{m,l}(M)$.  Then for $O'\subset O,$
$u$ is coisotropic of order $k$ on $O'$ relative to $\eH{m,l}(M)$ if for
each multiindex $\alpha,$ with $|\alpha|=k,$ $\alpha_N= 0$, there exists
$Q_\alpha\in\ePs{0,0}(M),$ elliptic on $O'$ such that $Q_\alpha
A_{\alpha,m,l}u\in L^2(M).$
\end{lemma}

\begin{proof}
By Lemma~\ref{lemma:coiso-gens} we only need to show that under our
hypotheses, for $\alpha$ with $\alpha_N\neq 0$, there exists
$Q_\alpha\in\ePs{0,0}(M)$ elliptic on $O'$ such that $Q_\alpha
A_{\alpha,m,l}u\in L^2(M)$. But such $A_{\alpha,m,l}$ is of the form
$Q_\alpha A_{\beta,m,l}A_Nu$ with $\beta_j=\alpha_j$ for $j\neq N$,
$\beta_N=\alpha_N-1$, and $A_N= W_{-1,2}\Box+B$, with $B\in\ePs{0,0}(M),$
$W_{-1,2} \in \ePs{-1,2}(M).$ By the hypotheses, $\Box u=0$ while for any
$Q_\alpha\in \ePs{0,0}(M)$ with $\eWF'(Q_\alpha)\subset O$, $Q_\alpha
A_{\beta,m,l}Bu\in L^2(M)$ since $u$ is coisotropic of order $k-1$ on $O$
relative to $\eH{m,l}(M)$ and $Q_\alpha A_{\beta,m,l}B\in\coiso^{k-1}$.
\end{proof}

Thus, below $\alpha,\beta$ will stand for {\em reduced multiindices},
with $\alpha_N=0$, $\beta_N=0$.

\begin{lemma}\label{HMV.104}(Special case of \cite[Lemma~6.5]{MR2020655} adopted to the present setting.)
Suppose $Q \in \ePs{0,0}(M)$, and let $C_0$ and $C_{ij}$ be given by
\eqref{HMV.103}, resp. \eqref{eq:A_i-comm}. Then, 
assuming \eqref{eq:A_i-comm}, and using the notation of
\eqref{eq:red-multi}, we have 
\begin{equation}\begin{split}
&\sum_{|\alpha|=k}i [A_{\alpha,m,l}^*Q^*QA_{\alpha,m,l},\Box]\\
&\quad=\sum_{|\alpha|,|\beta|=k}A_{\alpha,m+1/2,l-1}^*Q^*C'_{\alpha\beta}
QA_{\beta,m+1/2,l-1}\\
&\quad\quad+\sum_{|\alpha|=k}
\left(A_{\alpha,m+1/2,l-1}^*Q^*E_{\alpha,m+1/2,l-1}
+E^*_{\alpha,m+1/2,l-1}Q A_{\alpha,m+1/2,l-1}\right)\\
&\quad\quad+\sum_{|\alpha|=k} A_{\alpha,m,l}^* i[Q^*Q,\Box]A_{\alpha,m,l}, 
\label{HMV.105}\end{split}\end{equation}
where
\begin{equation*}
E_{\alpha,m,l}=W_{m,l}E_{\alpha},\ E_{\alpha}\in \coiso^{k-1}+\coiso^{k-1}A_N,
\ \eWF'(E_{\alpha})\subset\eWF'(Q),
\label{HMV.106}
\end{equation*}
and for all $\alpha$, $\beta$,
\begin{equation}
C'_{\alpha\beta} \in \ePs{0,0}(M), \ 
\sigma(C'_{\alpha\beta})|_\fcaldot=-(m+l)\xis\delta_{\alpha\beta},
\label{eq:C'>0}\end{equation}
with $\delta$ denoting the Kronecker delta function.
\end{lemma}

\begin{remark}
The first term on the right hand side of \eqref{HMV.105} is the principal term
in terms of $\coiso$ order; both $A_{\alpha,m,l}$ and $A_{\beta,m,l}$ have
$\coiso$ order $k$. Moreover,
\eqref{eq:C'>0} states that it has non-negative principal symbol near
$\pa\fcaldot$. The terms involving $E_{\alpha,m,l}$ have $\coiso$ order $k-1$,
or include a factor of $A_N$, so they can be treated as error terms.
On the other hand, one does need to arrange that $i[Q^*Q,\Box]$ is positive,
as discussed below.
\end{remark}

\begin{proof}
The commutator with $\Box$ distributes over the factors in
$A_{\alpha,m,l}=W_{m,l}A_1^{\alpha_1}\ldots A_{N-1}^{\alpha_{N-1}}$. Using
\eqref{eq:A_i-comm} for each individual commutator $i[A_i,\Box]$,
\eqref{HMV.103} for $W_{m,l}$, and rearranging the factors (with error,
i.e.\ commutator, terms arising from rearrangement included in one of
the $E_{\beta,m,l}$ terms as $\module$ is a Lie algebra) gives the
conclusion. See \cite[Proof of Lemma~6.5]{MR2020655} for a more
leisurely discussion.
\end{proof}

We now consider the operators from Lemma~\ref{HMV.104}
$C'=(C'_{\alpha\beta}),$ $|\alpha|=k=|\beta|,$ as a matrix of operators, or
rather as an operator on a trivial vector bundle with fiber $\RR^{|M_k|}$
over a neighborhood of $\fcaldot$, where $|M_k|$ denotes the number of
elements of the set $M_k$ of multiindices $\alpha$ with $|\alpha|=k.$ Let
$c'=\sigma(C')|_{\pa\fcaldot}.$

Then $c'=\sigma(C')|_{\pa\fcaldot}$ is
positive or negative definite with the sign of $\sigma(C_0)|_{\pa\fcaldot}$.
The same is therefore true microlocally near $\pa\fcaldot.$ For the sake of
definiteness, suppose that $\sigma(C_0)_{\pa\fcaldot}>0.$ Then there exist a
neighborhood $O_k$ of $\pa\fcaldot,$ depending on $|\alpha|=k,$ and
$B\in\ePs{0,0}(M),$ $G\in\ePs{-1,0}(M),$ with
$\sigma(B)>0$ on $O_k$ such that
\begin{equation}\label{eq:Q-O_m}
Q\in\ePs{m,l}(M),\ \eWF(Q)\subset O_k\Rightarrow
Q^*C'Q=Q^*(B^*B + G)Q.
\end{equation}

With $O_k$ as above,
we assume that $\eWF'(Q)\subset O_k$,
\begin{equation}\begin{split}
&i[Q^*Q,\Box]= \sum\tilde B_j^*\tilde B_j + \tilde G
+\tilde F, \text{ where }\\
&\tilde B_j\in\ePs{1/2,-1}(O),\ \tilde F\in\ePs{1,-2}(O),
\ \tilde G \in \ePs{0,-2}(O).
\end{split}
\label{eq:[Q,P]-0}\end{equation}
In the actual application to solutions $u$ of the wave equation, stated below,
$\tilde F$ will be such that $u$ is {\em a priori}
regular on $\eWF'(\tilde F)$, namely it is coisotropic of the order
that we wish to propagate.

In fact, due to the two step nature of the proof below, we also need
another microlocalizer $Q'\in\ePs{0,0}(M)$ satisfying analogous
assumptions with $\tilde B_j$, etc., replaced by $\tilde B_j'$, etc.,
\begin{equation}
i[(Q')^*Q',\Box]= \sum(\tilde B_j')^*\tilde B_j' + \tilde G'
+\tilde F',
\label{eq:[Q,P]-p}\end{equation}
with properties analogous to \eqref{eq:[Q,P]-0}, except that
$\eWF'(Q')\subset O'_k$, etc., where $O'_k$ is contained in the elliptic
set of $Q$.

\begin{proposition}\label{prop:HMV1-mod}
(Cf.\ \cite[Appendix]{Hassell-Melrose-Vasy2};
modified version of \cite[Proposition~6.7]{MR2020655}.)
Suppose that $\xis<0$ on $O_k$, $u$ is coisotropic of order $k-1$ on $O_k$
relative to $\eH{m,l}(M)$, $m>l+\frac{f}{2}$,
$\eWF{}(\Box u)\cap O_k=\emptyset$ and that there exist
$Q,Q'\in\ePs{0,0}(O_k)$
that satisfies \eqref{eq:[Q,P]-0}-\eqref{eq:[Q,P]-p} with
$u$ coisotropic of order $k$ on a neighborhood of $\eWF'(\tilde F)\cup
\eWF'(\tilde F')$ relative to $\eH{m,l}(M)$.
Then $u$ is coisotropic of order $k$ on $O''$ relative to $\eH{m,l}(M)$
where $O''$ is the elliptic set of $Q'.$

The same conclusion holds if $\xis>0$ on $O_k$ and $m<l+\frac{f}{2}$.
\end{proposition}

\begin{proof}
For the reader's convenience, we sketch the argument. We first prove
a weaker statement, namely that  $u$ is coisotropic of order $k$ on
$O'$ relative to $\eH{m-1/2,l}(M)$
where $O'$ is the elliptic set of $Q.$ Then we strengthen the result
to complete the proof of the proposition.

First consider $u$ coisotropic of order $k$ on $O_k$ relative to $\eH{m,l}(M)$.
Set
\begin{equation*}
r=m-1\Mand s=-l-\frac{f-1}{2}.
\end{equation*}
Let $Au'=(QA_{\alpha,r+1/2,s-1} u')_{|\alpha|=k},$
regarded as a column vector of length $|M_k|.$ Now consider
\begin{multline}
\sum_{|\alpha|=k}
\langle u',i[A_{\alpha,r,s}^*Q^*QA_{\alpha,r,s},\Box]u'\rangle
=\|BAu'\|^2 + \langle Au', GAu' \rangle +\\
\sum_{|\alpha|=k}
\Big(\langle QA_{\alpha,r+1/2,s-1} u',E_{\alpha,r+1/2,s-1} u'\rangle
+\langle E_{\alpha,r+1/2,s-1} u',QA_{\alpha,r+1/2,s-1} u'\rangle \Big)\\
+\sum_{|\alpha|=k} \Big(\|\tilde B A_{\alpha,r,s}u'\|^2
+\langle A_{\alpha,r,s}u', \tilde F A_{\alpha,r,s}u'\rangle
+\langle A_{\alpha,r,s}u', \tilde G A_{\alpha,r,s}u'\rangle \Big).
\label{HMV.108e}\end{multline}
Dropping the term involving $\tilde B$ and applying the Cauchy-Schwarz
inequality to the terms with $E_{\alpha,r+1/2,s-1}$, we have for any
$\ep>0$,
\begin{equation*}\begin{split}
\|BAu'\|^2\leq &\sum_\alpha
\Big|\langle u',i[A_{\alpha,r,s}^*Q^*QA_{\alpha,r,s},\Box]u'\rangle 
\Big|+\ep\| Au' \|^2\\
&\qquad+\ep^{-1} \sum_\alpha \|E_{\alpha,r+1/2,s-1} u'\|^2
+|\langle Au', GAu' \rangle|\\
&\qquad+\sum_{|\alpha|=k} \Big(|\langle A_{\alpha,r,s}u', \tilde F
A_{\alpha,r,s}u'\rangle|
+|\langle A_{\alpha,r,s}u', \tilde G A_{\alpha,r,s}u'\rangle| \Big).
\end{split}\end{equation*}
Choosing $\ep>0$ small enough, the second term on the right can
be absorbed in the left hand side (since $B$ is strictly positive), 
and we get 
\begin{equation}\begin{split}
\frac1{2} \|BAu'\|^2\leq &\sum_\alpha
|\langle u',i[A_{\alpha,r,s}^*Q^*QA_{\alpha,r,s},\Box]u'\rangle|\\
&\qquad
+\ep^{-1} \sum_\alpha \|E_{\alpha,r+1/2,s-1} u'\|^2
+|\langle Au', GAu' \rangle|\\
&\qquad+\sum_{|\alpha|=k}
|\langle A_{\alpha,r,s}u', \tilde G A_{\alpha,r,s}u'\rangle|
+\sum_{|\alpha|=k}|\langle A_{\alpha,r,s}u', \tilde F
A_{\alpha,r,s}u'\rangle|
\end{split}\label{eq:B'}\end{equation}
Now, all but the first and last terms on the right hand side are bounded by
the square of a
coisotropic order $k-1$ norm of $u'$ relative to $\eH{m,l}(M)$.

For the $E_{\alpha,r+1/2,s-1}$ term this is immediate,
for $r+1/2=m-1/2\leq m$, $s-1=-l-\frac{f+1}{2}$, and $E_{\alpha,r+1/2,s-1}=
W_{r+1/2,s-1}E_\alpha$, $E_\alpha\in\coiso^{k-1}$, and as the pairing
is relative to the Riemannian density $|dg|$, while the
weighting of $\eH{m,l}(M)$ is relative to a smooth non-degenerate b-density,
$\|E_{\alpha,r+1/2,s-1} u'\|$ being bounded by a
coisotropic order $k-1$ norm of $u$ relative to $\eH{m,l}(M)$
requires that $m-(r+1/2)\geq 0$, $l+\frac{f+1}{2}+s\geq 0$
(and these are $1/2$, resp.\ $0$).

For the $G$ and $\tilde G$ terms, this can be seen by factoring out some
$A_j$ from $A_{\alpha,r+1/2,s-1}=W_{r+1/2,s-1}A_{\alpha}$, i.e.\ writing
it as $A_{\alpha,r+1/2,s-1}=A_{\beta,r+1/2,s-1}A_j$ with $|\beta|=|\alpha|-1
=k-1$, and regarding $A_j$ simply as an operator in $\ePs{1,0}(M)$.
As $G\in\ePs{-1,0}(M)$, the claimed boundedness requires
$2m-(-1+2+2(r+1/2))\geq 0$, $2l+(f+1)+2(s-1)\geq 0$, and both are
in fact $0$. For $\tilde G\in\ePs{0,-2}(M)$, the claimed boundedness
requires (now we have factors of $A_{\alpha,r,s}$)
$2m-(0+2+2r)\geq 0$ and $2l+(f+1)+(-2+2s)\geq 0$, and again both of these
are $0$.

The first term on the right hand side vanishes if $\Box u'=0$, while
the last term is bounded by the square of a
coisotropic order $k$ norm of $u'$ relative to $\eH{m,l}(M)$ in a neighborhood
of $\eWF'(\tilde F)$.

We apply this with $u'$ replaced by $u_\delta=(1+\delta|D_t|)^{-1}u$
where now $u$ is coisotropic of
order $k-1$ relative to $\eH{m,l}$ and still solves $\Box u=0$.
Then letting $\delta\to 0,$ using the strong convergence of
$(1+\delta|D_t|)^{-1}$ to the
identity, and the
assumption that $u$ is coisotropic of order $k$ relative
to $\eH{m,l}(M)$ in a neighborhood of $\eWF'(\tilde F)$, shows that
$BAu\in L^2(M),$ hence that $u$ is coisotropic of order $k$ relative
to $\eH{m-1/2,l}(M)$ on the elliptic set of $B$,
finishing the proof of the first part.

For the second part, we carry through the same argument, but replacing
$Q$ by $Q',$ $\tilde B$ by $\tilde B',$ etc., and taking
\begin{equation*}
r=m-\frac{1}{2},\ s=-l-\frac{f-1}{2}.
\end{equation*}
The claim is that all but the first and last terms on the right hand side
of \eqref{eq:B'} (with the changes just mentioned) are bounded by
the square of either a coisotropic norm of order $k-1$ of $u'$ relative to
$\eH{m,l}(M)$ (namely, the $E_{\alpha,r+1/2,s-1}$ term) or a
coisotropic norm of order $k$ of $u'$ relative to $\eH{m-1/2,l}(M)$
(the $G$ and $\tilde G'$ terms).

For the $E_{\alpha,r+1/2,s-1}$ term, as above, the claimed boundedness
requires $m-(r+1/2)\geq 0$, $l+\frac{f+1}{2}+s\geq 0$, and now
these both vanish.

For the $G$ and $\tilde G'$ terms the conclusion is also immediate,
without factoring
out an $A_j$ this time. Indeed,
as $G\in\ePs{-1,0}(M)$, the claimed boundedness requires
$2(m-1/2)-(-1+2(r+1/2))\geq 0$, $2l+(f+1)+2(s-1)\geq 0$, and both are
in fact $0$. For $\tilde G\in\ePs{0,-2}(M)$, the claimed boundedness
requires (now we have factors of $A_{\alpha,r,s}$)
$2(m-1/2)-(0+2r)\geq 0$ and $2l+(f+1)+(-2+2s)\geq 0$,
and again both of these are $0$.

The regularization argument as above proves that $u$
is coisotropic of order $k$ relative
to $\eH{m,l}(M)$ on the elliptic set of $B$,
finishing the proof of the proposition.
\end{proof}

\def\cprime{$'$} \def\cprime{$'$} \def\cdprime{$''$} \def\cprime{$'$}
  \def\cprime{$'$} \def\ocirc#1{\ifmmode\setbox0=\hbox{$#1$}\dimen0=\ht0
  \advance\dimen0 by1pt\rlap{\hbox to\wd0{\hss\raise\dimen0
  \hbox{\hskip.2em$\scriptscriptstyle\circ$}\hss}}#1\else {\accent"17 #1}\fi}
  \def\cprime{$'$} \def\ocirc#1{\ifmmode\setbox0=\hbox{$#1$}\dimen0=\ht0
  \advance\dimen0 by1pt\rlap{\hbox to\wd0{\hss\raise\dimen0
  \hbox{\hskip.2em$\scriptscriptstyle\circ$}\hss}}#1\else {\accent"17 #1}\fi}
  \def\cprime{$'$} \def\bud{$''$} \def\cprime{$'$} \def\cprime{$'$}
  \def\cprime{$'$} \def\cprime{$'$} \def\cprime{$'$} \def\cprime{$'$}
  \def\cprime{$'$} \def\cprime{$'$} \def\cprime{$'$} \def\cprime{$'$}
  \def\polhk#1{\setbox0=\hbox{#1}{\ooalign{\hidewidth
  \lower1.5ex\hbox{`}\hidewidth\crcr\unhbox0}}} \def\cprime{$'$}
  \def\cprime{$'$} \def\cprime{$'$} \def\cprime{$'$} \def\cprime{$'$}
  \def\cprime{$'$} \def\cprime{$'$} \def\cprime{$'$} \def\cprime{$'$}
  \def\cprime{$'$} \def\cprime{$'$} \def\cprime{$'$} \def\cprime{$'$}
  \def\cprime{$'$} \def\cprime{$'$} \def\cprime{$'$} \def\cprime{$'$}
  \def\cprime{$'$} \def\cprime{$'$} \def\cprime{$'$} \def\cprime{$'$}
  \def\cprime{$'$} \def\cprime{$'$} \def\cprime{$'$} \def\cprime{$'$}
  \def\cprime{$'$} \def\cprime{$'$}
\providecommand{\bysame}{\leavevmode\hbox to3em{\hrulefill}\thinspace}
\providecommand{\MR}{\relax\ifhmode\unskip\space\fi MR }
\providecommand{\MRhref}[2]{%
  \href{http://www.ams.org/mathscinet-getitem?mr=#1}{#2}
}
\providecommand{\href}[2]{#2}


\end{document}